%% file: main.tex
\documentclass[final]{siamltex}
\usepackage{latexsym,enumerate}
\usepackage{graphicx,color}
\usepackage{amssymb,amsmath}
\usepackage{enumerate}

\usepackage{vuk,a4wide,url,xspace,bbm}
\usepackage[american]{babel}
\usepackage{graphicx, psfrag}
\usepackage{subfigure}
\usepackage{wrapfig}
\usepackage{hyperref}
\usepackage[notref,notcite]{showkeys}

\title{Blood-flow modelling along and trough a braided multi-layer metallic stent
\thanks{
This research was partially funded by Cardiatis ({\tt www.cardiatis.com}), 
an industrial partner designing and commercializing metallic wired stents. This work
was supported by a grant from Institut des Syst{\`e}mes Complexes (IXXI, {\tt www.ixxi.fr})
}
}

\author{Vuk Mili{{\v s}}i{{\'c}}\thanks{Wolfgang Pauli Institute (WPI), UMI CNRS 2841,Vienna, AUSTRIA, ({\tt vuk.milisic@imag.fr})}
}

\begin{document}
\renewcommand{\labelenumi}{(\roman{enumi})}

\headheight -1.0cm

\maketitle

\begin{abstract}
In this work we study the hemodynamics in a stented artery 
connected  either to a  collateral artery or to an aneurysmal sac.
The blood flow is driven by the pressure drop.
Our aim is to characterize the flow-rate and the pressure
in the contiguous zone to the main artery: using boundary layer theory 
we construct a homogenized first order approximation with respect to  $\e$, the
size of the stent's wires. This provides an explicit expression 
of the velocity profile through and along the stent. The profile
depends only on the input/output pressure data of the problem
and some homogenized constant quantities: it is explicit.
In the collateral artery this gives
the flow-rate. In the case
of the aneurysm, it shows that : (i) the 
zero order pressure  inside the sac is equal to  the averaged 
pressure along the stent in the main artery, (ii) the presence of the stent 
inverses the rotation of the vortex.
Extending the tools set up in \cite{BrBoMi,MiVws} we prove
rigorously that our asymptotic approximation of velocities and pressures is first order accurate
with respect to  $\e$. We derive then new implicit interface conditions that our approximation formally satisfies, generalizing our analysis to other possible
geometrical configurations.
 In the last part
we provide numerical results that illustrate and validate
the theoretical approach.
\end{abstract}
\begin{AMS}
{76D05, 35B27, 76Mxx, 65Mxx}
\end{AMS}
%
\begin{keywords} 
{wall-laws, porous media, rough boundary, Stokes equation, 
multi-scale modeling, boundary layers, pressure driven flow,
error estimates, vertical boundary  correctors, blood flow, stent, artery, aneurysm}
\end{keywords}

%

\section{Introduction}

\begin{figure}[ht!] 
\begin{center}
\input{stents_pbm}
\includegraphics[width=0.25\textwidth,angle=67]{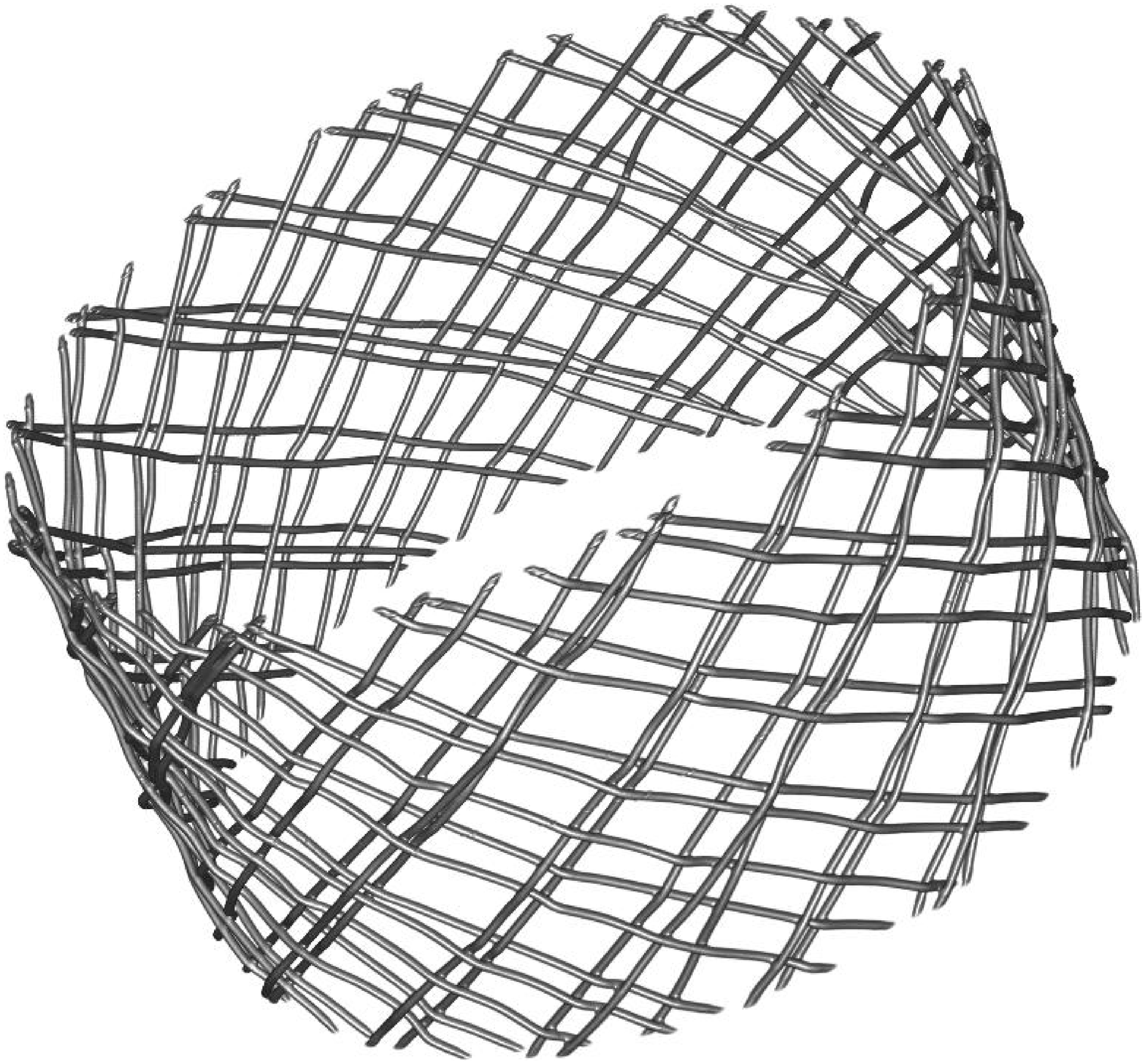}
\caption{A sketch of stented arteries: with a collateral artery (left),  an aneurysmal sac (middle) and a 3D example of a real metallic multi-wired stent (right)}\label{stent}
\end{center}
\end{figure}

Atherosclerosis and rupture of aneurysm are lethal pathologies of the
cardio-vascular system. A possible therapy consists in introducing a metallic
multi-layered stent (see fig. \ref{stent} right). This device
slows down the vortices in the aneurysm and doing so  favors coagulation
of the blood inside the sac. This, in turn, avoids possible rupture of the sac.

In this study we aim to investigate the fluid-dynamics
of blood in the presence of a stent. We  focus on two 
precise configurations in this context:
(i) a stented artery is connected to the collateral artery but the
aperture of the latter is partially occluded by the presence of the
stent (see fig. \ref{stent} left),
(ii) a sacular aneurysm is present behind a stented artery (fig. \ref{stent} middle).
From the applicative point of view these two situations are of interest since they 
represent a dual constraint that a stent should optimize somehow: the grid generated by the wires should be coarse enough to provide blood to the collateral arteries
(for instance iliac arteries in  the aorta), at the same time the wires should 
be close enough to have a real effect in terms of velocity reduction in the aneurysm.

Multi-layer metallic wired stents seem to satisfy both the constraints at 
the same time. Although experimentally exhibited \cite{Aus.08,LiLiChou},  these facts needed a better mathematical understanding.
We give here results in this sense, setting a common framework for both phenomena in the case of the Stokes flow. 

Inspired by homogenization  techniques applied to the case of rough
boundaries \cite{AchPiVaJCP.98,JaMi.03, NeNeMi.06} we construct a first-order 
multi-scale
approximation of the velocity and the pressure.
By averaging, we get a first order accurate macroscopic 
description of the
fluid flow. Indeed, we compute an explicit expression of the velocity
through the fictitious interface supporting the stent and separating the
main artery from the contiguous zone. This 
formula only depends on the input data of the problem
and some homogenized constants obtained solving microscopic  cell problems.
In the case of the aneurysmal sac we show rigorously
that the zero order pressure in the sac is constant and averaged with respect to
the pressure in the main artery, which was not known.
Then we show that formally this leads also to redefine the
problem in a new and implicit way in the domain decomposition flavor. 
Actually
we obtain a new set of interface conditions along the fictitious interface:
while for the normal velocity they look similar to those presented in \cite{CioMu.97,AllairePassoire,BrillardPassoire}, 
the tangential conditions are new to our knowledge.
They express a slip velocity in the main artery (as in \cite{JaMiSIAM.00}), 
but a discontinuous homotetic relationship between horizontal velocities across the interface of the stent (see system \eqref{homo_fo}). Our results concern the steady Stokes equations, as in \cite{AllairePassoire}, the same interface conditions are valid in the case unsteady Navier-Stokes case.

From the mathematical point of view this paper introduces several
novelties. The case of a sieve has been widely studied in a 
 different setting in \cite{conca1,conca2,AllairePassoire,BrillardPassoire,BouGiSiam.01}. 
In these works, the authors considered no-slip obstacles set on a surface 
with various dimensionalities but with a common point: 
the velocity was completely 
imposed at the inlet/outlet boundaries of the fluid domain. 
Although this could seem a technicality,  
it influences drastically the limiting regime of the flow. 
Indeed a complete velocity profile is imposed as a Dirichlet condition
at the inlet/outlet of the domain, so that
the total flow-rate through the sieve remains  constant whatever $\e$, the size  
of the obstacles:
a resistive term appears as a zeroth order limit in the fluid equations.
In the context of blood flow such  a regime seems hard to reach:
experiences show that when the wires are to dense
no transverse flow crosses the stent. This suggests
that through the porous interface, 
blood flow should be driven by a pressure drop 
more that a fixed flow-rate. 

In this direction, J{\"a}ger and Mikeli{\'c}
considered a pressure driven fluid in \cite{JaMiFilter}. But 
they studied an interface whose thickness was independent on $\e$, 
 which  seemed useless for our purpose : the diameter of the wires of 
the stent are dependent on the radius of the artery where the
stent should be implanted. It appears natural to consider 
roughness size that varies wrt $\e$ in any direction.
Moreover in this paper we introduce 
both a tangential and a transverse flow along and trough the stent.
Indeed, in  the limiting regime considered by J{\"a}ger and Mikeli{\'c} \cite{JaMiFilter},
the velocity  is zero. Here when the collateral artery or a sac
are completely closed by the stent, we still expect a Poiseuille profile in 
the main artery. 

At a more technical level, this work improves  the approach developed
in \cite{BrBoMi,MiVws} in order to  correct
edge oscillations introduced by periodic boundary layers. At the same time,
we give  an appropriate framework to 
deal with this problem in the case of Stokes equations. 
Indeed, due to the presence of the obstacles, the divergence operator is singular wrt $\e$, 
this implies degradation of convergence results when lifting the non-free divergence terms and estimating the pressure.
In this frame, 
we decompose the corrections of the superfluous boundary layer oscillations
in two parts~:
\begin{itemize}
\item on the microscopic side we use weighted Sobolev spaces to describe 
the behaviour at infinity of the vertical corner correctors, defined on a half plane. 
This provides accurate decay rates with respect to 
$\e$ at the macroscopic level near the corner. Indeed, using onto mappings 
between weighted Sobolev spaces we improve decay estimates
already derived in the scalar case in \cite{BrBoMi,MiVws}. Then in the spirit of \cite{AllairePassoire} we 
construct a microscopic lifting operator that allows the vertical correctors to fullfil the Dirichlet condition on the obstacles,
\item a complementary 
macroscopic corrector is added in a second step, that handles exponentially 
decreasing  errors far from the corners. 
\end{itemize}
An attempt to break the periodicity at the inlet/outlet of the domain 
was done in \cite{JaMiSIAM.00} by using a vertical corrector localized 
in a tiny strip near the vertical interface. 
But, decay estimates claimed in formula (77) p. 1123 \cite{JaMiSIAM.00} 
seem to work, to our knowledge, only 
for {\em a priori} estimates of the error and are not accurate
enough to be used in the very weak estimates.


We underline as well that in the  literature \cite{JaMiSIAM.00,JaMiJDE.01, JaMi.03,JaMiNe.01}  error estimates
between the direct rough solution and the approximations constructed thanks to boundary layer arguments concerned the $\bL^2$ norm of the velocity. 
In this paper we provide error estimates  of the same order for the pressure as well in the negative Sobolev $H^{-1}$ norm. 
This is obtained using the  microscopic nature of the  pressure correctors and in particular thanks to the very precise control of lateral 
correctors. We stress that these vertical correctors play a crucial part in our error analysis at several steps of this work.

The paper is organized as follows: in the two next sections, 
after some basic notations and definitions,  we give 
a detailed review of the results obtained
either in the case of a collateral artery or a sacular
aneurysm. 
We give in section \ref{sect.techn.prelim} the abstract results that are used in
section \ref{proofs} in order to prove
the claims. We provide 
numerical results showing  a first order
accuracy also in the discrete case in section \ref{num}. 
In Appendix \ref{app.1},
we give  proofs of existence, uniqueness and {\em a priori} estimates
for vertical correctors in  the weighted Sobolev spaces, while in 
Appendix \ref{periodic.bl} we detail the results claimed 
for the periodic boundary layers throughout the paper.

\section{Geometry and problem settings}

\subsection{Geometry}\label{geo}
\begin{figure}[ht]
\hspace{-0.5cm}
\input{domains}
\caption{The macroscopic domains $\Oe$ (left) and  $\Omega$
(middle) and the microscopic infinite strip $Z$ (right)}\label{domaines}
\end{figure}
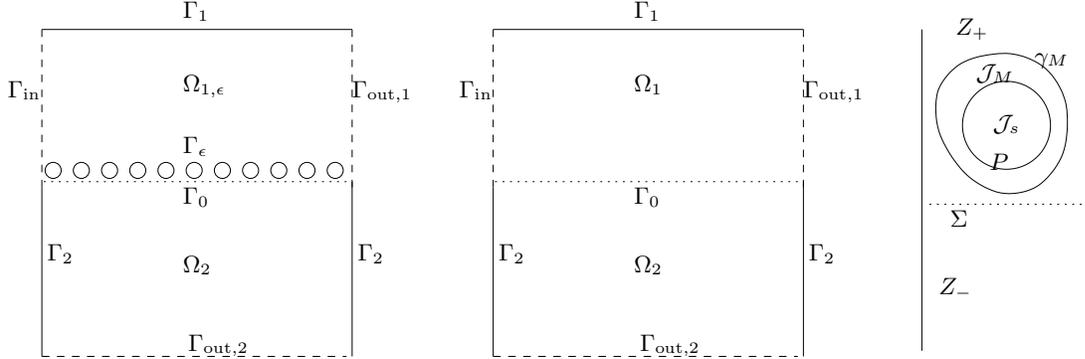

In this study we consider two space dimensions.
Let us define by $\js$ one or more solid obstacles included in $\cJ:=]0,1[^2$ of Lipschitz boundaries denoted $P$ in the sense of the definition  p. 13-14 of Chap. 1 in \cite{Ne.Book.67}.
We denote by $\jf:=]0,1[^2 \setminus \js$ the complementary fluid part of $\js$ in $]0,1[^2$. Also, we consider a smooth surface $\gamma_M$ strictly contained in $\js$ and enclosing 
$P$ and we denote $\jm$ the domain contained between $\gamma_M$ and $P$.
Then we define:
\begin{enumerate}[(i)]
\item Macroscopic domains:

The $\e$-periodic repetition of $\jf$ is denoted by $\leps$ and reads:
$$
\leps:= \cup_{i=0}^{m} \e ((i,0)+\jf), \text{ where } m:= \ue,
$$
the real $\e$ is always chosen such that $m$ is an integer. Then we set:
$$
\begin{aligned}
&\Ou := ] 0,1 [^2 ,  & \Gin := \{ 0 \} \times ]0,1[,\\
&\Oup  := ]0,1[ \times ]\e ,1[, & \Goutu := \{ 1 \} \times ]0,1[,\\
&\Oeu  := \Oup \cup (]0,1[ \times \{\e \} ) \cup \leps, & \Goutd := ]0,1[ \times \{ -1 \}, \\
&\Ode  := ]0,1[\times ]-1,0[, & \Gun := ]0,1[ \times \{ 1 \}, \\
&\Gz  := ]0,1[ \times \{ 0 \}, & \Gde := \{ 0 \} \times ]-1,0[ \cup \{ 1 \}\times  ]-1,0[ , \\
&\Omega  := \Ou \cup \Gz \cup \Ode, & \Gd := \Gun \cup \Gde \cup \Geps,\\
& \Oe := \Oeu \cup \Gz \cup \Ode, & \Gn := \Gin \cup \Goutu \cup \Goutd,\\
& \Op := \Oup \cup \Ode,& \rlay:=]0,1[\times]0,\e[.
\end{aligned}
$$
The spatial variable giving the position of a point in domains above is a vector called $x$.
\item The microscopic cell domain:

As the problem contains a solid interface surrounded by a fluid, the microscopic cell problems are set on an infinite strip $Z$ defined as follows
$$
\begin{aligned}
& Z^- := ]0,1[ \times ]-\infty ,0[, \\
& \Sigma := ]0,1[ \times \{ 0 \} , \\
& Z^+ :=]0,1[\times \rr \setminus \js, \\
& Z := Z^+ \cup \Sigma \cup Z^-, \\ 
& Z_{\gamma,\nu}:= Z \cap ]0,1[\times ]\gamma,\nu[ , & (\gamma,\nu) \in \RR^2  \text{ s.t. } \gamma<\nu.
\end{aligned}
$$
The microscopic position variable is denoted by $y:=x/ \e$.
\item The ``corner''  microscopic domain:

In order to handle periodic perturbations on the lateral boundaries $\Gin\cup\Gde \cup \Goutu$ one needs to define a microscopic zoom near the corners $O:=(0,0)$ and $\ov{x}:=(1,0)$ of $\Oe$. This leads to set the half-plane $\Pi$ and the corresponding boundaries as
$$\begin{aligned}
& \Pi := \rr \times \RR \\
& N := \{ 0 \} \times ] 0 , + \infty [, \\
& D := \{ 0 \} \times ]-\infty , 0 [, \\
\end{aligned}
$$
\end{enumerate}
If we choose the obstacle $\js$ to be a single disk, then a graphical 
illustration depicts the definitions above in fig. \ref{domaines} for $\e=1/11$.

The exterior normal vector
to any domain is denoted by $\bfn$, if not stated explicitly 
$\bfn$ is orientated from $\Ou$ towards $\Ode$ on the fictitious interface $\Gz$.
The tangent vector is defined as ${\bft}$. 

\subsection{Notations and definitions}\label{cut.off}
\begin{enumerate}[(i)]
\item Any two-dimensional vector is denoted by a bold symbol:
$\bfu := (u_1, u_2)$,
and single components are scalar and are not bold.
The same holds for the function spaces these vectors belong to:
bold letters denote vector spaces, for instance
$\bfL^2(\Omega) := (L^2(\Omega))^2$.
\item If $\eta \in H^1_{\loc}(Z)$ then we set
$$
\oeta(y_2) := \int_0^1 \eta(y_1,y_2) dy_1 ,\quad y_2 \in \RR,
$$
to be the horizontal average of a function defined on the infinite periodic strip $Z$.
Moreover by the double bar we denote a piecewise constant function defined on $Z$ 
as
$$
\ooeta(y):= \oeta(+\infty) \chiu{Z^+}(y) +  \oeta(-\infty) \chiu{Z^-}(y) , \quad y \in Z,
$$
whenever the function $\oeta(\cdot)$ admits finite limits when $|y_2| \to \infty$.
We need the values of the above function near the origin, thus  we set also:
$$
\ooeta^\pm := \ooeta(0^\pm).
$$
\item For any pair $(\bfu,p) \in \bfL^2(\Omega) \times H^{-1} (\Omega)$ we denote by $\sigma_{\bfu,p}$ the $2\times 2$ distributional matrix reading
$$
\sigma_{\bfu,p}(x) := \nabla \bfu - p \id  ,\text{ a.e. } x \in \Omega ,
$$
where $\id$ is the identity matrix in $\RR^2$. The tensor $\sigma_{\bfu,p}$ looks like the stress tensor but it is not symmetric. 
This is due to the incompressibility constraint: the Stokes problem can still be put in the divergence form with the definition of $\sigma_{\bfu,p}$ above.
\item The brackets $[\cdot]$ denote throughout the whole paper the jump of the quantity enclosed across  fictitious interfaces:
across $\Gz$ on the macroscopic scale, or across $\Sigma$ on the microscopic scale, so that for instance
$$
[ \sigma_{\bfu,p} ] := \sigma_{\bfu,p}(x_1,0^+) - \sigma_{\bfu,p}(x_1,0^-) ,\text{  while  } [\ooeta]:= \ooeta^+- \ooeta^-.
$$
\item For every microscopic function $\eta$ defined on either $Z$ or $\Pi$, we denote by
$$
\eta_\e (x) = \eta\lrxe, \quad \forall x \in \Oe.
$$
\newcounter{enumi_saved}
\setcounter{enumi_saved}{\value{enumi}}
\end{enumerate}
We also need cut-off functions that we define here:
\begin{enumerate}[(i)]
\setcounter{enumi}{\value{enumi_saved}}
\item The cut-off $\phi$ is a scalar function $\phi:\rr \to [0,1]$ s.t.
$\phi$ is a $C^\infty(\rr)$ monotone decreasing function and
$$
\phi(z):= \left\{
\begin{aligned} 
1 \text{ if } z \leq 1,\\
0 \text{ if } z \geq  2, \\
\end{aligned}
\right.
$$
for any positive real $z$.
\item The ``corner'' cut-off functions :
set $\psi_1:= \ov{\psi}(x)$ and $\psi_2:= \ov{\psi}(x-\ov{x})$ 
and $\ov{\psi}$ is a radial monotone decreasing cut-off function such that
$$
\ov{\psi}(x):=
\left\{
\begin{aligned} 
1 \text{ if } |x| \leq \frac{1}{3}, \\
0 \text{ if } |x| \geq  \frac{2}{3}, \\
\end{aligned}
\right.  \quad \forall x \in \RR^2. 
$$
Finally set
$\psi(x):=\psi_1(x)+\psi_2(x)$.
Note that with this definition $\ddn{\ov{\psi}}=0$ on $\giou$. 
\item The ``far from the corner'' cut-off function :
$\Phi$ is defined in a complementary manner on $\gio\cup\Gde$
such that
$$
\left\{
\begin{aligned}
& \psi+\Phi=1,\quad \text{ on } \gio\cup\Gde,  \\
& \ddn{\Phi}=0 \quad \text{ on } \gio  ,
\end{aligned}
\right. 
$$
and one shall take for instance $\Phi(x):=1-\psi(0,x_2)$ for all $x$ in $\Omega$.
\item For regularity purposes we set $\ldd$ to be a cut-off function
in the $\e$-neighborhood of the corners $O$ and $\ov{x}$ (p. 1122 \cite{JaMiSIAM.00}). First we set at the microscopic level:
$$
\lambda(y) := |y_2 | \chiu{B(O,1)}(y) + \frac{y_2}{|y|} \chiu{\{\Pi  \setminus B(O,1)\}}(y) , \quad \forall y  \in \Pi, 
$$
then we define 
$$
\ldd(x):= \lambda\lrxe \psi_1(x) + \lambda\left( \frac{1-x_1}{\e},\frac{x_2}{\e} \right)  \psi_2 ( x)+ \Phi (x), \quad \forall x \in \Omega.
$$
and an easy computation shows that
\begin{equation}\label{eq.lift.h1}
\nrm{\ldd}{H^1(\Omega )} \leq k \{| \log(\e) |^\ud   +1 \} ,
\end{equation}
where the constant $k$ does not depend on $\e$.
\end{enumerate}

\section{Main results}
\subsection{The case of a collateral artery}

We study the problem : find $(\ueps,\peps)$ solving the stationary Stokes equations 
\begin{equation}
\label{exact}
\left\{
\begin{aligned}
& -\Delta \ueps + \nabla \peps = 0 & \text{ in } \Oe,\\ 
& \dive \ueps = 0 & \text{ in } \Oe,\\
& \peps = \pin  \text{ on } \Gin, \quad \peps = \poutu  \text{ on } \Goutu, & \peps= \poutd  \text{ on }\Goutd, \\
&  \ueps \cdot \bft = 0 & \text{ on } \Gin\cup\Goutu\cup\Goutd, \\
&  \ueps = 0 & \text{ on } \gud\cup\Geps, \\
\end{aligned}
\right.
\end{equation}
Because of the microscopic structure of $\Geps$,
the solution of such a system is complex and expensive from the numerical point of view. For this reason throughout this article we use homogenization
 in order to construct approximations of $(\ueps,\peps)$.
This technique decomposes 
in two steps : 
\begin{enumerate}[1)]
\item the derivation of a multi-scale asymptotic expansion and 
the construction of an averaged macroscopic approximation. 
The first part can be seen as an iterative algorithm with respect to   powers 
of $\e$~: 
\begin{enumerate}
\item pass to the limit with respect to  $\e$ and obtain a macroscopic zero order approximation. In our case, because of the straight geometry of the main artery and of the boundary conditions, the Poiseuille profile
is obtained in $\Ou$ and a trivial solution in $\Ode$~:
\begin{equation}\label{poiseuille.col}
\left\{ 
\begin{aligned}
& \bfuz(x) = \frac{\pin-\poutu}{2}(1-x_2) x_2 \eu \chiu{\Ou},\\
& \pz(x) = (\pin (1-x_1) + \poutu x_1 ) \chiu{\Ou} + \poutd \chiu{\Ode}, 
\end{aligned}
\right.
\quad \forall x \in \Omega ,
\end{equation}
\item \label{ii} construct microscopic boundary layers  correcting errors
made by the zeroth order approximation on $\Geps$ and $\Gz$: we set up in the next section
three boundary layers $(\bfbeta,\pi),(\bfups,\varpi)$ and $(\bfchi,\eta)$ 
to this purpose. These functions solve periodic microscopic problems \eqref{beta.cell}, \eqref{eq.upsi} and \eqref{chi.cell}
on the strip $Z$. 
\item compute the constants that 
these correctors reach at $y_2=\pm\infty$: $(\oobeta^\pm,0)$, $(\ooupsi^\pm,0)$ and $(\oochi,\ooeta^\pm)$. Then subtract them to  the correctors.
Physically, $\oobeta^\pm , \ooupsi^\pm $  provide a microscopic 
feed-back relative to the horizontal velocity (see the wall-law
framework  in \cite{NeNeMi.06,JaMiSIAM.00} and references therein) 
whereas the pressure difference $[\ooeta]$ 
represents a microscopic resistivity 
in the flavor of \cite{AllairePassoire,BrillardPassoire,BouGiSiam.01}.
\item take into account the homogenized constants on the limit interface $\Gz$ 
by solving a macroscopic problem: find  $(\bfuu,\pu)$ s.t.
\begin{equation}\label{fo.so.macro}
\left\{
\begin{aligned}
-&\Delta \bfuu + \nabla \pu = 0  &\text{ in } \Omega_1\cup \Omega_2,\\ 
& \dive \bfuu = 0& \text{ in } \Omega_1\cup \Omega_2,\\
&  \bfuu = 0&  \text{ on } \gud,\\ 
& \left.
\begin{aligned}
&  \bfuu \cdot \bft = 0\\ 
&  \pu = 0%
\end{aligned} 
\right\} & \text{ on }  \Gn,\\ 
& \bfuu(x_1,0^\pm) = \left(  \duzu(x_1,0^+) \oobetas^\pm_1 + \lrduzu \ooupsis^\pm_1   \right ) \eu +  \dpde \oochid \ed & \text{ on } \Gz^\pm 
\end{aligned}
\right.
\end{equation}
This macroscopic
corrector depends on the zeroth order approximation and the homogenized
constants. Due to the explicit form of the Poiseuille profile, the Dirichlet data is explicit on both sides of $\Gz$, (nevertheless the solution $(\bfuu,\pu)$ is not explicit inside $\Omega_1\cup\Omega_2$).
\item go to \eqref{ii}  and correct, on a micrscopic scale,  errors made by $(\bfuu,\pu)$ on
$\Geps\cup \Gz$ in order to get higher order terms in the asymptotic ansatz.
\end{enumerate}
\item The second step consists then in  averaging this ansatz and obtaining
an expansion of the macroscopic solutions only. This gives, for instance, at first order~:
$$
\ovueps(x) := \bfuz(x) + \e \bfuu(x), \quad \ovpeps(x) := \pz(x) + \e \pu(x), \quad \forall x \in \Ou \cup \Ode.
$$
In particular as $\ovueps \cdot \bfn = \e [\pz]/[\ooeta]$ on $\Gz$, one gets an explicit first order velocity
profile across $\Gz$. As a consequence, we obtain  a new result~:
\begin{prop}\label{prop.fow}
The flow-rate  in the collateral artery $\Ode$ can be computed explicitly and reads
$$
Q_{\Gz}:= \int_\Gz \ovueps \cdot \bfn dx_1 = \frac{\e}{\left[ \ooeta \right]} \int_\Gz [\pz] dx_1 = \frac{\e}{\left[ \ooeta \right]} \int_\Gz (\poutu + (\pin -\poutu) (1-x_1)-\poutd )\, dx_1 
$$
\end{prop}
As stated above $[\ooeta]$ depends
only on the geometry of the microscopic obstacle $\js$ and is independent of any other parameter.
In the last section of this paper we give some numerical examples that  illustrate  the accuracy of this result (see fig. \ref{fig.profile} and \ref{fig.error}).
Note that the zeroth order approximation does not provide any transverse flow
through $\Gz$. 
Although our results provide  a first order correction, we underline that in the physiological
context the pressures $(\pin,\poutu)$ present in the main artery can be very important compared to $\poutd$~: the first order flow rate $Q_{\Gz}$ can thus
be quantitatively significant as well.
\end{enumerate}

In this work we constructed an suitable mathematical framework in order to analyse
the error made in the two main steps of the construction above. This allows to state the main result of this paper:
%
\begin{thm}\label{main.thm.paper} There exists a unique pair $(\ueps,\peps) \in \bH^1(\Oe)\times L^2(\Oe)$ solving problem \eqref{exact}.
The averaged asymptotic ansatz $(\ovueps,\ovpeps)$ belongs to 
$\bL^2(\Omega_j)\times H^{-1}(\Omega_j)$ for $j\in \{ 1,2 \} $ and satisfies the convergence result
$$
\nrm{\ueps - \ovueps}{L^2(\Omega_1\cup \Ode)} + \nrm{\peps - \ovpeps}{H^{-1}(\Ou'\cup \leps\cup \Ode )} \leq k \e^{\td^-},
$$
where $\td^-$ represent any real number strictly less then $\td$ and
the constant $k$ is independent on $\e$.
\end{thm}

Expressing  interface conditions satisfied by $(\ovueps,\ovpeps)$ on $\Gz$ 
in an implicit way and neglecting higher order rests, we show formally that
in fact $(\ovueps,\ovpeps)$ solve at first order a  new interface problem~: 
\begin{equation}\label{homo_fo}
\left\{
\begin{aligned}
& -\Delta \ovueps+ \nabla \ovpeps= 0 & \text{ in } \Ou\cup\Ode, \\ 
& \dive \ovueps = 0 & \text{ in } \Ou\cup\Ode,\\
& \ovueps  = 0 & \text{ on } \gud,\\ 
& \ovueps \cdot \bft  = 0 & \text{ on } \Gn, \\ 
& \ovpeps= \pin,\text{ on } \Gin, \quad \ovpeps =\poutu \text{ on } \Goutu, \quad \ovpeps =\poutd  & \text{ on }\Goutd, \\ 
& 
\left.
\begin{aligned}
& \ovueps^+\cdot \bft  = \e (\oobetas^+_1 +\ooupsis^+_1 ) \dd{\ovuepsu}{x_2}^+, \quad  
\frac{\ovueps^+\cdot \bft    }{\oobetas^+_1 +\ooupsis^+_1 }= \frac{\ovueps^-\cdot \bft  }{\oobetas^-_1 +\ooupsis^-_1}\\
&\ovueps^+ \cdot\bfn = \ovueps^- \cdot\bfn = - \frac{\e}{[\ooeta]} ([ \sigma_{\ovueps,\ovpeps} ]  \cdot\bfn  ,\bfn)
\end{aligned}
\right\} 
& \text{ on } \Gz.
\end{aligned}
\right.
\end{equation}

The horizontal velocity on $\Gz^+$ is related to the shear rate trough a kind of mixed boundary condition 
alike to the Beaver, Joseph and Saffeman condition \cite{JaMiSIAM.00}. This implicit relationship accounts for the friction effect due to the obstacles that ``resist'' to the flow in the main artery.
Nevertheless because the interface separates two domains $\Ou$ and $\Ode$, we 
obtain a second expression between the upper and the lower horizontal velocities $\ov{u}_{\e,1}^+$ and $\ov{u}_{\e,1}^-$: they are proportional and thus discontinuous. To our knowledge this is  new.

On the other hand, the interface condition on the normal velocity 
could be integrated in the Stokes equations as a kind of ``strange
term'' in the spirit of \cite{CioMu.97,AllairePassoire}, 
but as we are at first order with respect to  $\e$:
(i) the strange term is divided by $\e$ (in \cite{CioMu.97,AllairePassoire} this is a zero order term independent on  $\e$ )
(ii) the derivation does not follow at all the same argumentation.
In a forthcoming work we study the well-posedness
of such a system as well as its consistency 
with respect to  $(\ovueps,\ovpeps)$ and $(\ueps,\peps)$. 
Because of the particular signs of the homogenized constants 
but also the discontinuous nature of the interface 
conditions in the tangential direction to $\Gz$,
this seems a challenging task.

\subsection{The case of an aneurysm}\label{main.res.ane}

The framework introduced above can be extended to the case of an aneurysm; considering the same domain 
$\Oe$ as above we define a new problem : find $(\ueps,\peps)$ solving
\begin{equation}
\label{exact_sac}
\left\{
\begin{aligned}
& -\Delta \ueps + \nabla \peps = 0 & \text{ in } \Oe,\\ 
& \dive \ueps = 0 & \text{ in } \Oe,\\
& \peps = \pin  \text{ on } \Gin, &  \peps = \poutu   \text{ on } \Goutu, \\
&  \ueps \cdot \bft = 0 & \text{ on } \Gin\cup\Goutu, \\
&  \ueps = 0 & \text{ on } \gud\cup\Goutd. \\
\end{aligned}
\right.
\end{equation}
The main difference resides in the boundary condition imposed on $\Goutd$ : here we impose a
complete adherence condition on the velocity ; this closes the output $\Goutd$ and transforms the collateral artery into an idealized square aneurysm.

Again, we construct a similar multi-scale asymptotic ansatz. We extract the macroscopic part to get a homogenized expansion $(\ovueps,\ovpeps)$ reading
$$
\ovueps(x) := \bfuz(x) + \e \bfuu(x), \quad \ovpeps(x) := \pz(x) + \e \pu(x), \quad \forall x \in \Ou \cup \Ode,
$$
where $(\bfuz,\pz)$ is again a Poiseuille profile but complemented by an unknown constant pressure $\pzm$ inside the sac:
\begin{equation}\label{pois.sac}
\left\{ 
\begin{aligned}
& \bfuz(x) = \frac{\pin-\poutu}{2}(1-x_2) x_2 \eu \chiu{\Ou}\\
& \pz(x) = \pzp(x)\chiu{\Ou}  +\pzm \chiu{\Ode}, \\
&  \pzp(x):= \poutu + (\pin -\poutu) (1-x_1), \quad \pzm \in \RR
\end{aligned}
\right.
,\quad \forall x \in \Omega 
\end{equation}
Then again $(\bfuu,\pu)$ solves a mixed 
Stokes problem \eqref{fo.so.macro}, the only difference 
being that $\bfuu=0$ on $\Goutd$.
This gives again a new result:
\begin{coro}\label{coro.pres} 
The zeroth order pressure is constant in $\Ode$, moreover
it  satisfies 
the following compatibility condition with respect to  the pressure in the main artery:
$$
\pzm = \frac{1}{|\Gz|}\int_\Gz \pzp(x_1,0)\, dx_1.
$$
This  gives an explicit velocity profile on $\Gz$ which reads:
$$
\ovueps \cdot \bfn = \frac{\e}{[\ooeta]} (\pzp(x_1,0)-\pzm) + O(\e^2).
$$
\end{coro}

\begin{figure}[ht!] 
\includegraphics[width=0.45\textwidth]{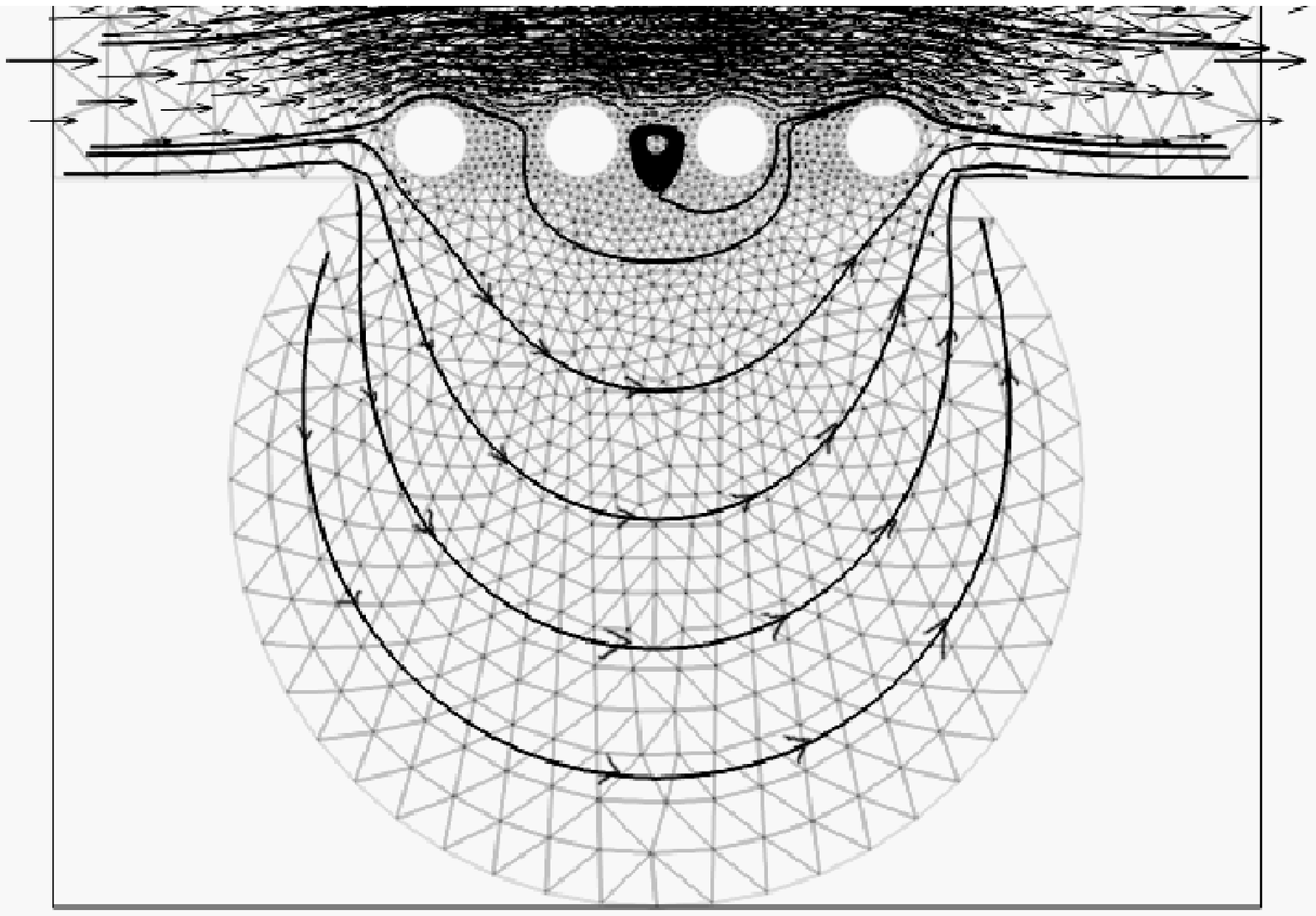}
\includegraphics[width=0.45\textwidth]{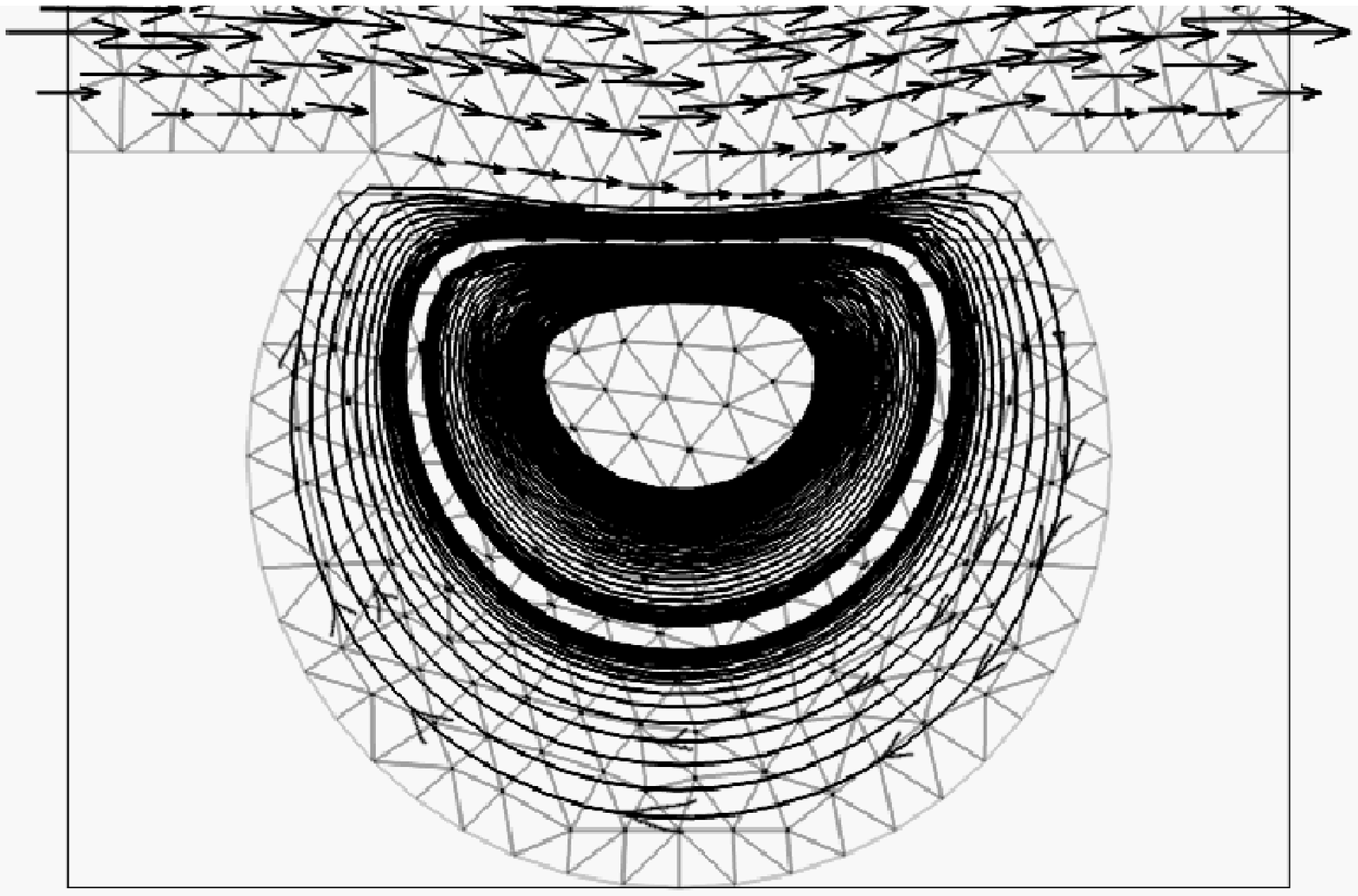}
\caption{Streamlines and velocity vectors in an aneurysmal sac,  with (left) and without a stent (right) }\label{sac_velo}
\end{figure}

The interface  condition exhibited on the normal velocity 
shows rigorously a phenomenon already observed experimentally \cite{Aus.08,LiLiChou}.
Set $x_{1,\max}:= \max_{x \in \Gz} x_1$ (resp. $x_{1,\min}:= \min_{x \in \Gz} x_1$)
and $\ov{x}_1:= (x_{1,\max}+x_{1,\min})/2$, 
when $x_1 < \ov{x}_1$ the pressure jump $[\pz]:=\pzp(x)-\pzm$ is positive,
otherwise it is negative. This implies that the first order 
flow trough the stent is entering $\Ode$ when $x_1<\ov{x}_1$ and
leaving it otherwise. Thus the prosthesis inverses the orientation 
of the cavitation in $\Ode$ with respect to  the non-stented artery (see fig. \ref{sac_velo}).

\begin{figure}[ht!] 
\includegraphics[width=0.45\textwidth]{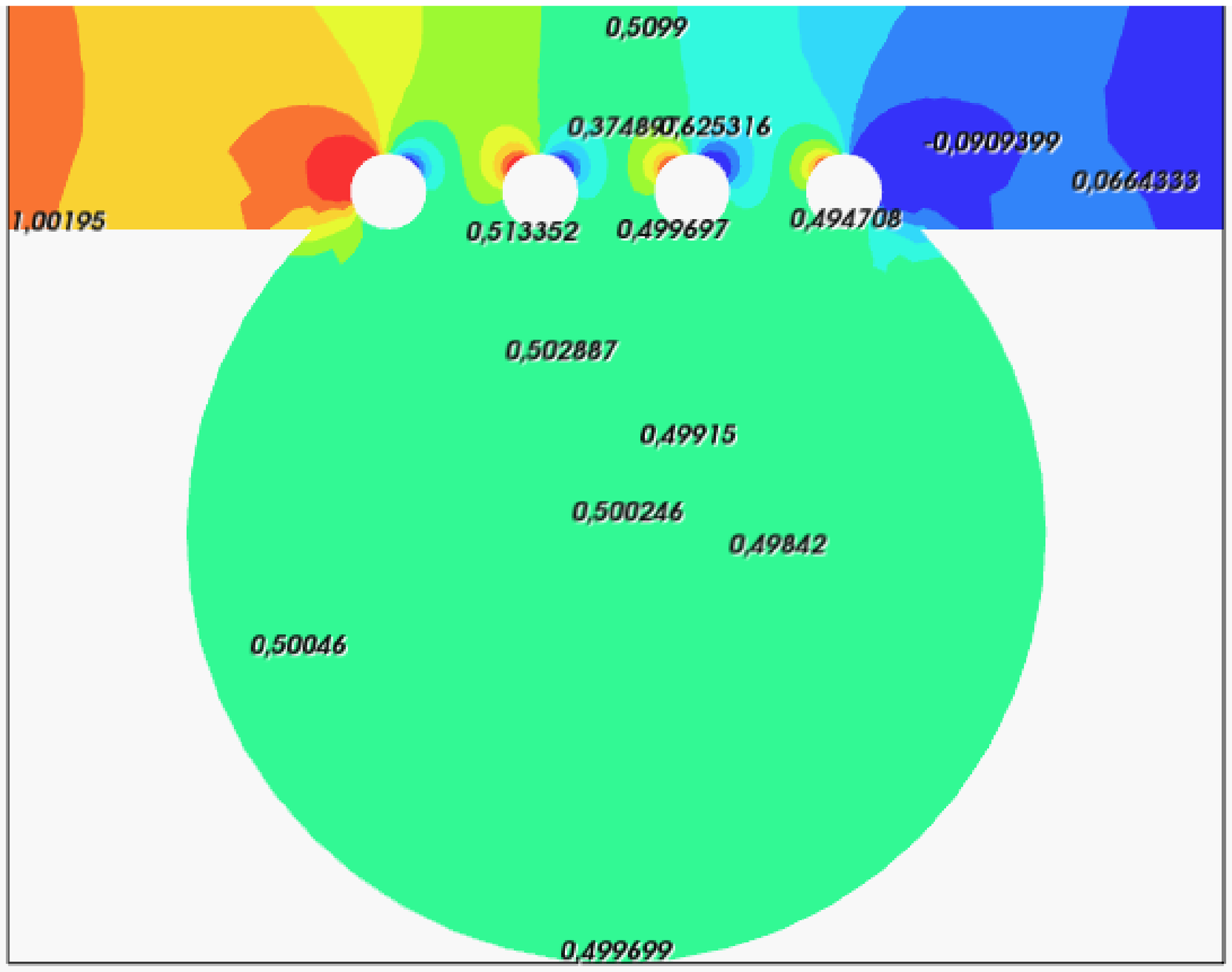}
\includegraphics[width=0.45\textwidth]{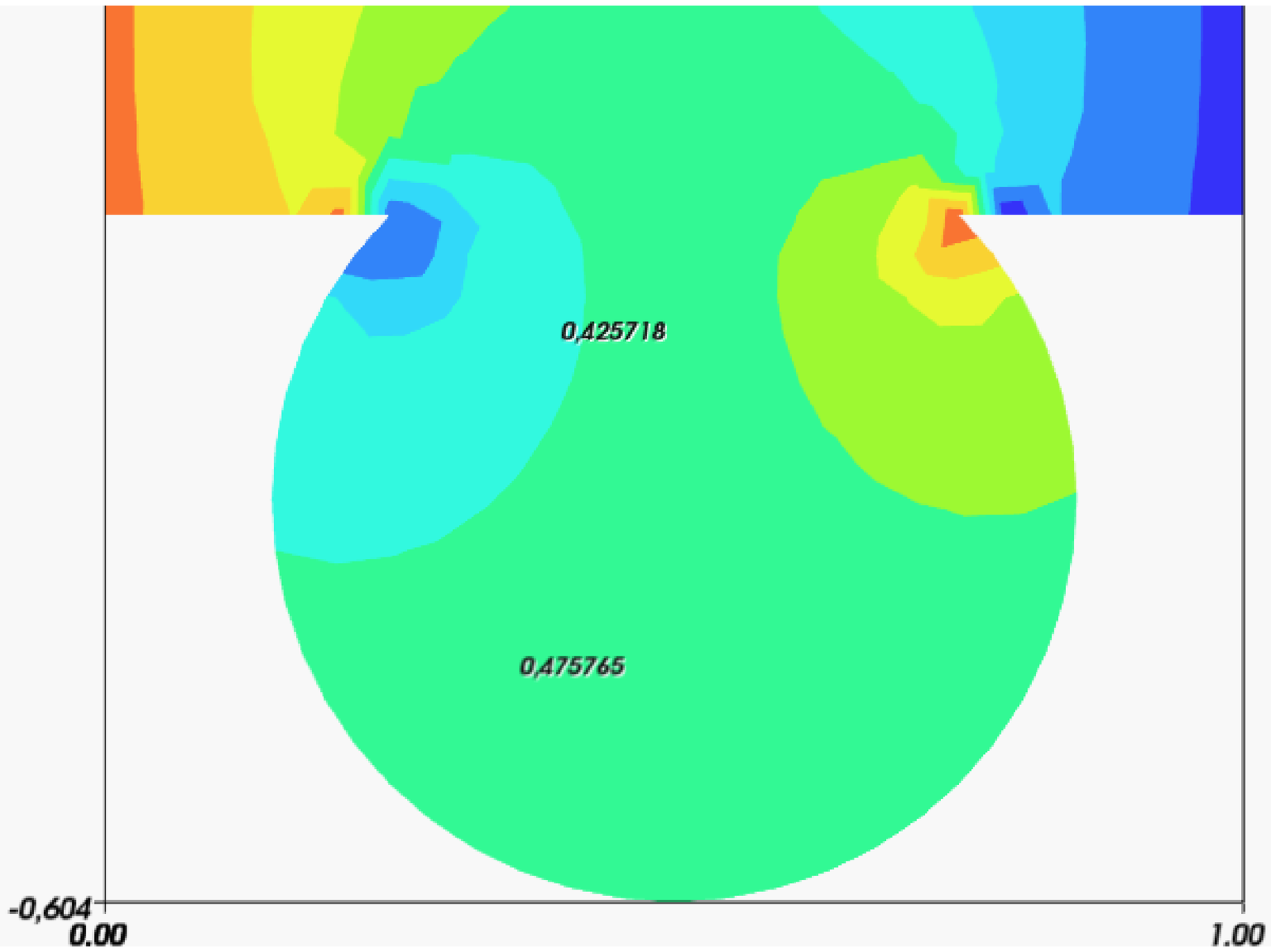}
\caption{Pressure in an aneurysmal sac,  with (left) and without a stent (right) }\label{sac_press}
\end{figure}

As stated in the corollary, we  show in the next section  that in fact 
the zero order pressure is the only constant that insures conservation 
of mass in $\Ode$. From the medical point of view the two claims on pressure and flow are of interest.
They quantify and confirm the stabilizing 
effect of a porous stent: besides reducing the stress on the wall of the
aneurysm,
the stent averages also the pressure inside the sac avoiding 
for instance corner singularities (see fig. \ref{sac_press}).

The geometry presented as an illustration in figures \ref{sac_velo} and \ref{sac_press} does not fit exactly in the hypotheses of section \ref{geo}: the main difference is the curved circular form of the boundaries of $\Ode$. Nevertheless the phenomenon observed when $\Ode$ is the square $]0,1[\times ]-1,0[$ still happens when $\Ode$ has this more physiological shape.

Again one has a mathematical validation of the formal multi-scale construction 
\begin{thm}\label{ane.main.thm}
There exists a unique pair $(\ueps,\peps)\in \bH^1(\Oe)\times L^2(\Oe)$ solving \eqref{exact_sac}.
The first order approximation $(\ovueps,\ovpeps)$ belongs to 
$\bL^2(\Omega_1\cup\Omega_2) \times H^{-1}(\Omega_1\cup \Omega_2)$, moreover we have a convergence result that reads
$$
\nrm{\ueps - \ovueps}{L^2(\Omega)} + \nrm{\peps - \ovpeps}{H^{-1}(\Op\cup \leps)} \leq k \e^{\td^-},
$$
where $\td^-$ represent any real number strictly less then $\td$,
the constant $k$ depends on the data of the problem and the domain but not on $\e$.
\end{thm}

We  show the same type of result as above :  $(\ovueps,\ovpeps)$  solve formally the same implicit problem \eqref{homo_fo} up to the second order error, but with a homogeneous Dirichlet condition on $\Goutd$.

\section{Technical preliminaries}\label{sect.techn.prelim}

In this section we introduce the basic results that allow to deal with the Stokes problem on a perforated domain $\Oe$ together with the specific boundary conditions as in \eqref{exact}.

\subsection{Weak solutions for sieve problems}

In  the spirit of Appendix in  \cite{SaPa.80} we start by the definition  a restriction operator $\Reps$ acting on functions defined in $\Omega$ and providing resulting functions defined on $\Oe$ and vanishing on $\Geps$.

\begin{definition}\label{def.restrict}
Let $\bV : = \{ \bfv \in \bH^1(\Oe) \text{ s.t. } \bfv=0 \text{ on } \Gd \text{ and } \bfv \cdot \bft = 0 \text{ on } \Gn \}$, and we endow it with the usual $\bH^1$ norm.
The tilde operator refers always to an extension by zero outside $\Oe$ i.e.
$$
\forall \bfu \in \bV, \ti{\bfu} \in \bH^1(\Omega) \text{ s.t. }\ti{\bfu}:= \left\{ 
\begin{aligned}
\bfu & \text{ if } x \in \Oe  \\
0 & \text{ otherwise }
\end{aligned}
\right.
$$ 
We define the restriction operator $\Reps \in \cL ( \bH^1(\Omega) ; \bH^1(\Oe))$ s.t.
\begin{enumerate}
\item  $\bfu \in \bV$ implies  $\Reps \ti{\bfu} \equiv \bfu$,
\item $\dive \bfu =0 $ in $\Omega$ implies that $\dive (\Reps \bfu )= 0$ in $\Oe $,
\item There exist three real constants $k_1,k_2$ and $k_3$ independent on $\e$ s.t.
$$
\left\{ 
\begin{aligned}
\nrm{\Reps \bfu}{\bL^2(\Oe)} & \leq k_1 \{ \nrm{\bfu}{\bL^2(\Omega)} + \e  \nrm{\nabla \bfu}{L^2(\Omega)^4} \} ,\\
\nrm{\nabla \Reps \bfu}{\bL^2(\Oe)} & \leq k_2 \left\{  \frac{1}{ \e } \nrm{\bfu}{\bL^2(\Omega)} + \nrm{\nabla \bfu}{L^2(\Omega)^4} \right\}, \\
\nrm{\nabla \Reps \bfu}{\bL^2(\Oe)} & \leq k_3  \frac{1}{ \sqrt{\e} } \nrm{\bfu}{\bH^1(\Omega)}. 
\end{aligned}
\right.
$$ 
\end{enumerate}
\end{definition}

\begin{lemma}\label{exist.restriction}
  There exists an opertor $\Reps$ in the sense of Definition \ref{def.restrict}
\end{lemma}

\begin{proof}
  The restriction operator $\Reps$ is constructed exactly as in Lemma 3 and 4 in the Appendix by L. Tartar in \cite{SaPa.80}, namely
for a given $\bfu \in \bH^1(\cJ)$ there exists a unique pair $(\bfv,q) \in \bH^1(\jm)\times (L^2(\jm)/ \RR)$ satisfying:
$$
\left\{ 
\begin{aligned}
& - \Delta \bfv + \nabla q = -\Delta \bfu  & \text{ in } \jm, \\
& \dive \bfv = \dive \bfu + \frac{1}{|\jm|} \int_{\js} \dive \bfu \; dy, & \text{ in } \jm\\
& \bfv = 0& \text{ on } P, \\
& \bfv = \bfu & \text{ on } \gamma_M, \\
\end{aligned}
\right.
$$
There exists a constant $k$ independent on $\bfu$ s.t. $\nrm{\bfv}{\bH^1(\jm)}\leq k \nrm{\bfu}{\bH^1(\cJ)}$. By setting
$$
R \bfu(y) := \left\{ 
\begin{aligned}
\bfu(y) &  \text{ if } y \in \cJ \setminus (\jm\cup \js)\\
\bfv(y) &  \text{ if } y \in \jm\\
0 &  \text{ if } y \in \js \\
\end{aligned}
\right.
$$
we evidently have $\nrm{ R \bfu }{\bH^1(\cJ)} \leq k \nrm{\bfu}{\bH^1(\cJ)}$. and $R \bfu$ coincides with $\bfu$ if $\bfu\equiv 0$ on $\js$ and $\dive \bfu=0$ implies $\dive (R \bfu)=0$.

Now let $\bfu \in \bH^1(\Omega)$. For any given $x\in \rlay$, we set $y_1 := x_1/ \e -  E(x_1/ \e)$, $y_2:=x_2/ \e$ and $i:=E(x_1/ \e)$ where $E(\cdot)$ is the lower integer part of its real argument.
For each $i$ we define a function $\underline{\bfu}_i:\cJ\to \RR$ s.t. $\underline{\bfu}_i(y):=\bfu(x)$. This allows us to set
 $\Reps \bfu$ as
$$
\Reps \bfu (x) := \left\{ 
\begin{aligned}
&\bfu(x) &  \text{ if } x \in \Op, \\
&\sum_{i=0}^{1/ \e} R \underline{\bfu}_i (x_1/ \e - i ,x_2/ \e) \chiu{\e (\cJ+i \eu)}(x)   &  \text{ if } x \in \rlay.\\
\end{aligned}
\right.
$$
This definition implies obviously that $\nrm{\Reps \bfu}{\bH^1(\Op)}\equiv \nrm{\bfu}{\bH^1(\Op)}$ and we focus on $\rlay$.

$$
\nrm{ \Reps \bfu }{\bL^2(\leps)}^2 = \sum_{i=0}^\ue \e^2 \int_{\jf + i \eu} | R \underline{\bfu}_i| dy \leq k \sum_{i=0}^\ue \e^2 \nrm{ \underline{\bfu}_i }{\bH^1(\cJ)}^2 \leq k \{ \nrm{\bfu}{\bL^2(\rlay)}^2+ \e^2 \nrm{\nabla_x \bfu}{L^2(\rlay)^4} \} .  
$$
The key point of the proof are now estimates on the gradient. Taking a regular function $\bfu \in \bm{\D}( \Omega)$, one obtains in a similar way as above:
\begin{equation}\label{est.rlay}
\nrm{ \nabla \Reps \bfu}{L^2(\leps)^4} \leq  k \left\{ \ue \nrm{ \bfu }{\bL^2(\rlay)} + \nrm{ \nabla \bfu }{L^2(\rlay)^4} \right\} 
\end{equation}
Now one writes 
$$
\bfu( x ) = \bfu(x_1,0) + \int_0^{x_2} \partial_{x_2} \bfu(x_1,s) ds
$$
which, after taking the square, integrating on $\rlay$ and using Cauchy-Schwartz gives 
$$
\nrm{\bfu}{\bL^2(\rlay)}^2 \leq k \left\{ \e \nrm{ \bfu}{\bL^2(\Gz)}^2 + \e^2 \nrm{\nabla \bfu}{L^2(\rlay)^4}^2 \right\} ,
$$
thanks to the continuity of the trace operator $\nrm{\gamma(\bfu)}{H^\ud(\Gamma)} \leq k \nrm{\bfu}{\bH^1(\Omega)}$, one has:
$$
\nrm{\bfu}{\bL^2(\rlay)}^2 \leq k' \left\{ \e \nrm{ \bfu}{\bH^1(\Omega)}^2 + \e^2 \nrm{\nabla \bfu}{L^2(\rlay)^4}^2 \right\}, 
$$
where the constant $k'$ does not depend on $\e$. Using this last inequality in \eqref{est.rlay} ends the proof.
\end{proof}

\begin{defi}\label{def.lift}
We define the corresponding lifting operator $\seps \bfu:= (\Reps - \id )\bfu $. 
For every $x$ in $\leps$ there exist a unique $i:=E(x/ \e) \in \{ 0,\dots,\ue \}$, $y_1 = x_1 / \e -i$ and $y_2=x_2 / \e$
s.t. if $\bfv_i$ solves 
$$
\left\{ 
\begin{aligned}
& - \Delta \bfv_i  + \nabla q = 0  & \text{ in } \jm+i \eu, \\
& \dive \bfv_i = \frac{1}{|\jm+ i \eu|} \int_{\js+i \eu} \dive \underline{\bfu}_i \; dy, & \text{ in } \jm+i \eu\\
& \bfv_i = \underline{\bfu}_i & \text{ on } P+i \eu, \\
& \bfv_i = 0 & \text{ on } \gamma_M+i \eu, \\
\end{aligned}
\right.
$$
then one sets $\seps \bfu(x) := \sum_{i =0}^{\ue} \bfv_i (y) \chiu{\jm+i \eu}$ for every $x \in \Oe$.
One has estimates similar to those of the restriction operator
$$
\begin{aligned}
\nrm{ \seps \bfu }{\bL^2(\Oe)}& \equiv \nrm{ \seps \bfu }{\bL^2(\leps)} \leq k \left\{ \nrm{ \bfu}{L^2(\rlay)} + \e \nrm{ \nabla \bfu}{L^2(\rlay)} \right \},   \\
\nrm{ \nabla \seps \bfu }{\bL^2(\Oe)}& \equiv \nrm{ \nabla \seps \bfu }{\bL^2(\leps)} \leq \left\{ \frac{1}{\e} \nrm{ \bfu}{L^2(\rlay)} + \nrm{ \nabla \bfu}{L^2(\rlay)} \right \} \leq \frac{1}{\sqrt{\e}} \nrm{ \bfu }{\bV }\quad .
\end{aligned}
$$
\end{defi}

\begin{proposition}\label{prop.surj}
  Let $g \in {\L}^2(\Oe)$ there exists at least one vector $\bfv \in \bV$ s.t.
$$
\dive \bfv = g \text{ in } \Oe,\quad \snrm{\bfv}{\bH^1(\Oe)} \leq \frac{ k}{\sqrt{\e}} \nrm{ g }{L^2(\Oe)}.
$$
\end{proposition}

\begin{proof}
  We extend $g$ by zero in $\Omega \setminus \Oe$ which we denote $\ti{g}$, we use Lemma III.3.1 and Theorem III.3.1 in \cite{Galdi}  stating that
there exists $\bfw \in H^1(\Omega)$ s.t. 
$$
\dive \bfw = \ti{g} \text{ in } \Omega ,\quad \snrm{\bfv}{\bH^1(\Omega)} \leq k \nrm{ \ti{g} }{L^2(\Omega)},
$$
where the constant $k$ does not depend on $\e$. Using the restriction operator $\Reps$ defined in the proof of Lemma \ref{exist.restriction}  one sets then
$$
\bfv := \Reps \bfw,
$$
thanks to the estimates that the restriction operator satisfies, one gets the desired result:
$$
\snrm{\bfv}{\bH^1(\Oe)} \leq \frac{k}{\sqrt{\e}} \snrm{\bfw}{\bH^1(\Omega)} \leq  \frac{k'}{\sqrt{\e}} \nrm{\ti{g}}{L^2(\Omega)} = \frac{k'}{\sqrt{\e}} \nrm{g}{L^2(\Oe)}.
$$
\end{proof}

Thanks to the latter proposition one easily gets by duality arguments as in p. 374 in the Appendix in \cite{SaPa.80}  
\begin{proposition}\label{prop.est.necas}
  There exists a constant $k$ independent on $\e$ s.t. for every distribution $p \in \D'(\Oe)$ s.t. $\nabla p \in \bV'$, one has
$$
\nrm{p}{L^2(\Oe)} \leq \frac{k}{\sqrt{\e}} \nrm{ \nabla p }{\bV'}\; .
$$ 
\end{proposition}

At this stage we can derive existence and uniqueness as well as {\em a priori} estimates for the solutions of the problem: given $(\bff,g,h) \in \bV'\times L^2(\Oe)\times H^{-\ud}(\Gn)$ find $(\bfu, p)$ s.t.
\begin{equation}\label{eq.type.pbm} 
\left\{  
\begin{aligned} 
& - \Delta \bfu + \nabla p = \bff &  \text{ in } \Oe, \\
& \dive \bfu =g &    \text{ in } \Oe, \\
& \bfu = 0  & \text{ on } \Gd, \\
& \left. 
\begin{aligned}
&  \bfu \cdot \bft = 0 \\
& p = h
\end{aligned}
\right\} & \text{ on } \Gn,
\end{aligned}  
\right. 
\end{equation}

\begin{theorem}\label{thm.gen.weak.sol}
If the data of problem \eqref{eq.type.pbm} are s.t.  
$(\bff,g,h) \in \bV'\times L^2(\Oe)\times H^{-\ud}(\Gn)$ then there exists a unique solution $(\bfu,p)\in \bV\times L^2(\Oe)$, moreover one has:
$$
\nrm{\bfu}{\bH^1(\Oe)}+ \nrm{p}{L^2(\Op)}+ \sqrt{\e}\nrm{p}{L^2(\leps)} \leq k\left\{ \nrm{\bff}{\bV'} + \frac{k}{\sqrt{\e}} \nrm{g}{L^2(\Oe)} + \nrm{ g-h}{H^{-\ud}(\Gn)} \right\} ,
$$
where the constant $k$ is independent on $\e$.
\end{theorem}

\begin{proof}
  The existence and uniqueness of $(\bfu,p)\in \bV\times L^2(\Oe)$ are standard results of the literature (see for instance \cite{ErGu.04.book} and references therein). We focus here on the control of the norms for this solution pair. Lifting the divergence source term provides easily {\em a priori} estimates on $\bfu$:
$$
\nrm{ \nabla \bfu }{L^2(\Oe)} \leq k\left\{ \nrm{\bff}{\bV'} + \frac{k}{\sqrt{\e}} \nrm{g}{L^2(\Oe)} + \nrm{ g-h}{H^{-\ud}(\Gn)} \right\} .
$$
Then we split $\Oe$ and restate problem \eqref{eq.type.pbm} on $\Op$, having for the pressure that
$$
 - \Delta \bfu + \nabla p = \bff \quad  \text{in} \quad \Op, \\
$$
which gives
$$
\nrm{\nabla p}{\bH^{-1}(\Op)} \leq k \nrm{\bff+ \Delta \bfu}{H^{-1}(\Op)} \leq k \left\{ \nrm{\bff}{\bV'} + \nrm{\bfu}{\bH^1(\Oe)} +\nrm{g}{H^{-\ud}(\Gn)} \right\} ,
$$
but on this domain there exists a constant independent on $\e$ s.t.
$$
\nrm{p}{L^2(\Op)} \leq k \nrm{ \nabla p }{\bH^{-1}(\Op)},
$$
which gives the same error estimates for the gradient of the velocity as well as for the pressure in $\Op$. Unfortunately because of the presence of the obstacles, in the rough layer one has only that
$$
\nrm{p}{L^2(\leps)} \leq \nrm{p}{L^2(\Oe)} \leq \frac{1}{\sqrt{\e}} \nrm{\nabla p}{\bV'} = \frac{1}{\sqrt{\e}} \nrm{ \bff + \Delta \bfu }{\bV'},
$$ 
which, by using again the {\em a priori} estimates of $\bfu$ in $\Oe$, gives the final estimate. 
\end{proof}

\subsection{Very weak solutions}

We recall here the framework 
of ``very weak'' solutions originally introduced in \cite{LioMagBook.1972,conca2,GaldiVws}.  

\begin{definition}\label{def.reg.dom} Let $\omega$ be an open bounded connected domain whose boundary $\partial \omega$ is split in  two disjoint parts $\partial \omega_D$ and $\partial \omega_N$.  It is said to satisfy the regularity property $\bfH^2\times H^1$ for the Stokes problem if for all $\bm{F} \in \bfL^2(\omega)$ and every $G\in H^1_0(\omega)$ the solutions of the problem 
\begin{equation}\label{reg} 
\left\{ 
\begin{aligned} 
& - \Delta \cT + \nabla \cX =   \bm{F} ,&  \text{ in } \omega, \\
&  \dive \cT  = G,& \text{ in } \omega, \\
& \cT = 0 , & \text{ on }\partial \omega_D,  \\
&\left. \begin{aligned}
& \cT\cdot \bft = 0 \\
&  \cX = 0 
\end{aligned} \right\} 
& \text{ on } \partial \omega_N,
\end{aligned}  
\right. 
\end{equation}
satisfy $\cT \in \bH^2 ( \omega) , \cX \in H^1(\omega)$ and if there exists $C_1=C_1(\omega)$ s.t.
$$
\nrm{\cT}{\bH^2 ( \omega) } + \nrm{ \cX}{H^1(\omega)} \leq C_1 \left\{  \nrm{\bm{F}}{\bfL^2(\omega)} + \nrm{\nabla G}{\bL^2(\omega)}\right\} .
$$
\end{definition}

Following exactly the same proof as in Appendix A in \cite{conca2} one shows
\begin{theorem}\label{thm.very.weak}
  If $\omega$ satisfies the regularity property of definition \ref{def.reg.dom} above, then there exists a unique solution $(\bfu,p) \in \bfL^2(\omega)\times H^{-1}(\omega)$ solving:
\begin{equation}\label{weak} 
\left\{ 
\begin{aligned} 
& - \Delta \bfu + \nabla p = \bff ,&  \text{ in } \omega, \\
& \dive \bfu = g,& \text{ in } \omega ,\\
& \left.\begin{aligned}
&  p = h \\
& \bfu\cdot \bft = \ell 
\end{aligned}\right\} 
&\text{ on } \partial \omega_N, \\
& \bfu = \bfm , & \text{ on }\partial \omega_D, 
\end{aligned}  
\right. 
\end{equation}
provided that the data satisfy:  $\bff \in \bV'(\omega)$, $g\in L^2(\omega)$, $h \in H^{-1}( \partial \omega_N )$, $\ell \in L^2(\partial \omega_N)$, $\partial_{\bft} \ell \in H^{-1}(\partial \omega_N)$  and $\bfm \in \bfL^2(\partial \omega_D)$.
Moreover there exists $C_2=C_2(\omega)$ s.t.
$$
\begin{aligned}
\nrm{\bfu}{\bfL^2 ( \omega) } + \nrm{p}{H^{-1}(\omega)} \leq  C_2 & \left\{  
\nrm{\bfm}{\bfL^2(\partial\omega_D)}
+ \nrm{\bff}{\bV'(\omega)} + \nrm{g}{L^2(\omega)} \right. \\
& \quad \left. + \nrm{g-\frac{\partial \ell}{\partial \bft} -h}{H^{-1}(\partial \omega_N)}
+ \nrm{\ell}{L^2( \partial \omega_N)}
\right\}. 
\end{aligned}
$$
where $\bV(\omega):=\{ \bfv \in \bH^1(\omega)\text{ s.t. } \bfv =0 \text{ on } \partial \omega_D \text{ and }  \bfv \cdot \bft =0 \text{ on } \partial \omega_N \}$ and $\bV'(\omega)$ is its dual. 
We denote by ``very weak'' solution such a pair $(\bfu,p)$.
\end{theorem}

\begin{theorem}\label{thm.very.weak.estimates}
If the pair $(\bfu,p)$ is a weak solution of problem \eqref{eq.type.pbm}, one has then the very weak estimates:
$$
\nrm{ \bfu }{\bL^2(\Op)} + \nrm{ p }{H^{-1}(\Op)} \leq k \left\{ \nrm{\bff}{\bV'}+ \sqrt{\e} \nrm{ \bfu }{\bH^1(\Oe)} + \nrm{ g - h }{H^{-1} (\Gn)} + \nrm{g}{L^2(\Oe)} \right\} ,
$$
where the constant $k$ is independent on $\e$. Moreover in the rough layer one has:
$$
\nrm{ \bfu }{\bL^2(\leps)}+ \nrm{ p }{H^{-1}(\leps)} \leq k \e \left\{ \nrm{ \nabla \bfu }{L^2(\Oe)^4} + \nrm{p}{L^2(\leps)} \right\}, 
$$
where again the generic constant $k$ is independent on $\e$.
\end{theorem}

\begin{proof}
  The first estimate follows by applying Theorem \ref{thm.very.weak} with $\omega:= \Op$. It is easy to show that actually in this doamin, the constants present in the very weak estimates are independent on $\e$: the obstacles are not part of $\Op$.
Thus one has:
$$\begin{aligned}  
\nrm{\bfu}{\bfL^2 ( \Op) } + \nrm{p}{H^{-1}(\Op)} & \leq k\left\{  
\nrm{\bfu}{\bfL^2(\{ x_2=0\} \cup\{ x_2=\e \})}
+ \nrm{\bff}{\bV'(\Op)} + \nrm{g}{L^2(\Op)}
+ \nrm{g-h}{H^{-1}(\Gn \cap\partial \Op)} 
\right\}, \\
& \leq k' \left\{  
\sqrt{\e} \nrm{\nabla \bfu}{\bfL^2(\leps)}
+ \nrm{\bff}{\bV'(\Oe)} + \nrm{g}{L^2(\Oe)}
+ \nrm{g-h}{H^{-1}(\Gn)} 
\right\}, \\
\end{aligned}
$$
where we used Poincar{\'e} estimates knowing that $\bfu$ vanishes on $\Geps$.
It remains to consider the rough layer $\leps$. There, we have
\begin{equation}\label{poincare_stokes}
\nrm{ \bfu}{\bL^2(\leps)} + \nrm{ p }{H^{-1}(\leps) } \leq \e \left\{ \nrm{\bfu}{\bH^1(\leps)} + \nrm{ p }{L^2(\leps) } \right\}  ,
\end{equation}
where the $\bH^1(\leps)$ regularity is  obtained using Theorem \ref{thm.gen.weak.sol}.
Indeed the estimates on the velocity come using Poincar{\'e} estimates at the microscopic level as in 
Lemma 3.2 in \cite{conca1}, 
the pressure estimate is obtained by duality:  by definition of the dual norm one has
$$
\nrm{ p}{H^{-1}(\leps) } = \sup_{\varphi \in H^1_0(\leps)} <p,\varphi>_{  H^{-1} ,H^1_0 },
$$
where $H^1_0(\leps)$ denotes the set of functions in $H^1(\leps)$ vanishing on $\Gz\cup]0,1[\times \{\e\}\cup \partial \leps $. As $p$ belongs to $L^2(\leps)$
the duality bracket can be transformed into an integral, namely
$$
\begin{aligned}
<p,\varphi>_{  H^{-1} ,H^1_0 } & = \int_{\leps} \, p \, \varphi \, dx \leq \nrm{ p }{L^2(\leps) } \nrm{\varphi}{L^2(\leps)}
\leq \e \nrm{ p }{L^2(\leps) } \nrm{\varphi}{H^1(\leps)},  
\end{aligned}
$$
 taking the sup over all functions in $H^1_0(\leps)$, one concludes the norm correspondence. 
\end{proof}

\section{Proof of the main results}\label{proofs}

\subsection{The case of a collateral artery}


In what follows we set both $\poutu$ and $\poutd$ to be zero for simplicity.
The results remain valid for any fixed  constants $\poutu$ and $\poutd$ as well.

\subsubsection{The zero order term}

When $\e$ goes to zero  we show in a first step that $(\ueps,\peps)$
converges to $(\bfuz,\pz)$ the  Poiseuille profile stated in \eqref{poiseuille.col},
which solves in $\Oe$~:
$$
\left\{
\begin{aligned}
& -\Delta \bfuz + \nabla \pz = [\sigma_{\bfuz,\pz}]\cdot\bfn \, \delta_{\Gz} & \text{ in } \Oe,\\ 
& \dive \bfuz = 0 & \text{ in } \Oe,\\
&  \bfuz = 0 & \text{ on } \gud, \\ 
&  \bfuz\cdot \bft  = 0 &\text{ on } \Gn,\\ 
&  \pz = \pin \text{ on }  \Gin, \quad  \pz = 0  & \text{ on }  \Goutu\cup \Goutd,\\ 
& \bfuz \neq 0 & \text{ on } \Geps.
\end{aligned}
\right.
$$

\begin{thm}\label{zo_col} For any fixed $\e$, there exists a unique solution $(\ueps,\peps)\in \bV \times L^2(\Oe)$ of the problem \eqref{exact}. Moreover, one has
$$
\nrm{\ueps-\bfuz}{\bH^1(\Oe)} + \nrm{\peps - \pz}{L^2(\Op)} + \sqrt{\e} \nrm{\peps - \pz}{L^2(\leps)} \leq k \sqrt{ \e },
$$
where the constant $k$ does not depend on $\e$. On the other hand one can prove that:
$$
\nrm{\ueps-\bfuz}{\bL^2(\Oe)} + \nrm{\peps - \pz}{H^{-1}(\Op\cup\leps)} \leq k \e 
$$
\end{thm}

\begin{proof}
Existence and uniqueness of the solutions of problem \eqref{exact}
come from the standard theory of mixed problems \cite{Galdi,ErGu.04.book}, Theorem \ref{thm.gen.weak.sol}
gives more precisely
$$
\nrm{\ueps}{\bH^1(\Oe)} + \nrm{\peps}{L^2(\Op)} + \sqrt{\e} \nrm{\peps}{L^2(\leps)}  \leq k \nrm{\pin}{H^{-\ud}(\Gin)},
$$
where the constant $k$ is independent on $\e$.
As $\bfuz$ does not satisfy homogeneous boundary conditions we use the restriction operator 
already presented in the section above. Namely we set:
$$
\hbfu := \ueps - \Reps \bfuz,\quad \hp := \peps - \pz,
$$
these variables solve:
$$
\left\{ 
\begin{aligned}
& - \Delta \hbfu + \nabla \hp = - \Delta ( \bfuz - \Reps \bfuz ) + [\sigma_{\bfuz,\pz}]\cdot\bfn \, \delta_{\Gz}  = \Delta (\seps \bfuz) + [\sigma_{\bfuz,\pz}]\cdot\bfn \, \delta_{\Gz}  & \text{ in } \Oe, \\
& \dive \hbfu = 0 & \text{ in } \Oe, \\
& \hbfu = 0 & \text{ on } \Gd, \\
& 
\left.
\begin{aligned}
\hbfu \cdot \bft = 0 \\
\hp = 0 
\end{aligned}
\right\} & \text{ on } \Gn,\\
\end{aligned}
\right.
$$
where the lifting operator $\seps$ is given in Definition \ref{def.lift}.
Thanks to Theorem \ref{thm.gen.weak.sol}, one has then directly:
$$
\nrm{ \nabla \hbfu }{\bL^2(\Oe)}+ \nrm{ \hp }{ L^2(\Op) } + \sqrt{ \e } \nrm{ \hp }{ L^2(\leps) } 
\leq  \nrm{ \Delta \seps \bfuz + [ \sigma_{\bfuz,\pz} ] \cdot \bfn \dgz}{\bV'} \\
$$
Thanks to the vicinity of $\Geps$, one  deduces easily some trace inequalities \cite{conca1}:
\begin{equation}\label{ineeq.trace.poincare}
\nrm{\Psi}{L^2(\Gz)} \leq \sqrt{ \e} \nrm{\nabla \Psi}{L^2(\Oe)^2}, \quad \forall \Psi \in H^1(\Oe) \text{ s.t. } \Psi=0 \text{ on } \Geps,
\end{equation}
this estimate allows us to conclude that
$$
\sup_{\bm{\Psi}\in \bV } \int_{\Gz} ([ \sigma_{\bfuz,\pz} ] \cdot \bfn, \bm{\Psi}) dx_1 \leq \nrm{ [ \sigma_{\bfuz,\pz} ]\cdot \bfn}{ \bfL^2(\Gz)} \nrm{\bm{\Psi}}{\bL^2(\Gz)} \leq  \sqrt{\e}  \nrm{ [ \sigma_{\bfuz,\pz} ]\cdot \bfn}{ \bfL^2(\Gz)} \nrm{\bm{\Psi}}{\bV}.
$$
The specific form of the lifting $\seps \bfuz$  allows to write:
$$
\nrm{ \Delta (\seps \bfuz) }{\bV'} = \nrm{ \nabla (\seps \bfuz) }{L^2(\leps)} \leq \left\{ \frac{1}{\e} \nrm{ \bfuz}{L^2(\rlay)} + \nrm{ \nabla \bfuz}{L^2(\rlay)} \right \} \leq k \sqrt{ \e },
$$
where we used the explicit form of the Poiseuille profile in the rough layer.
Using Theorem \ref{thm.very.weak.estimates} one has then
$$
\nrm{ \hbfu }{\bL^2(\leps)} + \nrm{ \hp }{H^{-1}(\leps)} \leq k \e .
$$
One has then easily also that 
$$
\nrm{ \ueps - \bfuz  }{\bL^2(\leps)} \leq \nrm{ \ueps - \Reps \bfuz  }{\bL^2(\leps)} + \nrm{ \Reps \bfuz  - \bfuz }{\bL^2(\leps)} \leq \nrm{ \hbfu}{L^2(\leps)}  + \nrm{ \seps \bfuz }{L^2(\leps)} \leq k \e .
$$
\end{proof}

Estimates above show a threefold error:
the Dirichlet error on $\Geps$, 
the jump of the gradient of the velocity 
in the horizontal direction across $\Gz$,
and the pressure jump across $\Gz$.
In order to correct these errors we solve three  microscopic 
boundary layer problems. 

\subsubsection{The Dirichlet correction}
The first boundary layer corrects the  Dirichlet error on $\Geps$. It is very alike to the one introduced 
in the wall-laws setting \cite{JaMiSIAM.00,AchPiVaJCP.98,BrMiQam}. Namely we solve the problem:
find $(\bfbeta,\pi)$ such that
\begin{equation}
\label{beta.cell}
\left\{
\begin{aligned}
& -\Delta \bfbeta + \nabla \pi = 0   &\text{ in } Z, \\ 
& \dive \bfbeta = 0 &\text{ in } Z,\\
&  \bfbeta = -y_2 \eu  &\text{ on } P, \\ 
& \beta_2 \to 0   &|y_2| \to \infty,\\
& (\bfbeta,\pi) \text{ are } 1-\text{periodic in the }  y_1 \text{ direction}.
\end{aligned}
\right.
\end{equation}
We define as in \cite{Galdi} p. 56, the homogeneous Sobolev space
$\bmdu:=\{ \bfv \in \D'(Z),\text{ s.t. } \nabla \bfv \in (L^2(Z))^4 \}$.
Moreover we denote by $\bmduz$ the subset of functions belonging to $\bmdu$
and vanishing on $P$.

\begin{prop}
\label{cell.pbm.beta.prop} There exists a unique solution $(\bfbeta,\pi) \in \bmdu\times\ldl$, $\pi$ being 
defined up to a constant.
Moreover, one has:
$$
\bfbeta(y) \to \obeta^\pm_1  \eu , \quad y_2 \to \pm \infty 
$$
the convergence being exponential with rate $\gamma_\beta$
and 
$$
\left\{ 
\begin{aligned}
& \obeta_2(y_2) = 0, &  \forall y_2 \in \RR \setminus ]0,\ydp[, \\
&\obeta_{1}(y_2) = -|\js| - |\nabla \bfbeta|_{L^2(Z)^4}^2 +  \obeta_1(0), & \forall y_2>y_{2,P} , \\
& \obeta_1(y_2)=\obeta_1(0), &  \forall y_2 < 0,
\end{aligned}
\right.
$$
where $y_{2,P}:= \max_{y\in P} y_2$ and $|\js|$ is the 2d-volume of the obstacle $\js$
\end{prop}
For sake of conciseness the proof is given in the Appendix  \ref{periodic.bl}.

\subsubsection{Shear rate jump correction}
The second boundary layer corrects the jump of the normal derivative
of the axial velocity: we introduce a 
source term that accounts for a unit jump in the horizontal component but on the microscopic
scale. Namely, we look for $(\bfups,\varpi)$ solving:
\begin{equation}\label{eq.upsi}
\left\{
\begin{aligned}
& -\Delta \bfups + \nabla \varpi = \delta_\Sigma \eu&  \text{ in } Z, \\ 
& \dive \bfups = 0 &  \text{ in } Z,\\
&  \bfups = 0 &\text{ on } P, \\ 
& \Upsilon_2 \to 0 &  |y_2| \to \infty ,\\
& (\bfups,\varpi ) \text{ are } 1-\text{periodic in the } y_1 \text{ direction}.
\end{aligned}
\right.
\end{equation}
Again we give some basic results and the behaviour at infinity 
of this corrector.
\begin{prop}\label{prop.upsi} There exists a unique $(\bfups,\varpi )\in \bmduz\times\ldl$, 
$\varpi$ being defined up to a constant.
Moreover, one has:
$$
\bfups (y) \to \ooupsis_1^\pm   \eu , \quad y_2 \to \pm \infty ,
$$
and 
$$
\left\{ 
\begin{aligned}
&\ov{\Upsilon}_2(y_2) = 0 &  \forall y_2 \in \RR, \\
&\ov{\Upsilon}_{1}(y_2) =\ov{\Upsilon}_1(0)+ \obeta_1(0) & \forall y_2>\ydp , \\
&\ov{\Upsilon}_1(y_2)=\ov{\Upsilon}_1(0)  = \nrm{\nabla \bm{\Upsilon}}{L^2(Z)^4} &  \forall y_2 < 0, \\
\end{aligned}
\right.
$$
where $\ydp:= \max_{y\in P} y_2$.
\end{prop}

The reader finds again the proof in Appendix \ref{periodic.bl}.

\subsubsection{The pressure jump}
In order to cancel the pressure jump $[\pz]$, we use a corrector similar to the one
introduced and widely studied for a flat sieve in \cite{conca1} p. 25:
\begin{equation}\label{chi.cell}
\left\{
\begin{aligned}
& -\Delta \bfchi + \nabla \eta =0 &\text{ in } Z,\\ 
& \dive \bfchi = 0 &\text{ in } Z,\\
&  \bfchi = 0 &\text{ on } P, \\
& \chi_2 \to - 1 , & |y_2| \to \infty\\
& (\bfchi,\eta ) \text{ are } 1-\text{periodic in the } y_1 \text{ direction}.
\end{aligned}
\right.
\end{equation}
As in the proof of Proposition \ref{cell.pbm.beta.prop}, 
one repeats the arguments of Appendix \ref{periodic.bl} in order to obtain similarly to \cite{conca1}~:
\begin{prop}
There exists a unique solution $(\bfchi,\eta)\in \bmdu\times\ldl$ of system \eqref{chi.cell}, $\eta$ being
defined up to a constant.
Moreover, one has
$$
\bfchi\to \ochi \equiv -\ed,\quad |y_2| \to \infty, 
$$
the convergence being exponential with rate $\gamma_\bfchi$ and there exists 
two constants $\oeta({+\infty})$ and $\oeta({-\infty})$  depending only on the geometry of $P$
such that
$$
\eta(y) \to \oeta(\pm\infty),\quad |y_2| \to \infty.
$$ 
One then proves:
$$
\left| \nabla \bfchi \right|^2_{L^2(Z)} = [\ooeta].
$$
\end{prop}
This corrector will be used in the sequel, but we already utilize
it to give a first result on the average of $\pi$ and $\varpi$
\begin{coro}\label{prop.press.infty}
The solutions $(\bfbeta,\pi)$ and $(\bfupsi,\varpi)$ solving respectively \eqref{beta.cell} and \eqref{eq.upsi} satisfy~:
$$
\opi(y_2)= 0 \text{ and } \ov{\varpi}(y_2)= 0, \quad \forall y_2 \in \RR_-\cup]\ydp,+\infty[
$$
\end{coro}
For the proof see again Appendix \ref{periodic.bl}.

As explained in Remark \ref{rmk.so.bl}  below, 
we  need in section \ref{sec.fo.approx} a higher order  corrector 
that solves the problem: find $(\bm{\varkappa},\mu)$ s.t.
$$
\left\{
\begin{aligned}
& -\Delta \bm{\varkappa}+ \nabla \mu = - 2 (\nabla \bfchi - (\eta-\ooeta)\id ).\eu & \text{  in } Z,\\ 
& \dive \bm{\varkappa} = 0 & \text{  in } Z,\\
&  \bm{\varkappa} = 0 &\text{ on } P,\\ 
&  \varkappa_2 \to 1 &\quad |y_2| \to \infty,\\
& (\bm{\varkappa}, \mu ) \text{ are } 1-\text{periodic in the } y_1 \text{ direction}.
\end{aligned}
\right.
$$
\begin{prop}There exists a unique solution $(\bm{\varkappa},\mu)\in\bmdu\times\ldl$, $\mu$ being defined
up to a constant. One has also exponential convergence towards constants with rate $\gamma_\varkappa$ :
$$
\bm{\varkappa} \to \oov{\bm{\varkappa}}, \quad \mu \to \oomu, \text{ when } |y_2| \to \infty.
$$
Moreover one has the relationships between values at $y_2=\pm \infty$
$$
[\oomu] =  [\ooeta]- 2 \int_{Z} \chi_1 (\eta-\ooeta)) dy , \quad [\ov{\ov{\varkappa}}_1] = - 2 \int_Z ( \sigma_{\bfchi,(\eta-\ooeta)}.\eu,\bfbeta + y_2 \eu ) dy.  
$$
\end{prop}
The proof is  exactly the same as for Propositions \ref{cell.pbm.beta.prop} and \ref{prop.upsi} and thus is left to the reader.

In what follows we  use the $\e$-scaling of all boundary layers above, namely 
we set:
$$
\beps (x) := \bfbeta\lrxe,\quad \bm{\Upsilon}_\e  (x) := \bm{\Upsilon} \lrxe,\quad \chieps (x):= \bfchi \lrxe,\quad \varkappaeps(x) := \bm{\varkappa} \lrxe, \quad \forall x \in \Oe.
$$
the same notation holds for pressure terms as well.
\subsubsection{Vertical correctors on $\gio\cup\Gde$}

Above boundary layers are  periodic; 
their oscillations perturb homogeneous Dirichlet %
as well as Neumann stress boundary conditions on $\gio\cup\Gde$. 
The perturbation on these boundaries is $O(1)$, 
due to the vicinity of these edges to the geometrical
perturbation $\Geps$.
We introduce  
vertical boundary correctors defined on a 
half-plane $\Pi$. 
Each of them accounts for  perturbations induced by  
the periodic boundary layers on $\gio\cup\Gde$ in the very vicinity of corners $O$ and $\ov{x}$. These correctors solve at the microscopic scale 
the problems:
\begin{equation}\label{vertical_micro_corr}
\left\{
\begin{aligned}
& -\Delta \wb+ \nabla \tb  = 0 & \text{ in } \Pi, \\ 
& \dive  \wb = 0 & \text{ in } \Pi, \\
& \wb  = -(\bfbeta-\oobeta\lambda) & \text{ on } D, \\ 
& \wbd = -\beta_2  & \text{ on } N, \\ 
& \tb = -\pi &\text{ on } N,\\ 
\end{aligned}
\right. \quad 
\left\{
\begin{aligned}
& -\Delta \wu+ \nabla \tu  = 0&\text{ in } \Pi, \\ 
& \dive  \wu = 0 &\text{ in } \Pi,\\
& \wu  = -(\bm{\Upsilon}-\oov{\bm{\Upsilon}}\lambda) & \text{ on } D, \\ 
& \wud = -\Upsilon_2 &\text{ on } N, \\ 
& \tu = - \mu&\text{ on } N,\\ 
\end{aligned}
\right.
\left\{
\begin{aligned}
& -\Delta \wc+ \nabla \tc  = 0& \text{ in } \Pi, \\ 
& \dive  \wc = 0 & \text{ in } \Pi,\\
& \wc  = -(\bfchi-\oochi \lambda) &\text{ on } D, \\ 
& \wcd = -(\chi_2 -\oochid \lambda) &\text{ on } N, \\ 
& \tc = -(\eta-\ooeta) &\text{ on } N,\\ 
\end{aligned}
\right.\quad 
\end{equation}
and  $(\wk,\tk)$ solves a similar system 
lifting  $(\oovarkappa \lambda -\bm{\varkappa},\oomu-\mu)$ on $D\cup N$.
Note that the domain $\Pi$ does not contain any obstacles, 
we should use a restriction operator on the velocity 
vectors $\bfw_i$ in order to handle this feature (see below \eqref{eq.def.bfw}).
We define the usual weighted Sobolev space \cite{Ha.71,AmGiGiI.94}, for all $(m,p,\alpha) \in \N\times [1,\infty[ \times \RR$:
$$
\ws{m}{p}{\alpha}{\Pi}:=\left\{ v \in \cD'(\Pi) \, \text{ s.t. } \, |D^\lambda v| \rho^{\alpha + |\lambda| - m{\nobreakspace}} \in L^p(\Pi),\, 0\leq |\lambda| \leq m \,\right\}, 
$$
where $\rho:=(1+|y|^2)^\ud$. We endow this space with the corresponding 
weighted norm. By density arguments one proves that dual spaces of $ \ws{m}{p}{\alpha}{\Pi}$ are distributions and we set in the rest of this work
$$
\ws{-m}{p}{-\alpha}{\Pi} :=( \ws{m}{p}{\alpha}{\Pi} )', \quad \quad \forall (m,p,\alpha) \in \N\times[1,\infty[   \times \RR.
$$
Here we extend results obtained for
mixed boundary conditions and the rough Laplace equation
in \cite{BrBoMi, MiVws} to the case of the Stokes equations. 
In  the appendix we give the extensive proof  of the crucial claim:
\begin{thm}\label{MiStokesUnbounded}
Thanks to the exponential decrease to zero of the boundary data in \eqref{vertical_micro_corr},
 there exists a unique solution $(\bfw_i,\theta_i) \in \bws{1}{2}{\alpha}{\Pi}^2\times \ws{0}{2}{\alpha}{\Pi}$ for $i\in\{ \bfbeta,\bfupsi,\bfchi,\bm{\varkappa}\}$, for every real $\alpha$ s.t.
$|\alpha| < 1$.
\end{thm}

\begin{rmk}
The  weight exponent $\alpha$ provided by this result on the microscopic scale is  important. It accounts for the behaviour  when $\rho$ goes infinity of the vertical correctors above. The decay properties so described are used in Lemma \ref{lem.decr.micro.vert} in order to quantify,  in terms of powers of $\e$, the impact of the  perturbation induced by the periodic correctors on the  macroscopic lateral Dirichlet and Neumann boundary conditions: the greater $\alpha$ the smaller the error in terms of powers of $\e$. So we assume $\alpha$ very close to 1.
\end{rmk}

\bigskip

The Poiseuille profile admits an explicit form \eqref{poiseuille.col}
and thus its  derivative wrt $x_2$ reads 
$ \partial_{x_2} u_{0,1}= (\pin - \poutu)(1-2 x_2)/2 \chiu{\Ou}$. For the rest of the paper we implicitly assume  $ \partial_{x_2} u_{0,1}$  to be evaluated at $x_2=0^+$: it is constant and reads
\begin{equation}\label{def.duzu}
 \duzu :=  \duzu (x_1,0^+) = \frac{\pin-\poutu}{2}= \frac{\pin}{2}.
\end{equation}
We set for $i\in\{ \bfbeta,\bfups, \bfchi,\bm{\varkappa}\}$,
\begin{equation}\label{eq.def.bfw}
\left\{
\begin{aligned}
\bfw_{\e,i}(x) &:= c_i \Reps \bfw_i \lrxe \psi_1(x_1)+ \ti{c}_i \Reps \ti{\bfw}_i\left(\frac{x-\ov{x}}{\e}\right)\psi_2(x), \\
 \theta_{\e,i}(x)&  := c_i \theta_i \lrxe \psi(x_1)+\ti{c}_i\ti{\theta}_i\left(\frac{{x}-\ov{x}}{\e}\right)\psi_2(x),\quad 
\end{aligned}
\right.
\end{equation}
where we used the restriction operator $\Reps$ of definition \ref{def.restrict}, 
while  $(\ti{\bfw}_i,\ti{\theta}_i)$ solve similar problems as \eqref{vertical_micro_corr}
but on the  halfspace $\Pi^-:= \RR_- \times \RR$, and
the constants $c_i$ (resp. $\ti{c}_i$) denote
$$
\begin{aligned}
& c_\bfbeta := \duzu(O)  ,\quad c_\bfups := \lrduzu (O),\quad c_\bfchi:=\dpde(O) ,\quad c_{\bm{\varkappa}} := \pin, \\
& \ti{c}_\bfbeta := \duzu(\ov{x})  ,\quad \ti{c}_\bfups := \lrduzu (\ov{x}),\quad \ti{c}_\bfchi:=\dpde(\ov{x}) ,\quad \ti{c}_{\bm{\varkappa}} := \pin.
\end{aligned}
$$
For the particular explicit zero order solution $(\bfuz,\pz)$ expressed in \eqref{poiseuille.col},
$c_\bfchi$ is the only constant for which $c_i\neq \ti{c}_i$. As the analysis carried below on  $\Gin\cup\Gde$ is exactly 
the same on  $\Goutu\cup\Gde$ we implicitly assume that when terms appear containing constants $c_i,\psi_1,\bfw_i$  and $\theta_i$
similar expressions with $\ti{c}_i,\psi_2,\ti{\bfw}_i$ and $\ti{\theta}_i$ are considered as well.

\begin{lemma}\label{lem.decr.micro.vert}
Defining vertical correctors $(\bfw_{\e,i},\theta_{\e,i})$ with $i \in \{ \bfbeta,\bfupsi,\bfchi,\bfvarkappa \}$ one has the  estimates:
\begin{itemize}
\item In the whole domain:
$$
\nrm{ \Delta (\e \bfw_{\e,i}) - \nabla \theta_{\e,i} }{\bV'} \leq k \e^{1^-} ,\quad \nrm{ \e \dive \bfw_{\e,i} }{L^2(\Oe)} \leq k \e^{1^-},
$$
where the constant $k$ is independent on $\e$, and $1^-$ is any constant strictly less than 1.
\item In $\Op$, one has:
$$
\nrm{ \Delta (\e \bfw_{\e,i}) - \nabla \theta_{\e,i} }{\bH^{-1}(\Op)} \leq k \e^{1+\alpha} ,\quad \nrm{ \e \dive \bfw_{\e,i} }{L^2(\Op)} \leq k \e^{1+\alpha},
$$
where $\alpha$ is any positive real smaller than 1, and $k$ is again independent on $\e$.
\item On $\Gn$  one has 
$$
\nrm{ \e \dive \bfwepsi }{ H^{-\ud}( \Gn )} \leq k \exp\left(-\ue\right).
$$
\item In $\Omega_1 \cup \Omega_2$ one has:
$$
\nrm{ \bfwepsi }{ \bL^2(\Omega_j) } \leq k \e^\alpha ,\quad j \in \{ 1,2 \}. 
$$
\end{itemize}

\end{lemma}

\begin{proof}
  An easy calculation shows that for any $i \in \{ \bfbeta, \bfupsi,\bfchi, \bfvarkappa\} $ one has:
$$
\begin{aligned}
J& := - \Delta (\e \bfwepsi) + \nabla \thepsi = - \Delta_x \left(\e \bfw_{i} \lrxe \psi \right) + \nabla_x  \left( \theta_i \lrxe \psi \right) + \Delta_x \left( \e \psi \seps \bfw_i \lrxe\right) \\
& = - \left( 2 \nabla_y \bfw_i\lrxe -  \theta_i\lrxe  \id \right) \nabla_x \psi - \left( \e \bfw_i \lrxe  \Delta_x \psi \right) +  \Delta_x \left( \e \psi \seps \bfw_i \lrxe \right) 
\end{aligned}
$$
One has then that
$$
\begin{aligned}
\nrm{J}{\bV'} & \leq k \nrm{ \left( 2 \nabla_y \bfw_i\lrpe -  \theta_i\lrpe  \id \right) \nabla_x \psi - \left( \e \bfw_i \lrpe  \Delta_x \psi \right) }{\bL^2(\Oe)} \\
& \quad + \nrm{\Delta_x \left( \e \psi \seps \bfw_i \lrpe \right) }{\bV'} =: I_1 + I_2
\end{aligned}
$$ 
We split $I_1$ in two parts $I_{1,1}$ and $I_{1,2}$. The first part is estimates as:
$$
\begin{aligned}
I_{1,1}^2 &\leq \e^2 \int_\Pi (|\nabla \bfw_i|^2 + \theta_i^2) |\nabla \psi_1|^2 dy \leq k \e^2 \int_{\Pi \cap \left]\frac{1}{3\e},\frac{2}{3\e}\right[\times ]0,\pi[} (|\nabla \bfw_i|^2 + \theta_i^2)  dy  \\
& \leq \e^2 \nrm{ |\nabla \bfw_i| + \theta_i}{\ws{0}{2}{\alpha}{\Pi}}^2 \sup_{r\in\left] \frac{1}{3 \e},\frac{2}{3 \e}\right]} \rho^{-2\alpha} \leq k \e^{2(1+\alpha)}.
\end{aligned}
$$
Similarly the second part is estimated as  well by
$$
\begin{aligned}
I_{1,2}^2 = & \e^2 \int_{\Oe}|\bfw_i(x/ \e)|^2 (\Delta_x \psi_1(x) )^2 dx  = \e^4 \int_\Pi |\bfw_i|^2 (\Delta_x \psi_1(\e y) )^2 dy \\
& \leq k \e^4 \int_{\Pi \cap \left]\frac{1}{3\e},\frac{2}{3\e}\right[\times ]0,\pi[} |\bfw_i|^2 r dr d\tilde{\theta} \\
&  \leq k \e^4 \left(\int_{\Pi \cap \left]\frac{1}{3\e},\frac{2}{3\e}\right[\times ]0,\pi[} \left(\frac{|\bfw_i|}{\rho}\right)^2 \rho^{2\alpha} r dr d\tilde{\theta}  \right)\cdot \left( \sup_{r \in \left[\frac{1}{3 \e},\frac{2}{3 \e}\right]} \rho^{2-2\alpha} \right)\\
& \leq \e^{2(1+\alpha)}  k \nrm{\bfw_i}{\bws{1}{2}{\alpha}{\Pi}}^2.\\
\end{aligned}
$$
Passing then to $I_2$ one has easily that
\begin{equation}\label{eq.final}
\begin{aligned}
\nrm{ \e \Delta \left( (\seps \bfw_{i}) \lrpe  \psi \right) }{\bV'} & = \nrm{ \e \nabla_x \left((\seps \bfw_{i})\lrpe  \psi\right) }{L^2(\Oe)^4} \\ 
& \leq \nrm{ \e (\seps  \bfw_{i})\otimes \nabla \psi}{L^2(\Oe)^4} + \nrm{ \psi \nabla_y (\seps  \bfw_{i})\lrpe }{L^2(\Oe)^4} 
\end{aligned}
\end{equation}
Now thanks to the estimates on the lift
$$
\begin{aligned}
\nrm{\seps \bfw_i\lrpe}{\bL^2(\Oe)}^2 &\leq \e^2 \nrm{\seps \bfw_i }{\bL^2(B(O,\ue))}^2 
 \leq \e^2 \left\{ \nrm{  \bfw_i }{\bL^2(B(O,\ue))}^2  + \nrm{ \nabla  \bfw_i }{L^2(B(O,\ue))^4}^2 \right\} \\
& \leq \e^{2 \alpha} \nrm{\bfw_i}{\bws{1}{2}{\alpha}{\Pi}}^2 + \e^2 \nrm{\bfw_i}{\bws{1}{2}{0}{\Pi}}^2
\end{aligned}
$$
and the same way one gets:
$$
\begin{aligned}
\nrm{ \nabla_y \seps \bfw_i \lrpe }{L^2(\Oe)}^2 & \leq k \e^2 \nrm{ \nabla_y \seps \bfw_i }{L^2(B(O,\ue))^4}^2  \leq \e^2 \left\{ \nrm{  \bfw_i }{\bL^2(B(O,\ue))}^2  + \nrm{ \nabla_y  \bfw_i }{L^2(B(O,\ue))^4}^2 \right\} \\
& \leq k\left\{ \e^{2 \alpha} \nrm{\bfw_i}{\bws{1}{2}{\alpha}{\Pi}}^2 + \e^2 \nrm{\bfw_i}{\bws{1}{2}{0}{\Pi}}^2 \right\} 
\end{aligned}
$$
Putting together last two estimates in \eqref{eq.final} one obtains the first result of the claim. The result on the divergence follows the same lines.

On  $\Op$ the result is more straightforward since $\Reps \bfw_i \equiv \bfw_i$ for all $i \in \{ \bfbeta,\bfupsi,\bfchi,\bfvarkappa \}$. Then using again the correspondance between macroscopic powers of $\e$ and microscopic weighted spaces one gets easily the result.

On $\Gn$ one has that
$$
\e \dive \bfwepsi = \e \dive \left( \bfw_i \lrxe \psi \right)= \e \nabla \psi \cdot \bfw_i \lrxe  = \e \partial_{\bft} \psi ( \bfw_i \cdot \bft )
$$
because $\partial_{\bfn} \psi \equiv 0$ on this boundary.
On the microscopic scale the support of $\partial_{\bft} \psi$ is located in $\frac{1}{3 \e}$ and $\frac{2}{3 \e}$, thus  one has
$$
\nrm{ \e  \dive \bfwepsi }{ H^{-\ud}(\Gn)} \leq \nrm{ \e  \dive \bfwepsi }{ L^2(\Gn)} \leq \e \exp\left(- \ue \right)
$$
\end{proof}

Then, we define the complete vertical corrector as 
$$
\left\{
\begin{aligned}
\cW(x) &:= \e \sum_{i\in\{ \bfbeta,\bfups, \bfchi\}} \bfw_{\e,i}(x) + \e^2 \bfw_{\e,\bm{\varkappa}} + \bfW(x),\\
\cZ(x) & := \sum_{i\in\{ \bfbeta,\bfups, \bfchi\}} \theta_{\e,i} + \e^2 \theta_{\e,\bm{\varkappa}} + S(x),
\end{aligned}
\right.\quad \forall x \in \Oe 
$$
where $(\bfW,S)$ solve the system of equations on  the macroscopic domain $\Oe$~:
\begin{equation}\label{vertical_macro}
\left\{
\begin{aligned}
&\Delta \bfW + \nabla S = 0, & \text{ in } \Oe, \\ 
&\dive \bfW = 0 & \text{ in } \Oe,\\ 
& \left.
\begin{aligned}
\bfW\cdot\bft   &= \e \left\{ c_\bfbeta (\beps - \oov{\bfbeta})  + c_\bfupsi (\bfupsieps - \oov{\bfupsi} )  + c_\bfchi (\chieps - \oov{\bfchi} ), 
+ \e c_{\bm{\varkappa}} (\varkappaeps-\oovarkappa) \right\}  \cdot\bft \, \Phi \\
 S &=  \left\{ c_\bfbeta \pi_\e  + c_\bfupsi \varpi_\e  + c_\bfchi (\eta_\e -\ooeta)  + \e c_{\bm{\varkappa}} (\mueps-\oomu) \right\}\, \Phi \\
\end{aligned} \right\} 
  \quad & \text{  on } \Gn,\\ 
& \bfW = \e \left\{ c_\bfbeta (\beps-\oobeta) + c_\bfupsi (\bfupsieps - \oov{\bfupsi} )  + c_\bfchi (\chieps -\oochi) + \e c_{\bm{\varkappa}} (\varkappaeps - \oovarkappa) \right \} \Phi & \text{ on } \Gd \\ 
\end{aligned}
\right.
\end{equation}
where one notes  that $\cW\equiv 0$ on $\Geps$ because of the support of $\Phi$.
\begin{prop}\label{macro_corr_prop}
There exists a unique solution $(\bfW,S)\in \bH^1(\Oe)\times L^2(\Oe)$ of system \eqref{vertical_macro}, moreover one has:
$$
\nrm{ \bfW }{\bH^1(\Oe)} + \nrm{S}{L^2(\Oe)} \leq k e^{-\frac{\gamma}{\e}}
$$
where the exponential rate $\gamma$ and the constant $k$ do not depend on $\e$.
\end{prop}

\begin{proof}
Setting 
$$
\cR := 	\left( \e \left\{ \duzu (\beps - \oobeta) 
	   + \lrduzu  (\bfupsieps - \oov{\bfupsi}) 
	   + \dpde (\chieps-\oochi) \right\} 
 	   + \e^2 \pinp (\varkappaeps-\oovarkappa) \right) \Phi, 
$$
and $\ti{\bfW}:= \bfW - \cR$, by the standard theory for mixed 
problems \cite{Galdi,ErGu.04.book}, there exists a
unique solution $(\ti{\bfW},S)$ of the lifted problem. 
One has also {\em a priori} 
estimates :
$$
\nrm{ \ti{\bfW} }{\bH^1(\Oe)} + \nrm{ S }{L^2(\Oe)} \leq \nrm{ \Delta \cR }{ \bH^{-1}( \Oe )} + \nrm{ \dive \cR }{L^2( \Oe )} + \nrm{S}{H^{\ud}(\Gn)} 
$$
Thanks to the crucial presence of the cut-off function $\Phi$ 
and the exponential decrease
of rate $\gamma:=\min(\gamma_\bfbeta,\gamma_\bfupsi,\gamma_\bfchi,\gamma_{\bm{\varkappa}})$ of all the microscopic correctors, 
one gets  the exponential decrease of the rhs in the previous estimates. 
Because it is also trivial to show that $\nrm{\cR}{\bH^1(\Oe)} \leq k e^{-\frac{\gamma}{\e}}$ one ends the proof.
\end{proof}


\bigskip

\subsubsection{The complete first order approximation}\label{sec.fo.approx}
Having introduced every single element, 
we built a complete first order approximation. 
We define the full boundary layer corrector:
\begin{equation}\label{fo.bl.cor}
\begin{aligned}
\cU := \bfuz &+ \e \left\{  \duzu (\beps - \oobeta ) + \lrduzu (\bfupsieps - \oov{\bf{\Upsilon}}) + \frac{[\pz]}{[\ooeta]}( \chieps  -\oochi ) + \bfu_1 \right\}  \\
 & + \e^2 \left\{ \pinp (\varkappaeps - \oovarkappa ) +  \bfud \right\} + \cW,\\
 \cP  := \pz &+ \left\{ \duzu \pi_\e  + \lrduzu \varpi_\e + \frac{[\pz]}{[\ooeta]}( \eta_\e - \ooeta ) + \e p_1 \right\} + \e \pinp ( \mu_\e  - \oomu ) + \e^2 \pd + \cZ,
\end{aligned}  
\end{equation}
where the normal derivative $\partial_{x_2} u_{0,1}$ is defined in \eqref{def.duzu} and where the first order and second order macroscopic correctors $(\bfuu,\pu)$ and $(\bfud,\pd)$ solve respectively \eqref{fo.so.macro} and
\begin{equation}\label{so.macro}
\left\{
\begin{aligned}
-&\Delta \bfud + \nabla \pd = 0 & \text{ in } \Omega_1\cup \Omega_2,\\ 
& \dive \bfud = 0& \text{ in } \Omega_1\cup \Omega_2,\\
&  \bfud = 0 &\text{ on } \Gamma_D,\\
& \left.
\begin{aligned} 
&  \bfud\cdot \bft = 0 \\
&  \pd = 0
\end{aligned}\right\}
& \text{ on }  \Gamma_N,\\ 
& \bfud = \pinp \oovarkappa^\pm & \text{ on } \Gz^\pm. 
\end{aligned}
\right.
\end{equation}
Problems \eqref{fo.so.macro} and \eqref{so.macro} are defined on two separate domains $\Ou$ and $\Ode$: 
two distinct values are given as Dirichlet boundary conditions to the horizontal component of the velocity on $\Gz$. 
This is due to the different values of the constants whom the boundary layer correctors $\bfbeta$ and $\bfupsi$ tend to at + and - infinity. 
Thus the velocity vectors $\bfuu$ and $\bfud$ are not only discontinuous across $\Gz$ but also multi-valued at the corners $O$ and $\ov{x}$. 
It appears then clearly that
$\bfuu$ and $\bfud$ cannot belong to $\bH^1(\Ou\cup\Ode)$. For this reason, we use the 
concept of very weak solution introduced in the section above.
Because $\Oup$ and $\Ode$ are  convex polygons in $\RR^2$, they fulfill  regularity conditions of definition \ref{def.reg.dom} (see example 2.1 p 53 in \cite{conca2}). This allows to use Theorem  \ref{thm.very.weak} in order to obtain
\begin{corollary}
The pairs of functions $(\bfuu,\pu)$ and $(\bfud,\pd)$ solving problems \eqref{fo.so.macro} and \eqref{so.macro} exist and are unique ``very weak'' solutions in $\bL^2(\Ou\cup\Ode) \times H^{-1}(\Ou\cup\Ode)$ .
\end{corollary}

For some technical reasons appearing later on, one needs to set up an intermediate pair of functions $(\tiuud, \tipud)$  solving~:
\begin{equation}\label{eq.tiuud.tipud} 
\left\{  
\begin{aligned} 
& - \Delta \tiuud + \nabla \tipud = 0,&\text{ in } \Ou\cup\Ode \\
& \dive \tiuud = 0& \text{ in } \Ou\cup\Ode \\
& \tiuud = - \left( \duzu \oobeta + \lrduzu \ooupsi +  \dpde \oochi \right) \ldd & \text{ on } \Gun \cup \Gde \\
& 
\left. 
\begin{aligned}
& \tiuud \cdot \bft = - \left( \duzu \oobeta + \lrduzu \ooupsi +  \dpde \oochi \right) \ldd \cdot \bft \\ 
& \tipud = 0 \\
\end{aligned} 
\right\} & \text{ on } \gioud \\
& \tiuud \equiv 0& \text{ on } \Gz
\end{aligned}  
\right. 
\end{equation}
We define implicitly the same second order problem whose solutions we denote in the same fashion $(\tiutd,\tiptd)$.
These are  regularized versions of problem \eqref{fo.so.macro} and \eqref{so.macro} where we lifted the constants from the interface $\Gz$. 
Indeed the function $\bfuu$ is multi-valued in the corners $0$ and $\ov{x}$ multiplying the specific cut-off function $\ldd$ on these corners insures that:
\begin{proposition}\label{prop.h1.macro.fo.so}
  There exists a unique solution $(\tiuud, \tipud)$  solving \eqref{eq.tiuud.tipud}. Moreover one has the estimates:
$$
\nrm{ \tiuud }{\bH^1(\Omega)} + \nrm{ \tipud }{L^2(\Omega)} \leq k \left\{\sqrt{|\log \e | }+ \e \right\} .
$$
Taking the restricition to $\Oe$ of $ \tiuud$ one has:
$$
\nrm{ \Reps \tiuud }{\bH^1(\Oe)} \leq k \nrm{ \nabla \seps \tiuud }{L^2(\Oe)^4} \leq k  \sqrt{|\log \e | },
$$
where the generic constant $k$ is independent on $\e$.
\end{proposition}

\begin{proof}
We denote 
$$
\bfg := \left( \duzu \oobeta + \lrduzu \ooupsi +  \dpde \oochi \right) \ldd.  
$$
  On each subdomain by lifting the Dirichlet data one obtains in a classical way 
$$
\begin{aligned}
\nrm{ \tiuud }{\bH^1(\Ou \cup \Ode)} &+ \nrm{ \tipud }{L^2(\Ou \cup \Ode)} \leq k \left\{ \nrm{ \bfg }{\bH^{-1}(\Ou \cup \Ode)} + \nrm{ \dive \bfg }{L^2(\Ou \cup \Ode)} + \nrm{\dive \bfg}{\bH^{-\ud}(\Gin \cup \Goutu)} \right \} \\
& = k \left\{ \nrm{ \nabla \bfg }{L^2(\Ou \cup \Ode)} + \nrm{ \dive \bfg }{L^2(\Ou \cup \Ode)}  + \nrm{\dive \bfg}{\bH^{-\ud}(\Gin \cup \Goutu)}\right \}   \\
& \leq  k \left\{ \nrm{ \nabla \bfg }{L^2(\Ou \cup \Ode)} + \nrm{ \dive \bfg }{L^2(\Ou \cup \Ode)}  + \nrm{\dive \bfg}{\bL^2(\Gin \cup \Goutu)} \right \}   \\
& \leq k \left\{ \nrm{ \ldd}{H^1(\Omega)} + 2 |\oochi_2| \nrm{1}{L^2(0,\e)} \right \} 
\end{aligned}
$$
where the constant $k$ does not depend on $\e$. 
Note that the latter equality is true since $\partial_{\bfn} \ldd$ vanishes on the boundaries of $\Omega$ and $\oochi \cdot \eu=0$.
Now using the $H^1$ estimate on the gradient of $\ldd$ \eqref{eq.lift.h1},  
one recovers the first {\em a priori} estimate.
Working on $\Oeu$, when writing the system that $(\Reps \tiuud ,\tipud)$ solve, one gets easily that
$$
\nrm{ \Reps \tiuud }{\bH^1(\Oeu)} \leq k \nrm{ \nabla \seps \tiuud }{L^2(\Oeu)^4} \leq k' \left\{ \ue \nrm{  \tiuud }{L^2(\rlay)} + \nrm{  \nabla \tiuud }{L^2(\rlay)} \right\} 
$$
and because $\tiuud \equiv 0$ on $\Gz$ one uses the Poincar{\'e} inequality on the first term in the rhs above in order to obtain:
$$
\nrm{ \Reps \tiuud }{\bH^1(\Oeu)} \leq k'' \nrm{  \nabla \tiuud }{L^2(\rlay)} \leq k''' \sqrt{ |\log \e|}
$$
On $\Ode$, there are no obstacles i.e.  $\Reps \tiuud \equiv \tiuud$ which ends the proof.
\end{proof}

\bigskip

\begin{remark}\label{rmk.so.bl}
In the error estimates developed in the next sections, one applies  the momentum operator to the term:
\begin{equation}\label{terme.croix}
\frac{[\pz]}{[\ooeta]}(\e ( \chieps  -\oochi ), \eta-\ooeta).   
\end{equation}
Because $[\pz]$ depends on $x$, 
the rest is not zero: among others a  $O(1)$ double product of gradients remains, it reads
$$
2 \nabla_y ( \chieps  -\oochi )\cdot \frac{\nabla [\pz]}{[\ooeta]} =  2 \frac{\partial_{x_1} \pz }{[\ooeta]} \nabla_y ( \chieps  -\oochi )\cdot \eu =  - 2 \frac{ \pin}{[\ooeta]} \nabla_y ( \chieps  -\oochi )\cdot \eu. 
$$
This term could be estimated directly in the $L^2$-norm, giving 
$$
\nrm{\nabla_y ( \chieps  -\oochi )\cdot \eu}{\bfL^2(\Oe)} \leq k \sqrt{\e}
$$
which is a zeroth order error.
Needing better error estimates, we add 
the second order term in the asymptotic ansatz \eqref{fo.bl.cor} that reads: 
$$
- \partial_{x_1} \pz ( \e^2 (\varkappaeps - \oovarkappa) , \e (\mueps-\oomu)) 
$$
The Stokes operator applied to this corrector cancels exactly the double product above. 
The divergence of the velocity part of \eqref{terme.croix} gives as well a cross term : 
this is already of order $\e^\td$ in the $L^2$ norm, so we should not  correct it. 
\end{remark}

\subsubsection{{\em A priori} estimates}

We consider here a complete boundary layer approximation containing a regularized macroscopic correctors $(\tiuud,\tipud)$ and $(\tiutd,\tiptd)$ reading:

\begin{equation}\label{eq.bdl.reg.delta}
\begin{aligned}
\cUd := \bfuz &+ \e \left\{  \duzu \beps + \lrduzu \bfupsieps  + \frac{[\pz]}{[\ooeta]} \chieps  + \tiuud \right\}  + \e^2 \left\{ \pinp \varkappaeps +  \tiutd \right\} + \cW,\\
 \cPd  := \pz &+ \left\{ \duzu \pi_\e  + \lrduzu \varpi_\e + \frac{[\pz]}{[\ooeta]}( \eta_\e - \ooeta ) + \e \tipud \right\} + \e \pinp \mu_\e  + \e^2 \tiptd + \cZ.
\end{aligned}  
\end{equation}

Note that the problems \eqref{vertical_micro_corr} that the vertical correctors $\bfwepsi$ solve, account for perturbations induced by the periodic boundary layers $\beps, \bfupsieps, \chieps,\varkappaeps$ and $\tiuud,\tiutd$ on $\gio\cup \Gde$, this explains the presence of $\lambda$ in the definition of boundary terms in \eqref{vertical_micro_corr}. We recall that these correctors are included in the global corrector $\cW$.

\begin{theorem}\label{thm.a.priori.est}
The full boundary layer approximation $(\cUd,\cPd)$ defined in \eqref{eq.bdl.reg.delta} is a first order approximation of the exact solution $(\ueps,\peps)$ i.e.
$$
\nrm{ \ueps - \cUd}{\bH^1(\Oe)}+ \nrm{ \peps - \cPd }{L^2(\Op)} + \sqrt{\e} \nrm{ \peps - \cPd }{L^2(\leps)}  \leq k \e^{1^-}
$$
where the constant $k$ is independent on $\e$.
\end{theorem}

\begin{proof}
  Set $\hbfu := \ueps - \cUd$  and $\hp := \peps - \cPd$, they satisfy:
$$
\left\{
\begin{aligned}
& -\Delta \hbfu + \nabla \hp=  \sum_{i \in \{ \bfbeta,\bfupsi,\bfchi \} } \e \Delta \bfwepsi - \nabla \thepsi  + O\left(\e^{\td}\right) & \text{ in } \Oe, \\ 
& \dive \hbfu = \sum_{i \in \{ \bfbeta,\bfupsi,\bfchi \} } \e \dive \bfwepsi + O\left( \e^{\td}\right) & \\
& \hbfu  = 0 & \text{ on } \Gun \cup \Gde, \\ 
&  \hbfu  = -\left(\bfuz - \duzu x_2 \eu + \e \tiuud + \e^2 \tiutd \right) =: \bfg& \text{ on } \Geps, \\
&  \left.\begin{aligned}
&  \hbfu \cdot \bft = 0 \\ 
&  \hp =  0 \\
\end{aligned} \right\} & \text{ on } \Gn.\\ 
\end{aligned}
\right.
$$
There are two kind of errors: the first is due to the localisation of vertical correctors and is treated thanks to Lemma \ref{lem.decr.micro.vert}, the second is due to the macroscopic approximations that do not satisfy the Dirichlet condition on $\Geps$.
For the latter,  setting 
$\hhbfu = \hbfu - \seps  \bfg$ one has that
$$
\left\{
\begin{aligned}
& -\Delta \hhbfu + \nabla \hp=  \sum_{i \in \{ \bfbeta,\bfupsi,\bfchi \} } \e \Delta \bfwepsi - \nabla \thepsi  - \Delta \seps \bfg+ O\left(\e^{\td}\right) & \text{ in } \Oe, \\ 
& \dive \hhbfu = \sum_{i \in \{ \bfbeta,\bfupsi,\bfchi \} } \e \dive \bfwepsi + O\left( \e^{\td}\right) & \text{ in } \Oe,\\
& \hhbfu  = 0 & \text{ on } \Gd, \\ 
&  \left.\begin{aligned}
&  \hhbfu \cdot \bft = 0 \\ 
&  \hp =  0 \\
\end{aligned} \right\} & \text{ on } \Gn.\\ 
\end{aligned}
\right.
$$
Thanks to the explicit form of $\bfuz - \duzu x_2 \eu$ and to Proposition \ref{prop.h1.macro.fo.so} one deduces that
$$
\nrm{  \Delta \seps \bfg }{\bH^{-1}(\Oe)} = \nrm{  \nabla  \seps \bfg }{L^2(\Oe)^2} \leq k \left\{  \e^\td + (\e + \e^2) |\log \e |^\ud  \right\} ,
$$
combining this with the results of Lemma \ref{lem.decr.micro.vert} and using Theorem \ref{thm.gen.weak.sol} one gets the desired result.
\end{proof}

\subsubsection{Very weak estimates}\label{sec.very.weak}
We use here the framework of very weak solutions introduced above.
The essential motivation comes from the
lack of regularity of the averaged approximation $(\ovueps,\ovpeps)$ across the interface $\Gz$ 
and the boundary layers' optimal cost   in the $\bL^2\times H^{-1}$ norm.
The roughness $\Geps$ is contained inside the 
limiting domain $\Ou$: we decompose our 
domain in three parts $\Oup$, $\leps$ and $\Ode$. 
\begin{thm}\label{full}
The full approximation $(\cU,\cP)$ satisfies the error estimates:
$$
\nrm{\ueps-\cU}{\bL^2(\Omega_1\cup\Omega_2)} + \nrm{\peps - \cP}{H^{-1}(\Omega_1'\cup \leps\cup \Ode )} \leq k \e^{\td^-},
$$
where the constant $k$ does not depend on $\e$ and
$\td^-$  is a real strictly less than $\td$.
\end{thm}

\begin{proof}
  One gets thanks to Theorem \ref{thm.very.weak}, very weak estimates on $\hbfu:= \ueps - \cUd$ and $\hp := \peps - \cPd$:
$$
\nrm{ \hbfu }{\bL^2(\Op)} + \nrm{ \hp }{H^{-1}(\Op)} \leq k \left\{ \sqrt{\e} \nrm{ \hbfu }{H^1(\Oe)} + \nrm{ \Delta \e \bfwepsi - \nabla \thepsi}{\bH^{-1}(\Op)} + \nrm{\e \dive \bfwepsi }{L^2(\Op)} + O(\e^2) \right\}  
$$
and 
$$
\nrm{ \hbfu }{\bL^2(\leps)} + \nrm{ \hp }{H^{-1}(\leps)} \leq k \e \left\{ \nrm{ \hbfu }{\bH^1(\leps)} + \nrm{ \hp }{L^2(\leps)} \right \} .
$$
Thanks to Lemma \ref{lem.decr.micro.vert} and Theorem \ref{thm.a.priori.est} one has finally when  gathering both inequalities:
$$
\nrm{ \hbfu }{\bL^2(\Oe)} + \nrm{ \hp }{H^{-1}(\Op\cup \leps)} \leq k \e^{\td^-} 
$$
By a triangular inequality one obtains:
$$
\begin{aligned}
\nrm{\ueps -\cU }{\bL^2(\Ou\cup \Ode)} + \nrm{ \peps - \cP}{H^{-1}(\Ou'\cup \leps \cup \Ode)} & \leq  
\nrm{\ueps -\cUd }{\bL^2(\Ou\cup \Ode)} + \nrm{\cUd - \cU}{\bL^2(\Ou\cup \Ode)} \\ 
& + \nrm{ \peps - \cPd}{H^{-1}(\Ou'\cup \leps \cup \Ode)} + \nrm{ \cPd-\cP }{H^{-1}(\Ou\cup\Ode) } 
\end{aligned}
$$
Next, setting $\cbfu:=\cUd - \cU= \e(\tiuud - ( \bfuu - \bfg ))$ and $\cpp:=\cPd-\cP = \e(\tipud - \pu)$, these variables solve:
$$
\left\{
\begin{aligned}
& -\Delta \cbfu + \nabla \cpp =   0 & \text{ in } \Omega_1 \cup \Omega_2 ,\\ 
& \dive \cbfu = 0& \\ 
& \cbfu  =  \e \bfg (1- \ldd) & \text{ on } \Gde ,\\ 
& 
\left. 
\begin{aligned}
& \cbfu  \cdot \bft  = \e \bfg (1- \ldd) \cdot  \bft \\ 
& \cpp  = 0 \\
 \end{aligned}
\right\} 
& \text{ on } \Gin \cup \Goutu \cup \Goutd,\\
& \cbfu = 0, & \text{ on } \Gun \cup \Gz, \\
\end{aligned}
\right.
$$
where again 
$$
\bfg := \duzu \oobeta + \lrduzu \ooupsi +  \dpde \oochi  .
$$
Using then again the very weak estimates of Theorem \ref{thm.very.weak.estimates}  one gets
$$
\begin{aligned}
\nrm{ \cbfu }{\bL^2(\Omega_1 \cap \Omega_2)} + \nrm{ \cpp }{H^{-1}(\Omega_1 \cap \Omega_2)} \leq & k'\e \nrm{ \bfg (1- \ldd) }{\bL^2(\Gin \cup \Goutu \cup \Gde)} + k'' \e \nrm{ \partial_{\bft} \ldd }{H^{-1}(\Gn)} 
\leq  k''' \e^\td.
\end{aligned}
$$
We detail here only the second ter of the rhs above. If there exists $\hbar \in H^1_0(\Gn)$ s.t.
\begin{equation}\label{pde.hmu.nrm.gn}
- \frac{\partial^2 \hbar}{ \partial \bft^2} = \frac{\partial \ldd }{ \partial \bft} = - \ue \chiu{[0,\e]}(x_2), \quad \forall x_2 \in ]0,1[
\end{equation}
then ones has easily that
$$
\nrm{\frac{\partial \ldd }{ \partial \bft}  }{H^{-1}(\Gn)}  \leq \nrm{ \frac{\partial\hbar  }{ \partial \bft}}{L^2(\Gn)}.
$$
As $\gio$ are  straight segments, an easy computation gives that if 
$$
\hbar(x_2) := \frac{ x_2^2}{2 \e } \chiu{[0,\e]}(x_2) + \frac{ \e( 1-x_2 )}{2(1 - \e)} \chiu{[\e,1]}(x_2), \quad \forall x_2 \in [0,1]
$$
then $\hbar \in H^1_0(\gio)$ and it solves \eqref{pde.hmu.nrm.gn}. An explicit computation gives that
$$
\nrm{ \frac{\partial\hbar  }{ \partial \bft}}{L^2(\Gn)} \leq k (\sqrt{\e} + \e ),
$$
which ends the proof.
\end{proof}

\begin{rmk} 
We are not allowed to apply the very weak framework to $\Oeu$: even for $C^\infty$ obstacles,  $\Oeu$ does not satisfy uniformly wrt $\e$ the regularity property of definition \ref{def.reg.dom}. Thus we applied the very weak estimates above the rough layer in $\Oup$, this latter domain satisfying the regularity requirement of definition \ref{def.reg.dom} uniformly in $\e$. In the $\leps$ zone we use the Poincar{\'e} inequality to obtain the desired convergence rate. This explains why at last we obtain convergence results for the pressure terms  in the $H^{-1}(\Ou'\cup \leps \cup \Ode)$ norm which is smaller that the $H^{-1}(\Ou\cup\Ode)$ norm used in the case of a flat sieve (cf. p.50-52 in \cite{conca2}). 
\end{rmk}

Here we consider the oscillating part of our approximation. We recall that
$\ovueps:=\bfuz+\e\bfuu$ and $\ovpeps:=\pz+\e\pu$, and we set
$$
\oveps := \cU - \ovueps ,\quad \ovqeps := \cP - \ovpeps.
$$
The functions $\oveps,\ovqeps$ are explicit sums   
 of all the correctors in \eqref{fo.bl.cor}.
In order to prove error estimates we need 
the following two  results 
\begin{proposition}\label{prop.h.moins.un}
  If a periodic function $\tp$ is harmonic on $Z_{-\infty,0}$ and on $Z_{1,+\infty}$ and tends to zero when $|y_2|$ goes to $\infty$,  then setting $\tpeps= \tp(x/ \e)$, one has
$$
\nrm{\tpeps}{H^{-1}(\Ou')} \leq k \e^{\td} ,\quad \nrm{\tpeps}{H^{-1}(\Ode)} \leq k \e^{\td} ,
$$
where the constant $k$ is independent on $\e$.
\end{proposition}

\begin{proof}
  We prove the result for $\Ou'$, the proof is the same for $\Ode$. 
As $\tp$ is periodic and harmonic in $Z_{1,\infty}$,  it is explicit in terms of Fourier series:
$$
\tp= \sum_{n\in Z^*} \tp_n e^{-2\pi|n|y_2 + i 2\pi l y_1},\quad \forall y \in Z_{1,+\infty}, \quad \tp_n = \int_0^1 \tp(y_1,0) e^{-i2\pi|n|y_1} dy_1.
$$
It is thus decreasing exponentially fast towards 0.
Then we solve the problem : find $q$ s.t.
\begin{equation}\label{eq.q} 
\left\{  
\begin{aligned} 
& - \Delta q = \tp & \text{ in }  Z_{1,+\infty}, \\
& q= 0 & \text{ on } \{ y_2=1 \}, \\
& q \text{ is 1-periodic in the }y_1 \text{ direction}.
\end{aligned}  
\right. 
\end{equation}
Thanks to the exponential decrease of $\tp$ it is easy to show that it belongs to  $\bmdud(Z_{1,+\infty})'$ and thus
by the Lax-Milgram theorem, there exists a unique $q \in \bmdud(Z_{1,+\infty})$ solving \eqref{eq.q}. One can even 
decompose $q$ in Fourier modes and obtain again that it is an exponentially decreasing  to zero at infinity.  
Then we set $q_\e:=q(x/ \e)$, and we have
$$
- \Delta_x (\e^2 q_\e) = \tpeps \quad  \text{ in } \Ou'.
$$
Given $\varphi \in H^1_0(\Ou')$, we aim at computing 
$$
J(\varphi):= \int_{\Ou'} \tpeps \varphi dx = \int_{\Ou'} -\Delta_x (\e^2 \qeps )  \varphi dx  =\e^2 \int_{\Ou'} \nabla_x \qeps \cdot \nabla_x \varphi dx,
$$
One has immediately because of the microscopic structure of $\qeps$
$$
J(\varphi) = \e \int_{\Ou'} \nabla_y \qeps \cdot \nabla \varphi dx \leq \e^\td \nrm{ \nabla_y q }{L^2( Z_{1,+\infty} )} \nrm{\varphi}{H^1(\Ou')}
$$
the result follows writing that $\nrm{\qeps}{H^{-1}(\Ou')} = \sup_{ \varphi \in H^{1}_0(\Ou') } (J(\varphi)/ \nrm{\varphi}{H^1(\Ou')})$
\end{proof}

For the vertical correctors of pressure terms one has in the same way:
\begin{proposition}\label{prop.h.moins.un.theta}
  For a given $\theta \in \ws{0}{2}{\alpha}{\Pi}$ s.t. $\alpha \in ]0,1[$, setting $\theta_\e:=\theta(x/ \e)$ one has that
$$
\nrm{\theta_\e \psi }{H^{-1}(\Ou')} \leq k \e^{1+\alpha} , \quad \quad \nrm{\theta_\e \psi }{H^{-1}(\Ode)} \leq k \e^{1+\alpha }, 
$$
where the constant $k$ is independent on $\e$.
\end{proposition}

\begin{proof}
  We restrict ourselves to the case of $\Ou'$ again. 
We solve at the microsopic level:
\begin{equation}\label{eq.q.theta} 
\left\{ 
\begin{aligned}
& - \Delta q = \theta  , & \text{ in }  \rr \times ]1,+\infty[, \\
& q= 0 & \text{ on } \{ 0 \} \times ]1,+\infty[  \cup \rr \times \{ 1 \}.
\end{aligned}
\right.
\end{equation}
With arguments similar to those of the proof of Proposition \ref{well_posed},
 one can show that if $\theta$ is in $\ws{-1}{2}{\delta}{\rr \times ]1,+\infty[}$
with $\delta \in ]-1;1[$ then there exists a unique solution $q \in   \ws{1}{2}{\delta}{\rr \times ]1,+\infty[}$ 
solving  \eqref{eq.q.theta}. An easy computation 
shows that if $\theta \in  \ws{0}{2}{\alpha}{\Pi}$ then $\theta \in  \ws{-1}{2}{\alpha-1}{\Pi}$, which implies setting $\delta=\alpha-1$ 
the existence of a solution $q \in \ws{1}{2}{\alpha-1}{\Pi}$ provided that $\alpha \in ]0,2[$. 
As, by the definition of $\theta$, $\alpha \in ]0,1[$, 
we restrict ourselves to solutions $q \in   \ws{1}{2}{\alpha-1}{\Pi}$ with $\alpha \in ]0,1[$.
 Again we set $\qeps=q(x / \e)$ which means  that
$$
- \Delta_x (\e^2 q_\e) = \tetaeps \text{ in } \Ou'.
$$
Given a test function $\varphi \in H^1_0(\Ou')$, we aim at computing 
$$
J(\varphi):= \int_{\Ou'} \tetaeps \psi \varphi dx = \int_{\Ou'} -\Delta_x (\e^2 \qeps ) \psi  \varphi dx  =\e^2 \int_{\Ou'} \nabla_x \qeps \cdot \nabla_x ( \psi \varphi ) dx.
$$
Because of the microscopic structure of $\qeps$ one has again
$$
J(\varphi) = \e \int_{\Ou'} \nabla_y \qeps \cdot  \nabla_x ( \psi \varphi ) dx \leq \e \nrm{\psi}{W^{1,\infty}(\Ou')} \nrm{ \nabla_y \qeps }{L^2( \Ou' )} \nrm{\varphi}{H^1(\Ou')} ,
$$
passing from the macro to the micro scale we have
$$
\begin{aligned}
\nrm{ \nabla_y \qeps }{L^2( \Ou' )} 
\leq \left( \e ^2 \int_0^\ue \int_1^\ue  | \nabla q|^2 \rho ^{2\alpha-2} dy \sup_{\rho \in B(0,\ue)} \rho^{2-2\alpha} \right)^\ud 
\leq \e^{\alpha} k \nrm{ q }{\ws{1}{2}{\alpha-1}{\rr \times ]1,+\infty[ }} 
\end{aligned}
$$
by similar arguments as in Lemma \ref{lem.decr.micro.vert}.
Again the result follows writing that $\nrm{\qeps}{H^{-1}(\Ou')} = \sup_{ \varphi \in H^{1}_0(\Ou') } (J(\varphi)/ \nrm{\varphi}{H^1(\Ou')})$
\end{proof}

\begin{thm}\label{oscil}
The rapidly oscillating rest $(\cU - \ovueps,\cP - \ovpeps)$ satisfies 
$$
\nrm{\cU - \ovueps}{\bL^2(\Omega)}+ \nrm{\cP - \ovpeps}{H^{-1}(\Ou'\cup \leps\cup\Ode)} \leq k \e^{\td}
$$
where the constant $k$ is independent on $\e$.
\end{thm}
\begin{proof}
Because $\oveps$ is explicit and reads~:
$$
\oveps =  \e \left\{ \duzu (\bfbeta_\e  -\oov{\bfbeta}) + \lrduzu (\bm{\Upsilon}_\e-\oov{\bfupsieps}) 
+   \frac{[\pz]}{[\ooeta]} (\bm{\chi}_\e-\oov{\bfchi}) + \e \pinp (\bm{\varkappa}_\e-\oov{\bm{\varkappa}}) \right\} + \e \bfw_{\e,i} + \bfW,
$$
a direct computation of the $L^2$ norm gives that
$$
\begin{aligned}
\nrm{\oveps}{\bL^2(\Omega_j)} \leq& \e k \left\{ \nrm{\bfbeta_\e-\oov{\bfbeta}}{\bL^2(\Omega_j)} 
+ \nrm{\bm{\chi}_\e-\oov{\bfchi}}{\bL^2(\Omega_j)} +  \nrm{\bm{\Upsilon}_\e-\oov{\bfupsi}}{\bL^2(\Omega_j)}
+ \e \nrm{\bm{\varkappa}_\e-\oov{\bm{\varkappa}}}{\bL^2(\Omega_j)} \right\} \\
& + \e  \nrm{\bfw_{\e,i}}{\bL^2(\Omega_j)} + \nrm{\bfW}{\bL^2(\Omega_j)} \leq k \e^{\td}.
\end{aligned}
$$
We use again the 
decomposition of $\Oe$ in subdomains $\Omega_1,\leps$ and $\Omega_2$. The $y_1$-periodic pressures $\pi_\e ,\varpi_\e  ,\eta_\e -\ooeta$ and $\mueps$ fulfill hypotheses of Proposition \ref{prop.h.moins.un}, the vertical correctors $\theta_i$ for $i\in \{ \bfbeta, \bfupsi,\bfchi \}$ satisfy hypotheses of Proposition \ref{prop.h.moins.un.theta} one then concludes 
$$
\nrm{\ovqeps}{H^{-1}(\Omega_1')} \leq k \e^\td , \quad \nrm{\ovqeps}{H^{-1}(\Ode)} \leq k \e^\td.
$$
where the pressure correctors as $S$ and $\e$ terms are implicitly treated by a direct estimates of the $L^2$ norm.  
In $\leps$ we use the dual estimate \eqref{poincare_stokes} based on the  Poincar{\'e} inequality,
to get 
$$
\nrm{\ovqeps}{H^{-1}(\leps)} \leq k \e \nrm{\ovqeps}{L^2(\leps)} \leq k \e \nrm{\ovqeps}{L^2(\Oe)} \leq k \e^\td .
$$
\end{proof}

Combining Theorems \ref{full} and \ref{oscil} above, one gets the proof of Theorem \ref{main.thm.paper}.

\bigskip

\subsubsection{Implicit interface conditions}

We start with the horizontal velocity. 
We call $\bfuu^\pm$ (resp $\partial_2 u_{0,1}^\pm$) the values above and below $\Gz$. The first
order interface condition derived above on $\Gz$ reads:
$$
\bfuu^\pm = \left\{ \duzu^+ \oobetas_1^\pm + \lrduzu \ooupsis_1^\pm \right\} \eu + \dpde \oochid \ed,
$$
assembling together normal derivatives of the velocity
on both sides and because $\partial_{x_2} u_{0,1}^-\equiv 0$, one has also :
$$
\left\{ 
\begin{aligned}
\bfuu^+ &= \left\{ \duzu^+ (\oobetas^+_1 +\ooupsis^+_1 ) \right\} \eu + \dpde \oochid \ed,  \\
\bfuu^- &= \left\{ \duzu^+ (\oobetas^-_1 +\ooupsis^-_1 ) \right\} \eu + \dpde \oochid \ed,
\end{aligned}
\right.
$$
which finally  gives 
$$
\frac{\bfuu^+ \cdot \eu }{\oobetas_1^+ +\ooupsis_1^+ }= \duzu^+  \text{ and }  \frac{\bfuu^+ \cdot \eu }{\oobetas_1^+ +\ooupsis_1^+ }= \frac{\bfuu^- \cdot \eu }{\oobetas_1^- +\ooupsis_1^- }.
$$
Setting $\ovueps:=\bfuz + \e \bfuu$ and because $\bfuz\equiv0$ on $\Gz$, one has
also
$$
\ovueps^+\cdot \bft  = \e (\oobetas_1^+ +\ooupsis_1^+ ) \dd{\ovuepsu}{x_2} + O(\e^2), \text{ and }  \frac{\ovueps^+\cdot \bft    }{\oobetas_1^+ +\ooupsis_1^+ }= \frac{\ovueps^-\cdot \bft  }{\oobetas_1^- +\ooupsis_1^- }.
$$
One recovers a  slip velocity condition in the main artery and a new
discontinuous relationship between the horizontal components of the velocity 
at the interface.

For the vertical velocity, thanks to the continuity
of $\oochi_2$ across $\Gz$, one  has that 
$$
\bfuud^+ = \bfuud^- = \bfuud = -\frac{ [ \pz ]}{[\ooeta]}=\frac{ ([\sigma_{\bfuz,\pz}] \cdot \ed , \ed )}{[\ooeta]}=\frac{ ([\sigma_{\ovueps,\ovpeps}] \cdot \ed , \ed )}{[\ooeta]} + O(\e) , 
$$
this in turn gives the implicit interface condition~:
$$
\ovueps \cdot\bfn = -\frac{\e}{[\ooeta]} ([ \sigma_{\ovueps,\ovpeps} ]  \cdot\bfn  ,\bfn) + O(\e^2).
$$
\subsection{The case of an aneurysmal sac }

When $\e$ goes to 0, the limit solution $(\bfuz,\pz)$ is explicit (we set  $\poutu=0$ in \eqref{pois.sac}):
$$
\left\{ 
\begin{aligned}
& \bfuz(x) = \frac{\pin}{2}(1-x_2) x_2 \eu \chiu{\Ou},\quad \forall x \in \Omega \\
& \pz(x) = \pin (1-x_1) \chiu{\Ou} + \pzm \chiu{\Ode},
\end{aligned}
\right.
$$
where $\pzm$ is any real constant. 
Following the same lines as in Theorem \ref{zo_col} one obtains  
\begin{thm}\label{zero.sac} 
For every fixed $\e$, there exists a unique solution $(\ueps,\peps)\in \bH^1(\Oe) \times L^2(\Oe)$ of the problem \eqref{exact}. Moreover, one has
$$
\nrm{\ueps-\bfuz}{\bH^1(\Oe)^2} + \nrm{\peps - \pz}{L^2(\Op)}+ \sqrt{\e} \nrm{\peps - \pz}{L^2(\leps)} \leq k \sqrt{ \e }
$$
where the constant $k$ depends on $\pzm$ but not on $\e$.
\end{thm}

\subsubsection{First order approximation}

Due to the presence of three kind of errors above, we construct a full 
boundary layer approximation $(\cU,\cP)$ exactly as in \eqref{fo.bl.cor}.
One has to make few minor changes in the definition of $(\cW,\cZ)$ that are left to the reader.
The only difference stands in the pressure jump:
$$
[ \pz ] = \pzp (x_1,0)- \pzm,
$$
where $\pzm$ is the constant pressure not yet fixed. The first order macroscopic corrector
$(\bfuu,\pu)$ should satisfy 
\begin{equation}\label{fo.so.macro.sac}
\left\{
\begin{aligned}
-&\Delta \bfuu + \nabla \pu = 0 & \text{ in } \Omega_1\cup \Omega_2,\\ 
& \dive \bfuu = 0 & \text{ in } \Omega_1\cup \Omega_2,\\
&  \bfuu = 0 & \text{ on } \gud\cup\Goutd,\\\ 
& \begin{aligned}
&  \bfuu \cdot \bft = 0, \\
&  \pu = 0, \\ 
\end{aligned}
&\text{ on } \gio, \\
& \bfuu = \duzu \oobeta^\pm + \lrduzu \oov{\bfupsi}^\pm   +  \dpde \oochi &\text{ on } \Gz^\pm .
\end{aligned}
\right.
\end{equation}
As we impose the velocity on every edge of $\Ode$ there is a compatibility 
condition between the Dirichlet data and the divergence free condition reading
$$
\int_{\Ode} \dive \bfuu  \, dx = \int_{\partial \Ode} \bfuu \cdot \bfn \, d\sigma = \int_{\Gz} \bfuu \cdot \bfn  \, dx_1 = 0 ,
$$
and this precisely identifies the pressure $\pzm$ giving
\begin{equation}\label{fix.pzm}
|\Gz| \pzm = \int_{\Gz} \pz^+ (x_1,0) dx_1.
\end{equation}
The first order constants are fixed in the definition of $(\cU,\cP)$.
Even if  $\pzm$ is now well defined, the first and second order pressures 
 $\pu$ and $\pd$  are
again computed in $\Ode$ up to a constant. This is why we
still  need norms on a quotient space $L^2(\Ode)/ \RR$ for the pressure in $\Ode$.
Following the same lines as in the section above but taking 
 into account the pressures in $\Ode$
up to a constant as in the proof of Theorem \ref{zero.sac},  
one proves Theorem \ref{ane.main.thm}.

\section{Numerical validation}\label{num}
We present in this section a numerical validation in the case of a collateral artery, 
as one obtains similar results in the case of an aneurysm we do not display these results. 
We solve numerically problem \eqref{exact} in 2D, for various values of $\e$.
For each $\e$, we confront the corresponding numerical quantities with 
the information provided by the homogenized  first-order explicit approximation : 
velocity profiles, pressure, flow-rate.
Numerical errors estimates are computed with respect to  the different norms evaluated 
above in a theoretical manner. 

We do not include in these sections approximations based on the implicit interface 
conditions presented in \eqref{homo_fo}: this will be done
in a forthcoming work that  investigates  
new theoretical and numerical questions that these conditions pose.

\subsection{Discretizing the rough solution $(\ueps,\peps)$}
 The domain $\Oe$ is discretized for $\e \in ]0,1]$
using a triangulation. To discretize the velocity-pressure variables,
a $(\PP_2,\PP_1)$ finite element basis is chosen. Because of the
presence of microscopic perturbations,  when solving
the  Stokes equations, 
the penalty method gave instabilities. For this reason we opted for
the Uzawa conjugate gradient solver (see p. 178 in \cite{ffm}, and references there in). The code is written in the {\tt freefem++} language\footnote{\url{http://www.freefem.org/ff++}}. On the boundary we impose the following data~: $\pin=2,\poutu=0, \poutd=-1$. 
\begin{figure}[ht!] 
\begin{center}
\includegraphics[width=0.3\textwidth]{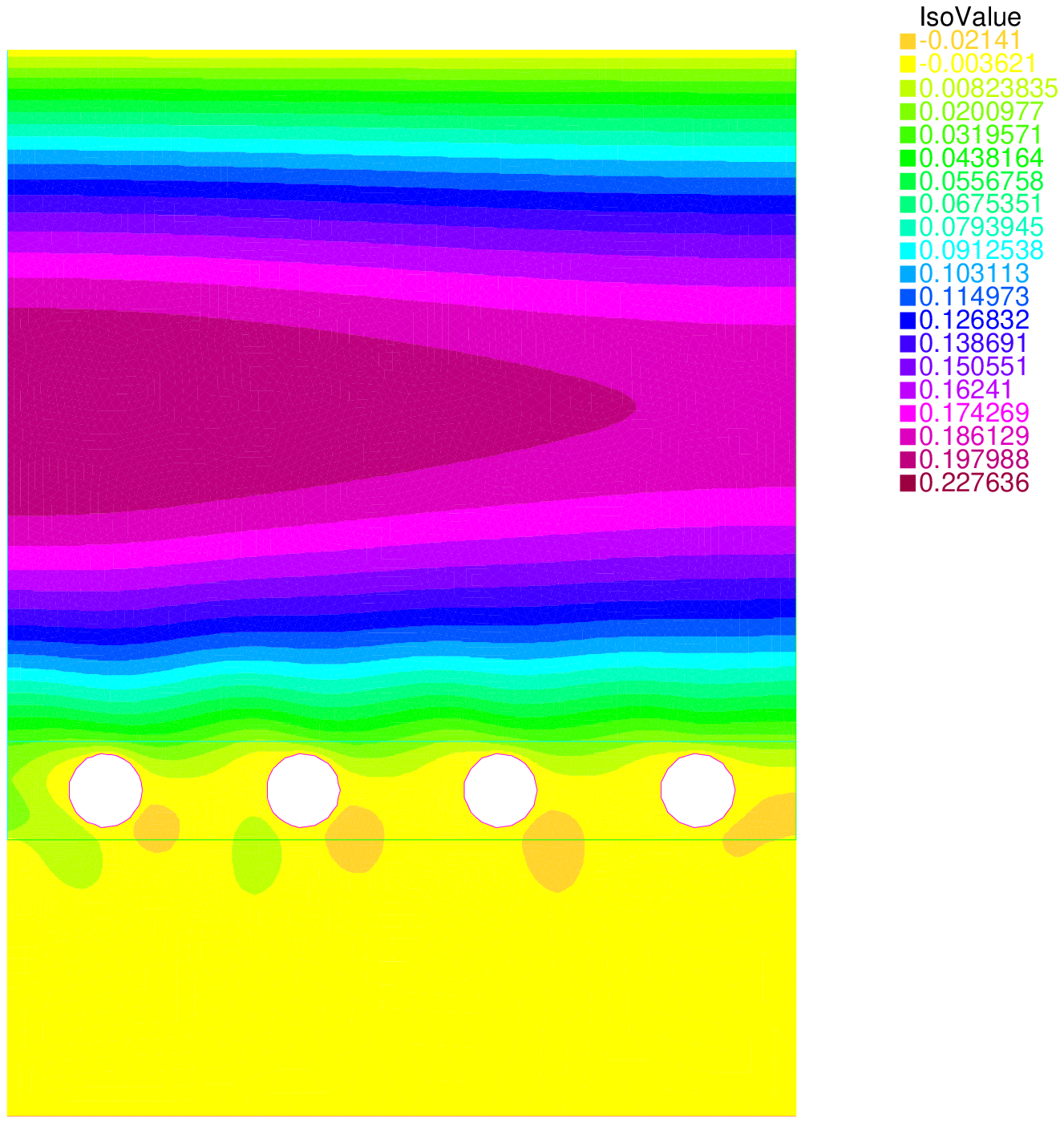}
\includegraphics[width=0.3\textwidth]{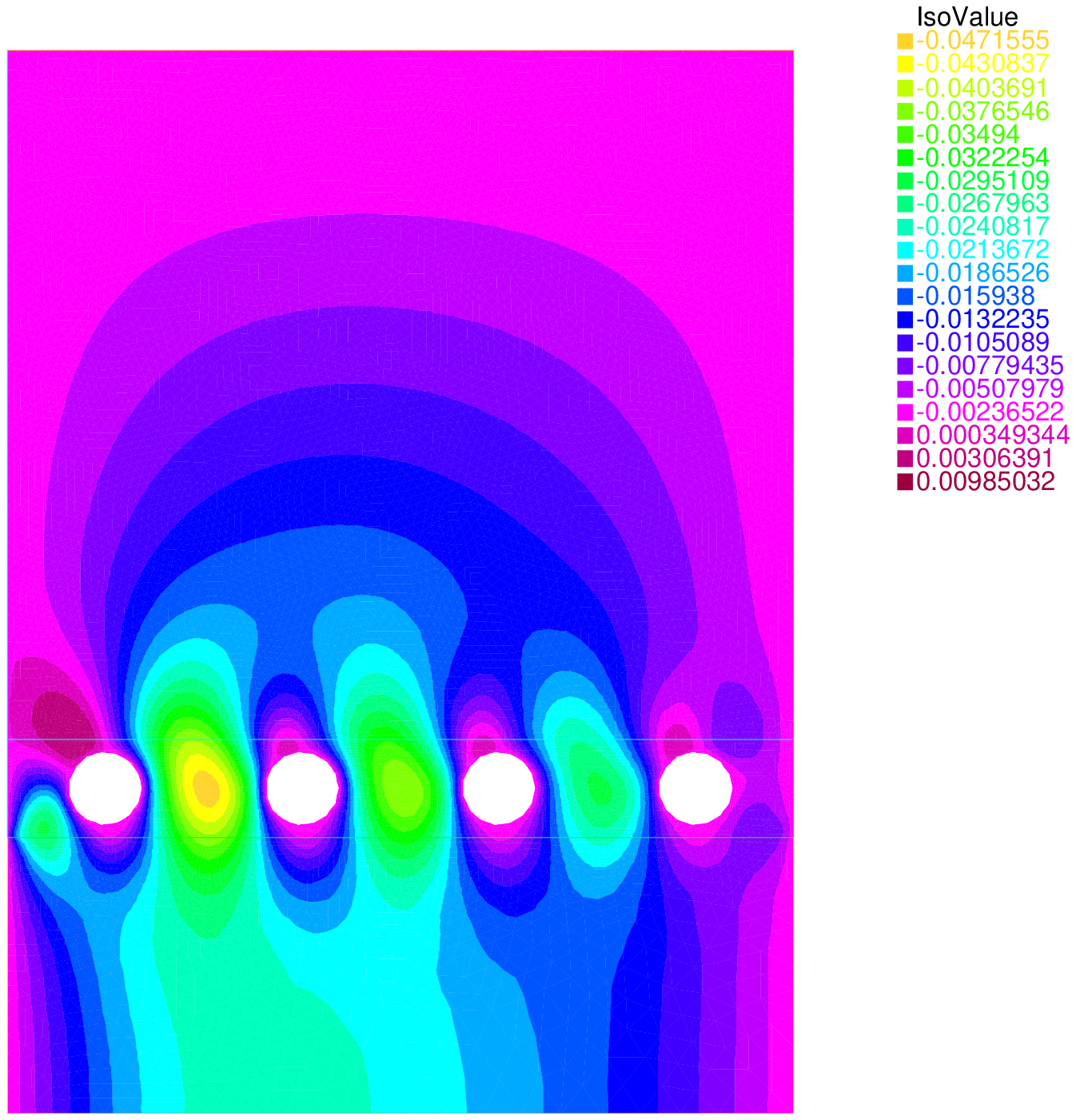}
\includegraphics[width=0.3\textwidth]{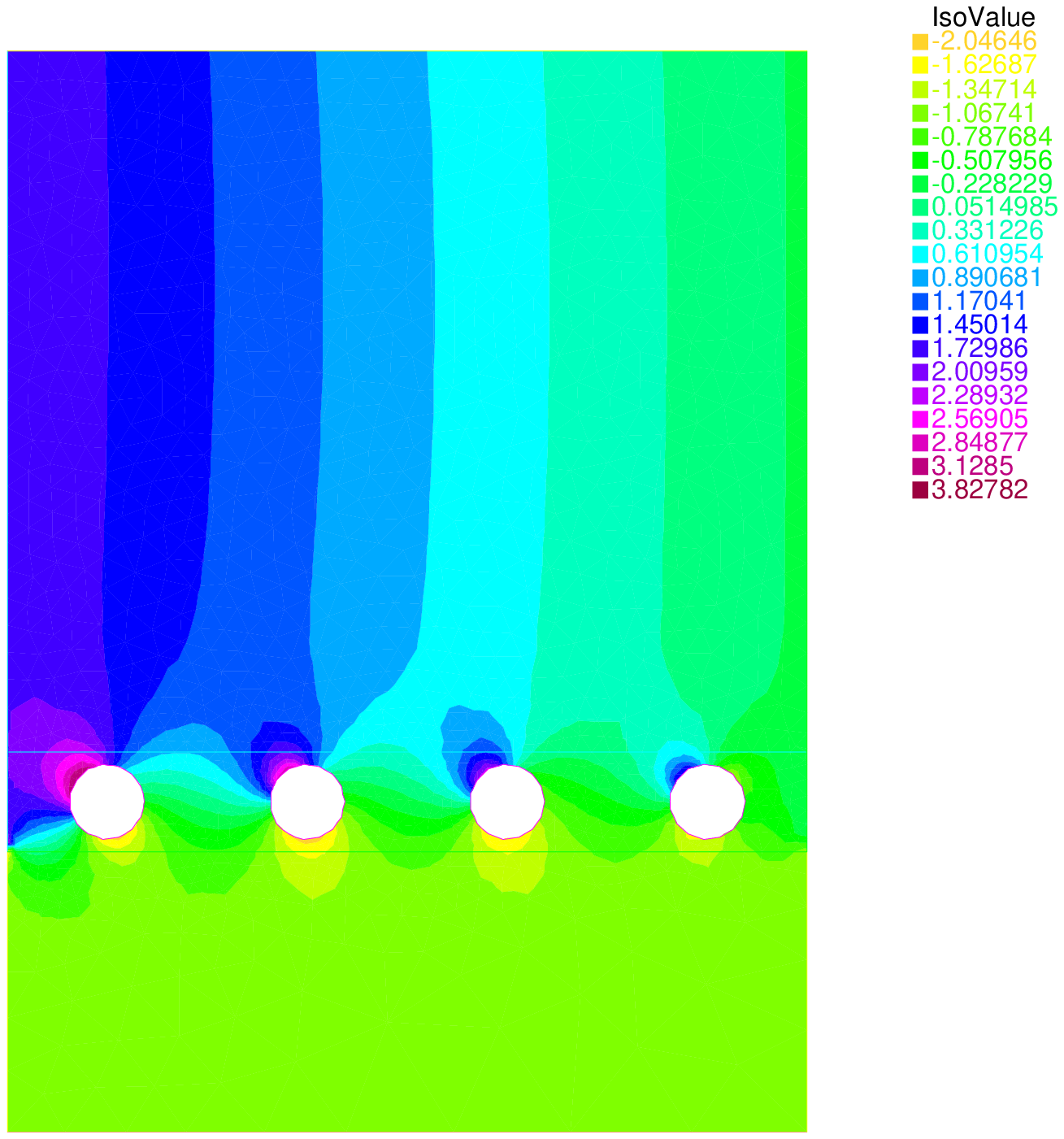}
\end{center}
\caption{Direct computation $u_{\e,1}$ (left) $u_{\e,2}$ (middle) and $\peps$ (right) }\label{fig.direct.sol}
\end{figure}
In order to improve accuracy of the direct simulations we 
use  mesh adaptation iterations as described p. 96-97 in \cite{ffm}: 
using the hessian matrices of components
of $\ueps$, one defines a metric that modifies 
the mesh (see fig. \ref{mesh} (middle) for a final shape of the mesh).

\subsection{The microscopic cell problems}
Using the same numerical tools, we solve the microscopic problems 
\eqref{beta.cell}, \eqref{eq.upsi} and \eqref{chi.cell}. 
These are defined on the infinite perforated strip $Z$: one is forced to 
truncate the  domain and works on $Z_{-L,L}$ with $L>0$ large.
We impose boundary data at the top and the bottom of $Z_{-L,L}$ namely
$$
\beta_2(y) = \Upsilon_2(y) = 0 , \quad \chi_2(y) = -1 \quad \text{ on }\{ y \in ]0,1[ \times \RR \text{ s.t. }y_2 = \pm L\}, 
$$
and we let natural boundary conditions on the other components. When $L$ goes to infinity it is proved in \cite{JaMiNe.01}  that the solutions of the truncated problem defined on $Z_{-L,L}$ converge exponentially with respect to $L$ to the solution of the unbounded problem.
We compute numerical values of $\oobetas_1^\pm$ and $\ooupsis_1^\pm$ and 
the  pressure drop $[\ooeta]$. 
If $\js$ is a sphere of radius $3/16$ in a period of size 1 centered at $(1/2,1/4)$ the numerical computations provide values listed in table \ref{table.homo}.
\begin{table}[ht!] 
\begin{center}
\begin{tabular}{|c|l|c|l|}
\hline
\hline
constants & values   & constants & values \\
\hline
$\oobetas^+_1$ & -0.377928 & $\oobetas^-_1$ & -0.122114  \\
\hline
$\ooupsis^+_1$ &-0.000371269 & $\ooupsis^-_1$ & 0.121744  \\
\hline 
$[\ooeta]$  &27.9435  & & \\
\hline
\end{tabular}
\end{center}
\caption{Homogenized numerical constants}\label{table.homo}
\end{table}
One can notice that contrary to the resistive matrix of \cite{AllairePassoire}
the tangential part of the coefficient are negative. This is due to the fact that the obstacles lie above the interface in the main flow. The horizontal first order slip velocity is thus negative (see below).

\begin{figure}[ht!] 
\begin{center}
\includegraphics[width=0.3\textwidth]{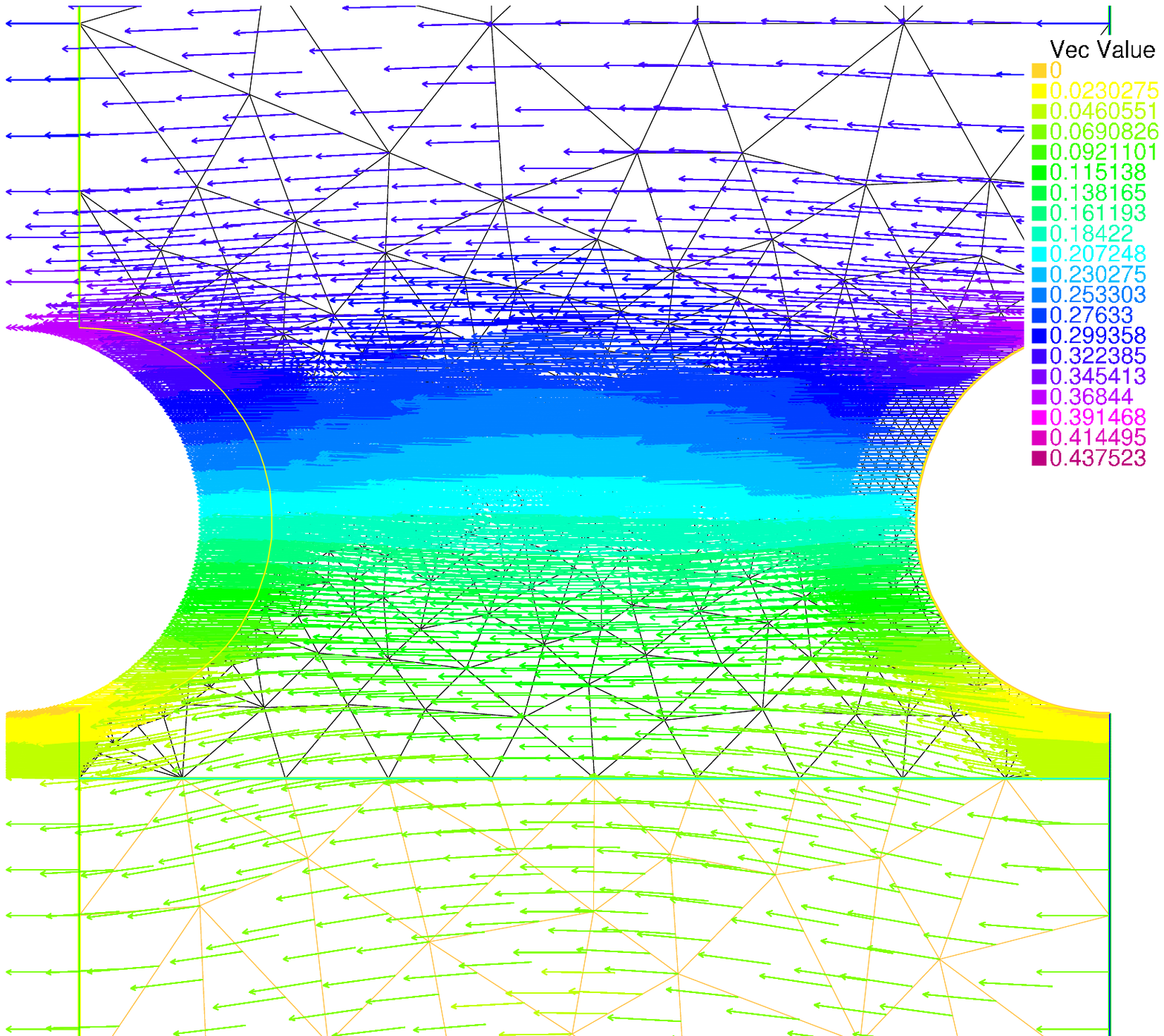}
\includegraphics[width=0.3\textwidth]{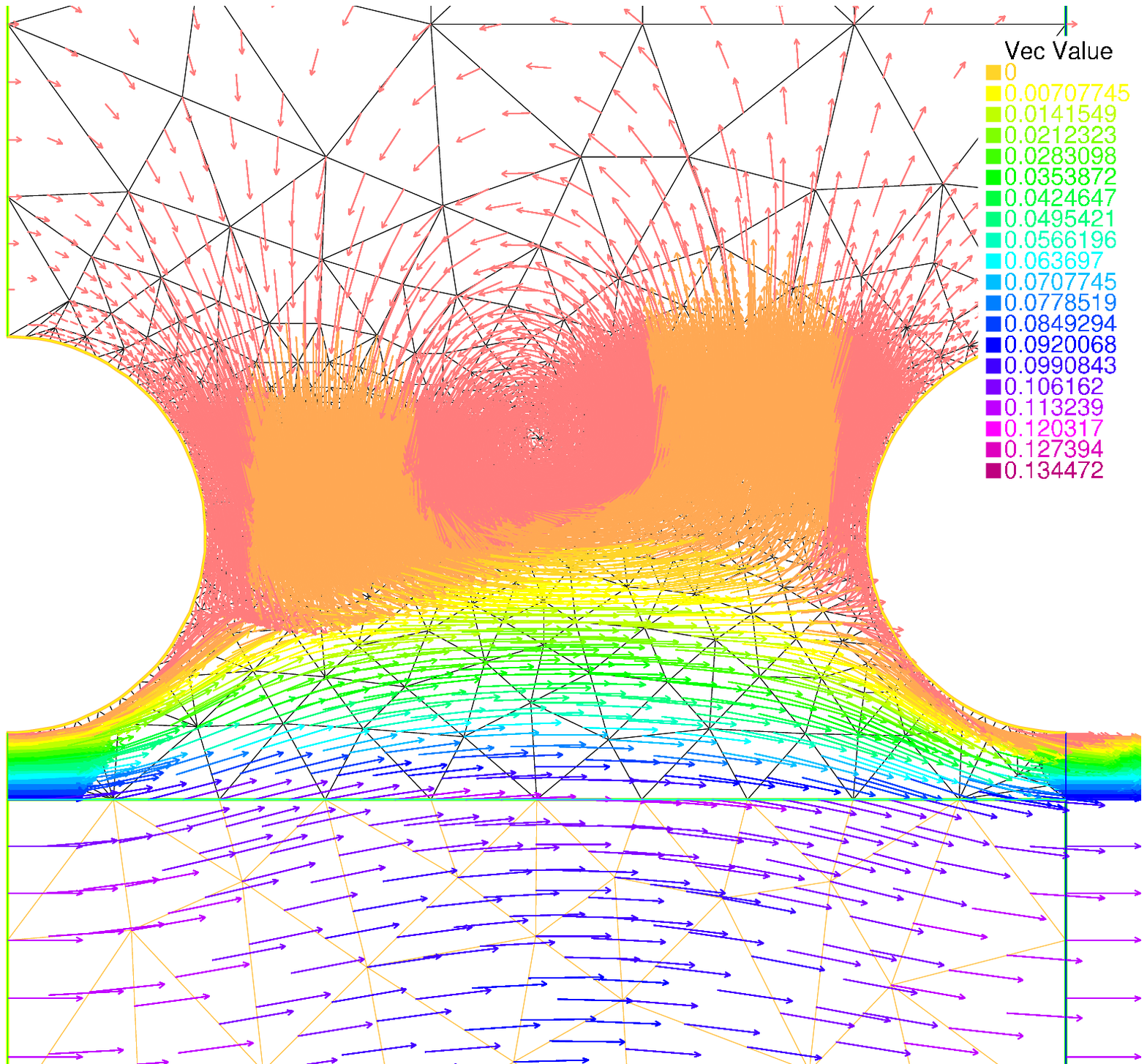}
\includegraphics[width=0.3\textwidth]{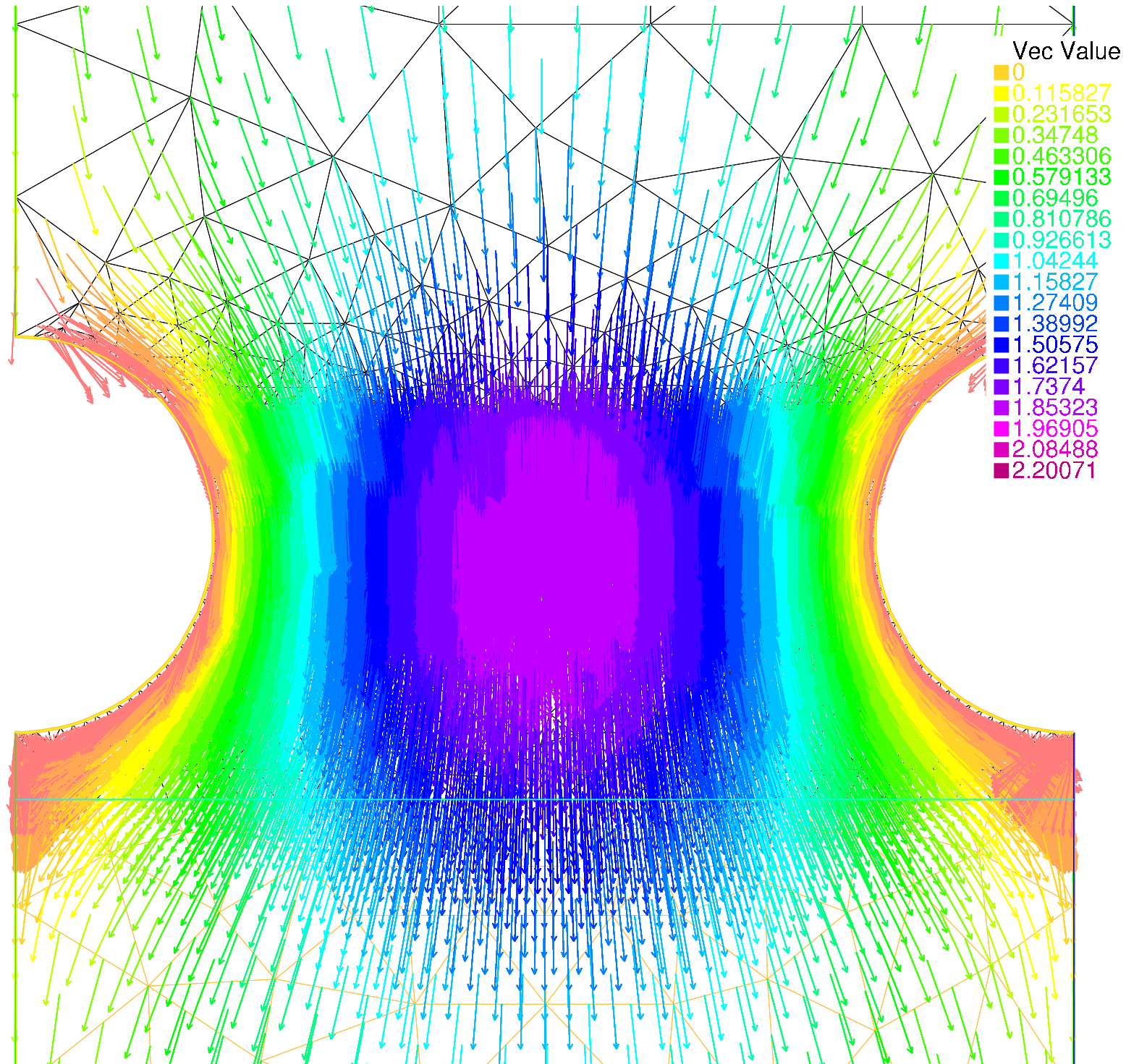}
\end{center}
\caption{Velocity vectors for micrscopic problems: $\bfbeta$ (left), $\bfupsi$(middle) and $\bfchi$ (right)}\label{fig.velo.micro}
\end{figure}

\begin{figure}[ht!] 
\begin{center}
\includegraphics[width=0.3\textwidth]{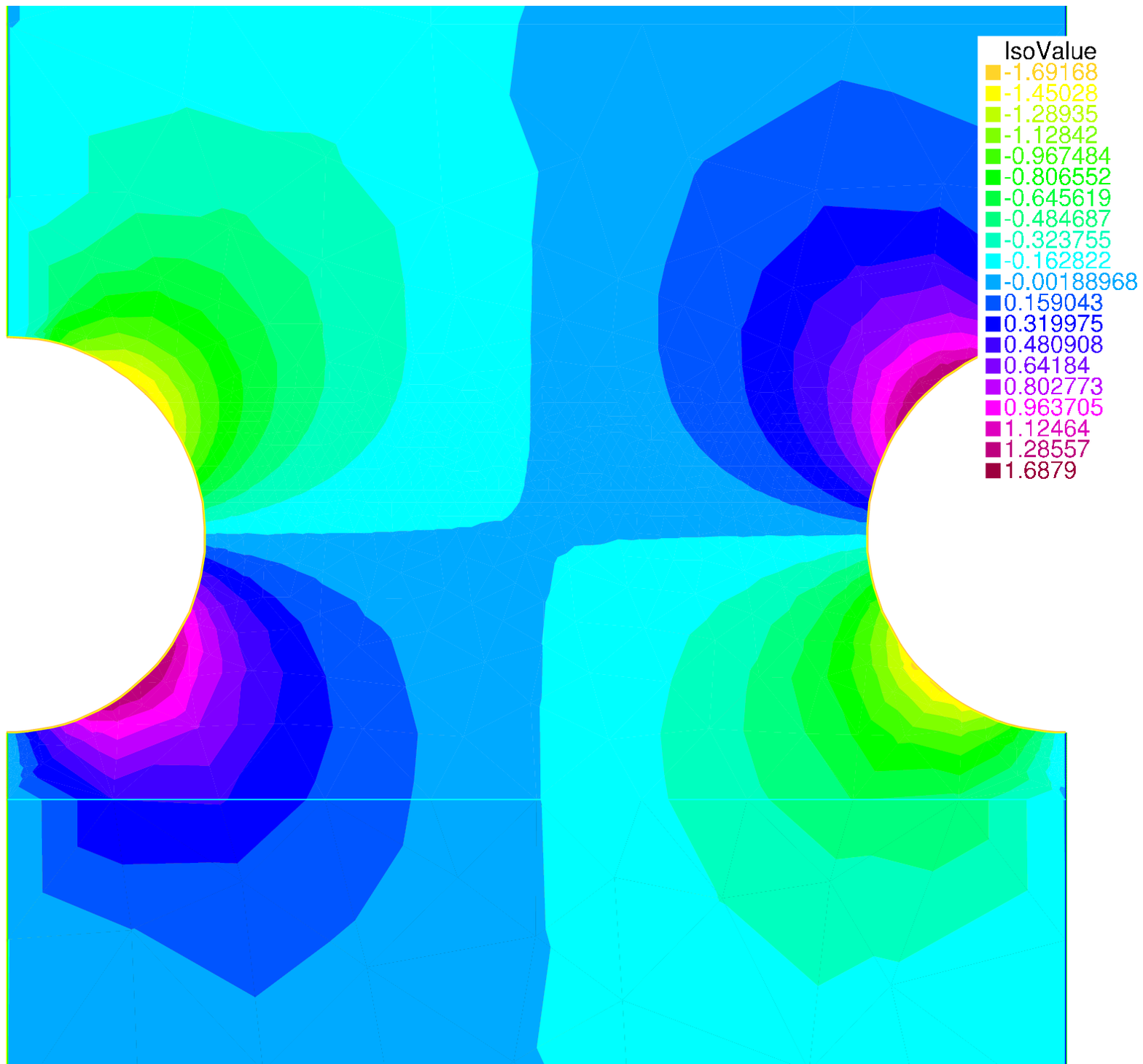}
\includegraphics[width=0.3\textwidth]{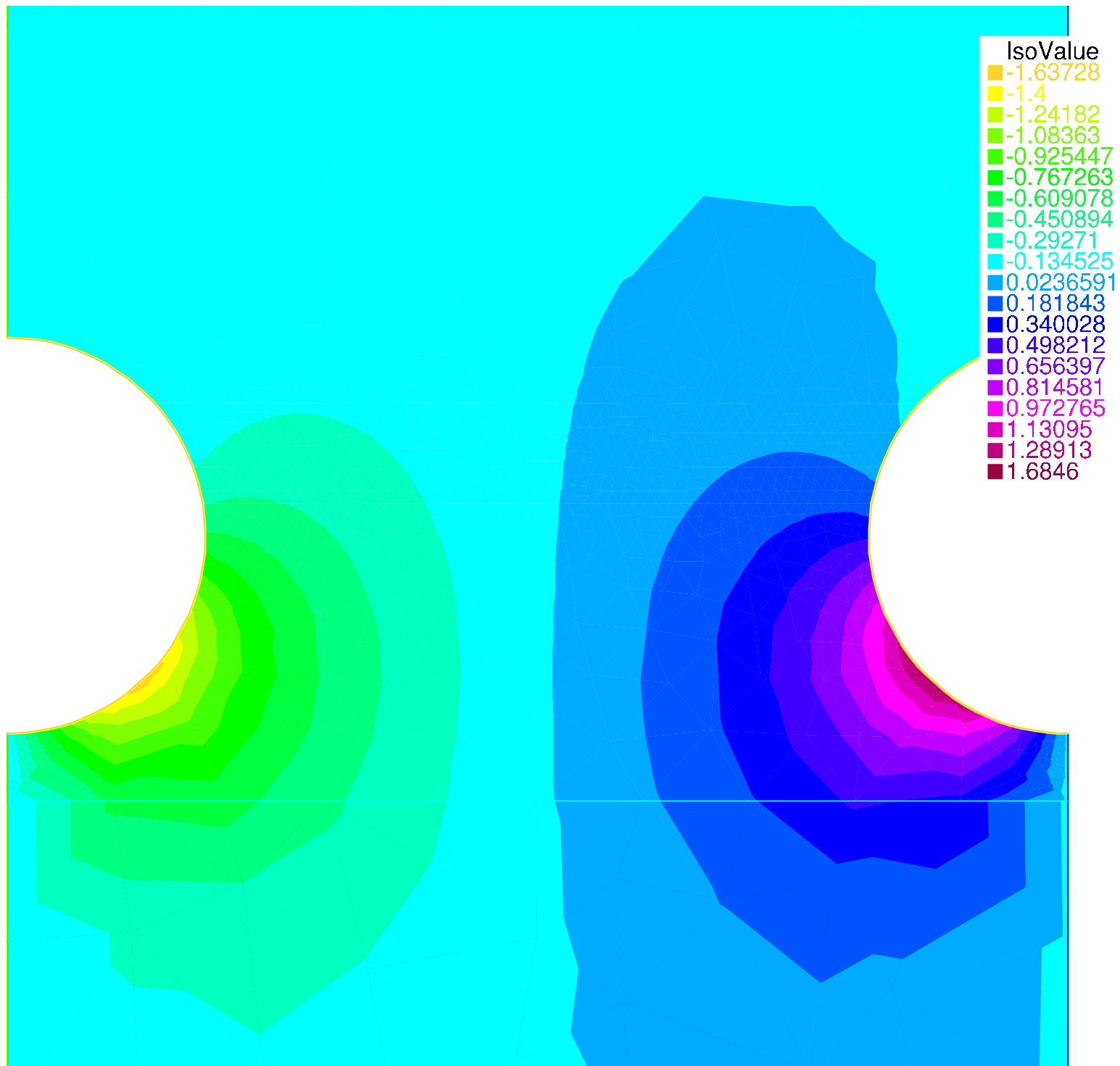}
\includegraphics[width=0.3\textwidth]{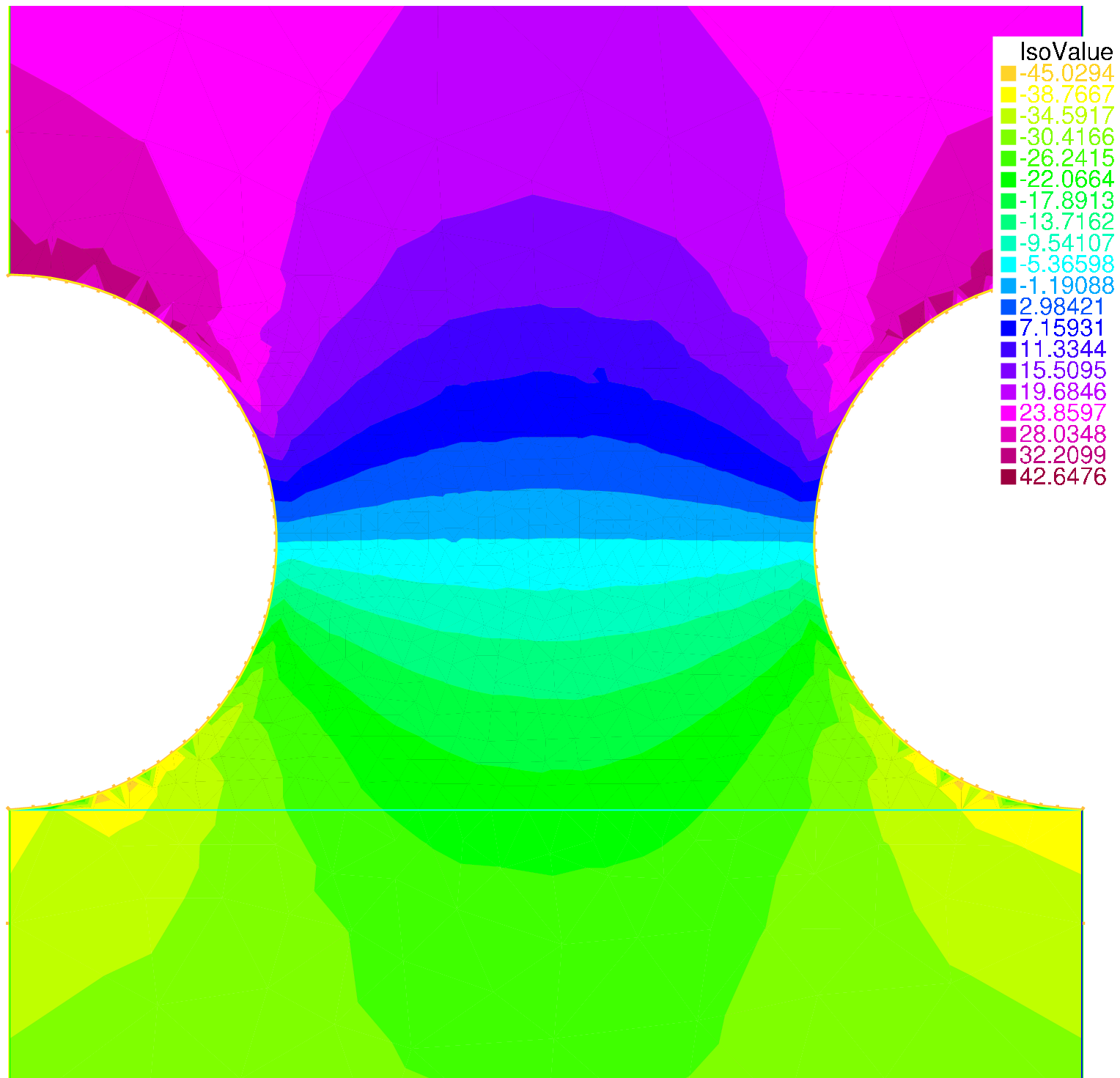}
\end{center}
\caption{Pressure for micrscopic problems: $\pi$ (left), $\mu$(middle) and $\eta$ (right)}\label{fig.press.micro}
\end{figure}

\subsection{Explicit first order problem}
\begin{figure}[ht!] 
\begin{center}
\includegraphics[width=0.3\textwidth]{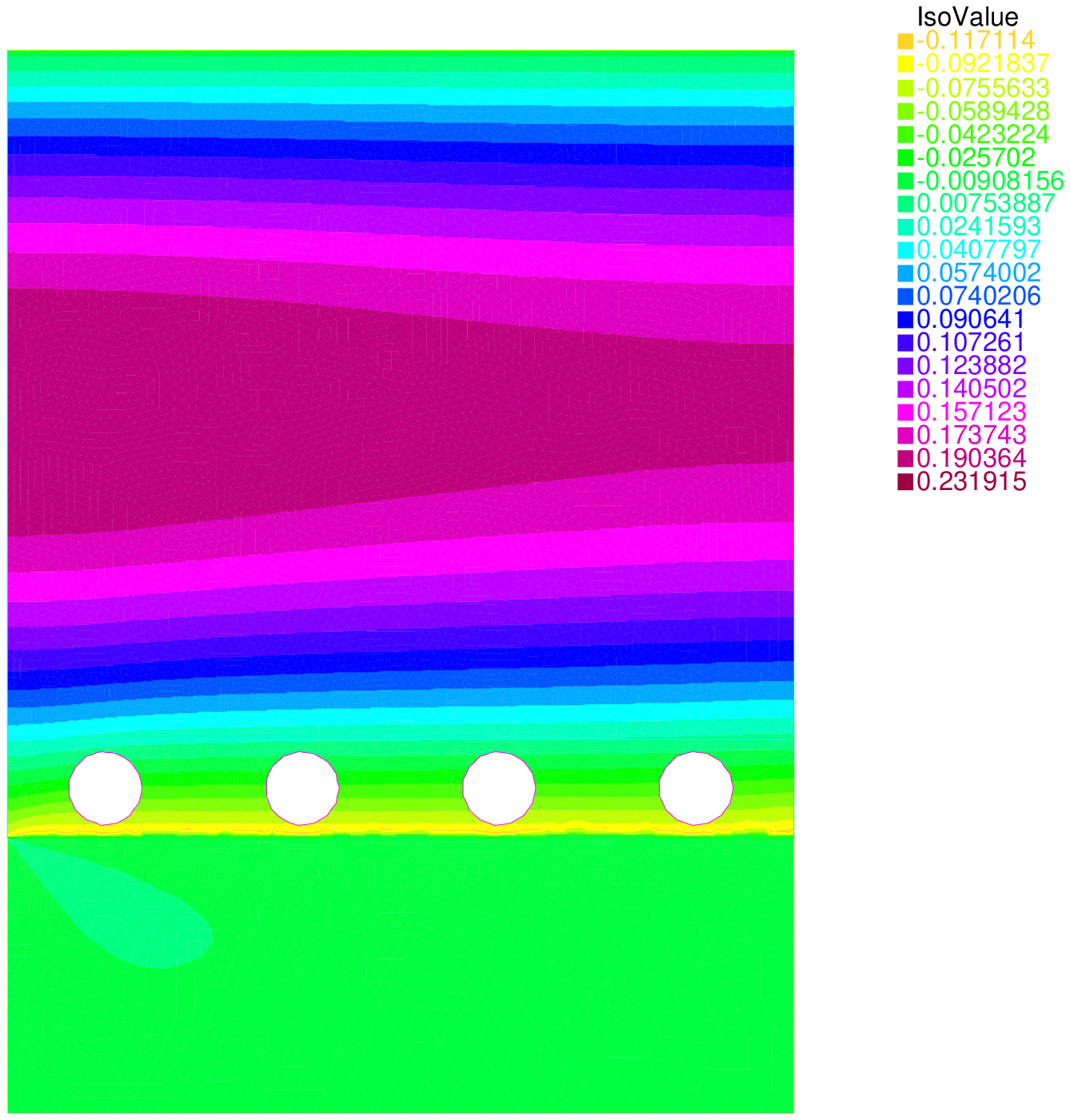}
\includegraphics[width=0.3\textwidth]{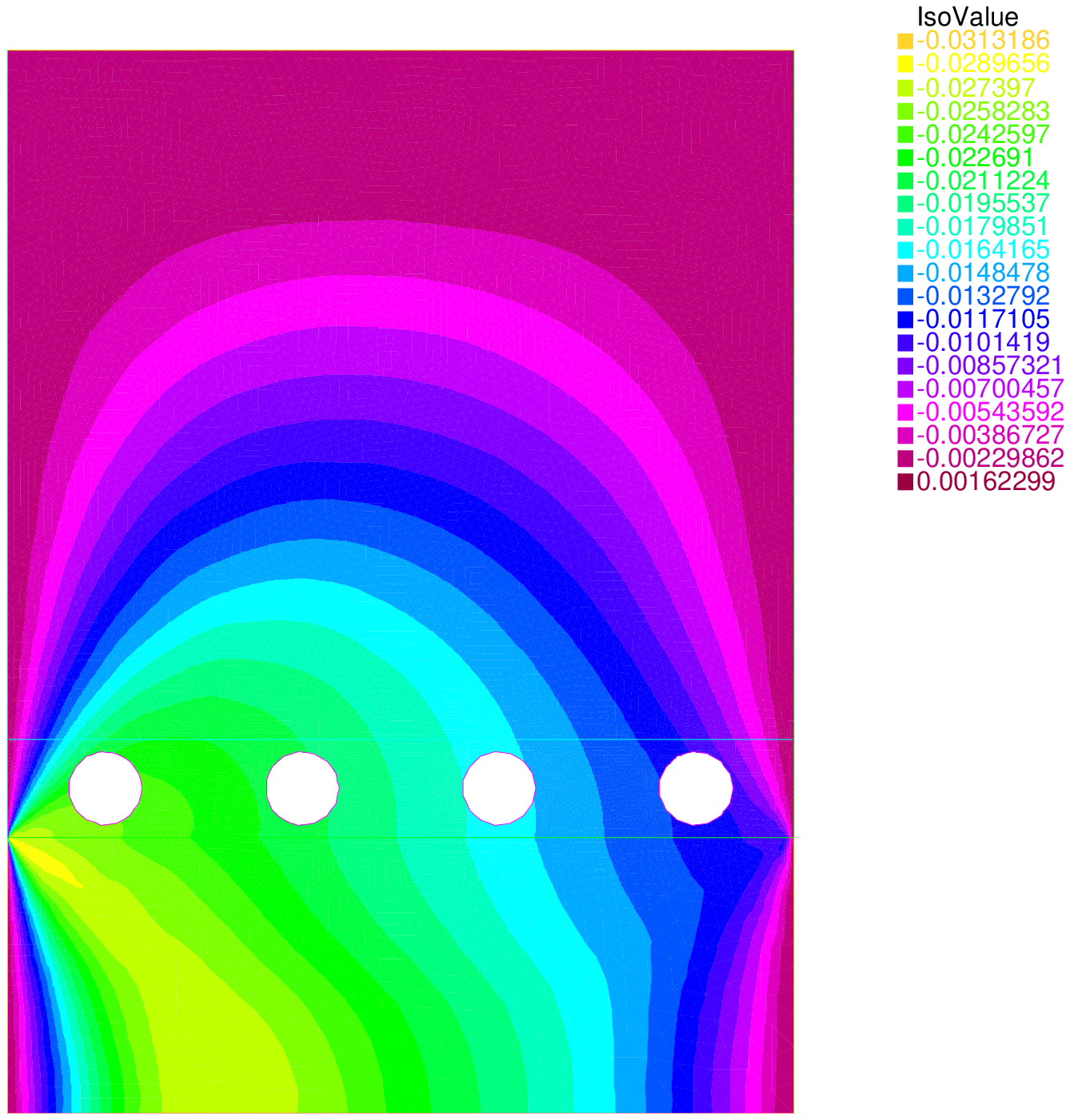}
\includegraphics[width=0.3\textwidth]{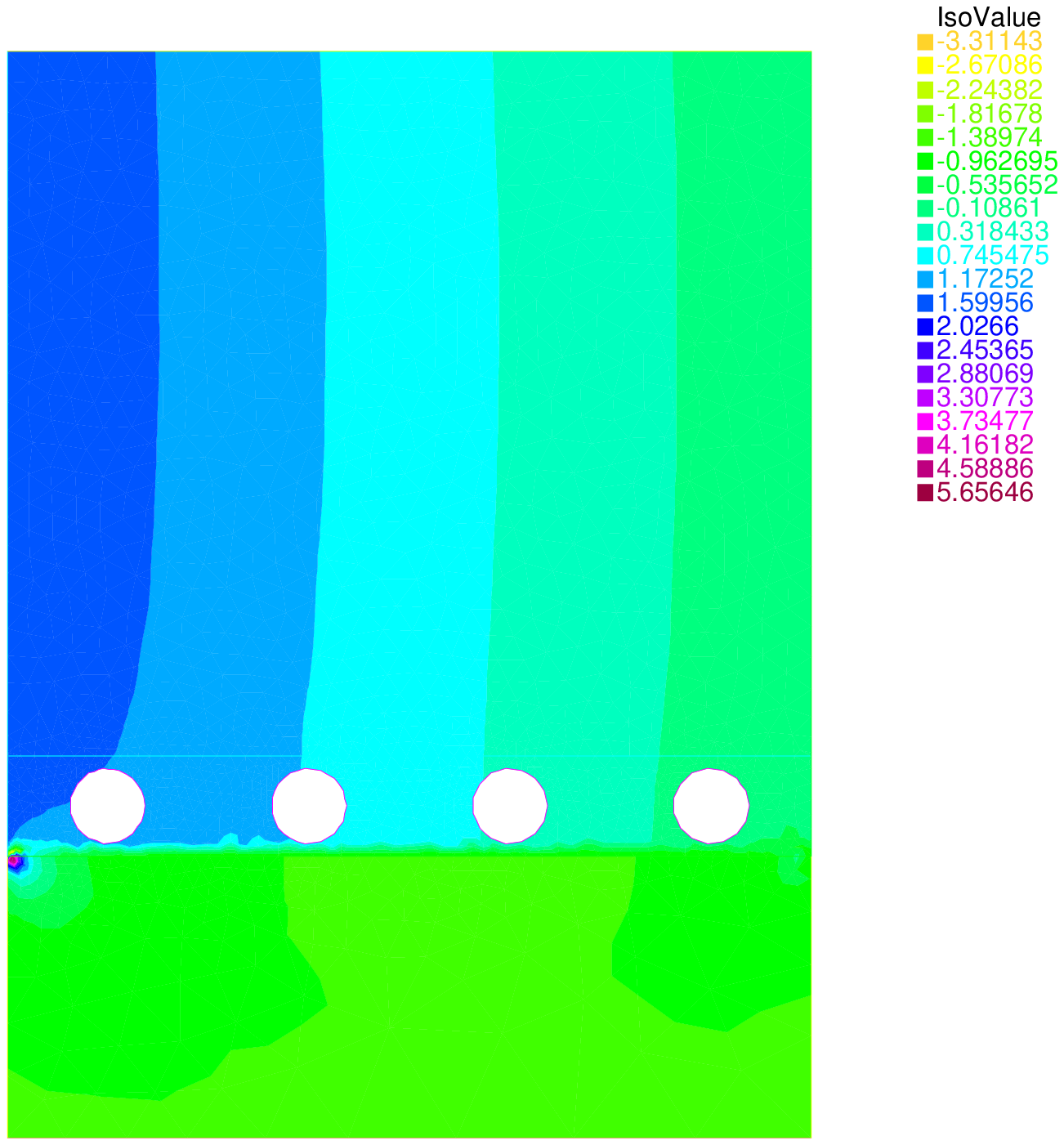}
\end{center}
\caption{
Explicit first order approximation $\ov{u}_{\e,1}$ (left), 
$\ov{u}_{\e,2}$ (middle) 
and $\ovpeps$ (right) }\label{fig.explicit.sol}
\end{figure}

We solve problem \eqref{fo.so.macro} on triangulations of $\Ou$ and $\Ode$. 
Because of the discontinuity of the Dirichlet data at the corners $O$ and $\ov{x}$, 
the solution $(\bfuu,\pu)$ does not belong to $\bH^1(\Ou\cup \Ode)\times L^2(\Ou\cup \Ode)$. 
Indeed a  pressure singularity occurs at $O$ and $\ov{x}$: refining the 
triangulation at the corners one get a point-wise explosion of the pressure near $O$ and $\ov{x}$. We add then the zeroth order explicit poiseuille profile to obtain a numerical approximation of $(\ovueps,\ovpeps)$. For $\e=0.25$, we display in fig. \ref{fig.explicit.sol} velocity components and pressure projected on $\Oe$ in order to be compared to $(\ueps,\peps)$ in the next paragraphs.
Since the pressure is not bounded (the numerical value is very high in a very small neighborhood of $O$ and $\ov{x}$  we display the ``regular part'': we cut-off the pressure function near the corners for visualisation purposes only.
\subsection{Comparisons and error estimates} 
\begin{figure}[ht!] 
\begin{center}
\input{velo_horiz}
\input{velo_vert}
\input{q2_colateral}
\end{center}
\caption{Velocity profiles for $\e=0.125$ close to $\Gz$: horizontal (left) and  vertical (middle). Rough versus homogenized flow-rates (right) wrt $\e$}\label{fig.profile}
\end{figure}
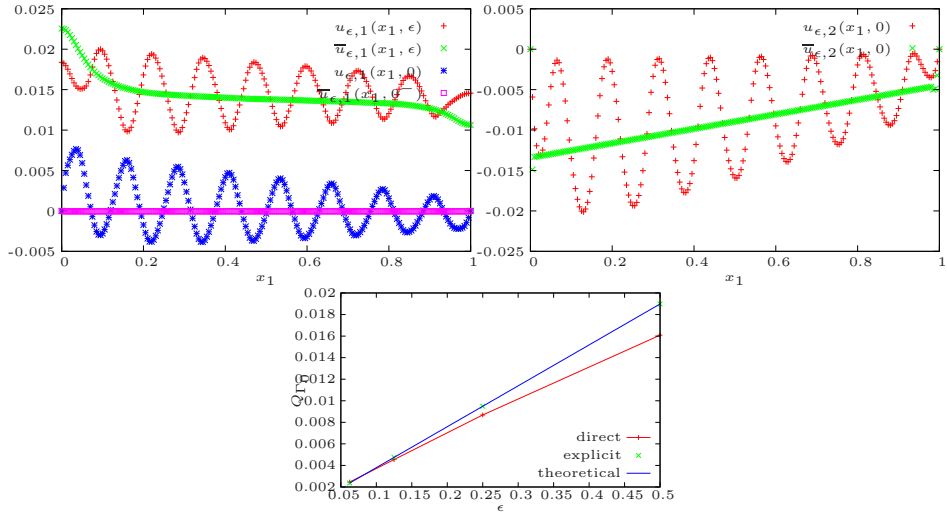

\begin{figure}[ht!] 
\begin{center}
\input{error_l2}
\input{error_hm1}
\end{center}
\caption{Numerical error estimates :  velocity (left) and  pressure (right) wrt $\e$}\label{fig.error}
\end{figure}
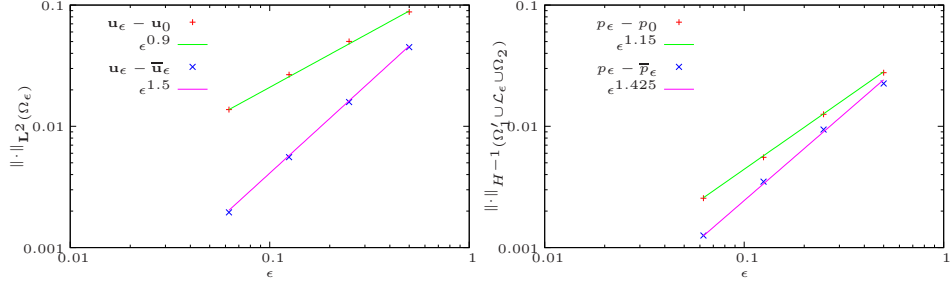

In fig. \ref{fig.profile}, we display the horizontal (left) velocity profile above and below the obstacles for the direct solution $\ovuepsu$ and our approximation $\ovuepsu$. In the middle we show the normal velocity in the same framework. On the right of the same figure, we plot various values of $\e$ on the $x$-axis and the flow-rate through $\Gz$ on the $y$-axis. One observes that the
asymptotic expansion gives the first order approximation of the flow-rate with respect to  $\e$ near $\e=0$ which was expected. One notices also that the actual rough flow-rate behaves as a square-root of $\e$. This seems difficult to prove using averaged interface conditions only \cite{BrMiCras,BrMiQam}. 
In fig. \ref{fig.error},  we plot numerical error estimates for the zero order approximation $(\bfuz,\pz)$ 
and for our explicit first order averaged approximation $(\ovueps,\ovpeps)$ wrt to the direct solution $(\ueps,\peps)$. 
Left we display the $\bL^2(\Oe)$ error for velocity vectors. On the right, we compute the pressure error estimates in the $H^{-1}(\Ou'\cup\leps\cup\Ode)$ norm: for $\pz$ (resp. $\ovpeps$) we solve numerically
$$
\left\{ 
\begin{aligned}
&  - \Delta q = \peps - \pz \, (\text{resp. } \ovpeps), & \text{ in } \Ou'\cup\leps\cup\Ode, \\
& q= 0& \text{ on } \partial \Ou'\cup P \cup \partial \Ode, 
\end{aligned}
\right.
$$
for each $\e$ then we display $\nrm{\nabla q}{L^2( \Ou'\cup\leps\cup\Ode )}$. One recovers theoretical claims 
of Theorems \ref{zo_col} and \ref{main.thm.paper}.  

\begin{figure}[ht!]  
\input{hmesh}
\includegraphics[width=0.2\textwidth]{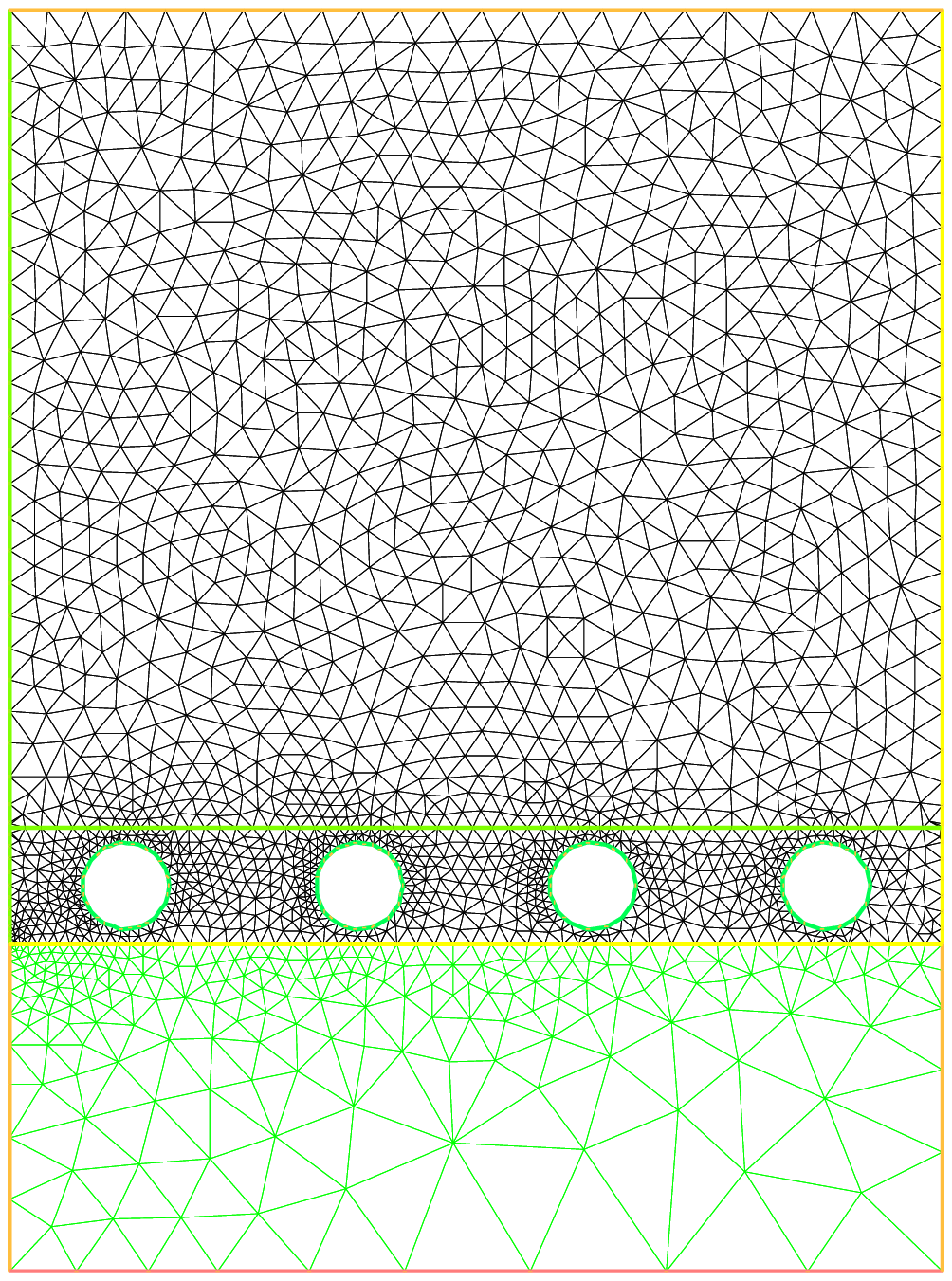}
\includegraphics[width=0.2\textwidth]{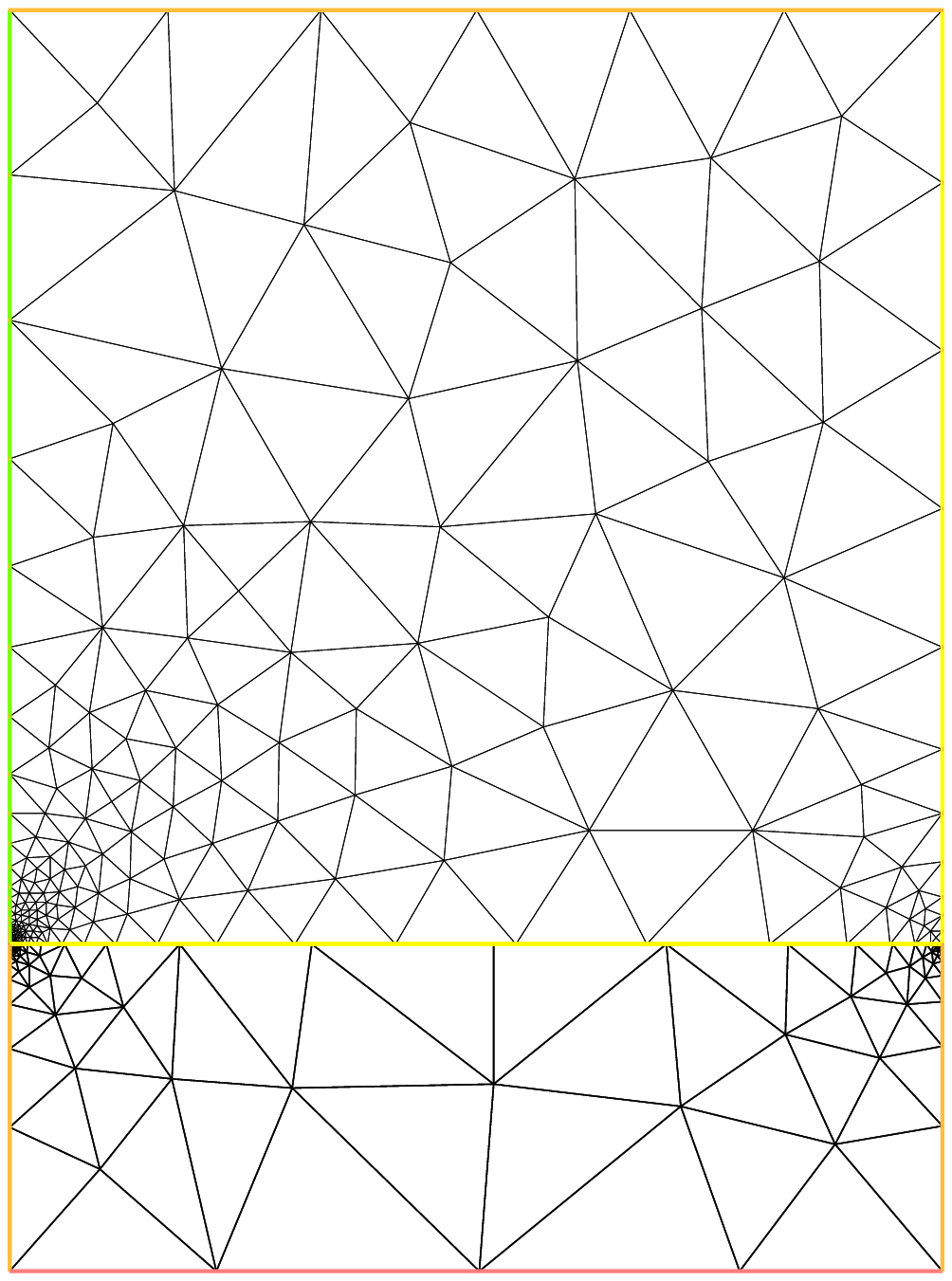}
\caption{The mesh size $h$ for direct computations (left), a direct mesh for $\e=0.25$ (middle) and the mesh used to compute $(\bfuu,\pu)$ (right)}\label{mesh}
\end{figure}

In fig. \ref{mesh} middle and right we display the meshes used for a single value of $\e=0.25$ and for the computations of $(\bfuu,\pu)$. On the left we display  the mesh size $h$ used for the direct simulations with respect to $\e$. 

\bigskip

\newpage 

\appendix
\section{Well-posedness in weighted Sobolev spaces of the vertical correctors}
\label{app.1}

Given the  data $({\bf f},h)$,
we study the problem: find $(\bfw,\theta) \in \bws{1}{2}{\alpha}{\Pi}\times\ws{0}{2}{\alpha}{\Pi}$ solving
\begin{equation}\label{half_space_stk_eq} 
\left\{ 
\begin{aligned} 
& -\Delta \bfw + \nabla \theta = 0 & \text{ in } \Pi, \\ 
& \dive \bfw =0 & \text{ in } \Pi, \\
&  \bfw = {\bf f} &\text{ on } D,\\ 
&  \bfw \cdot \bft = {\bf f} \cdot \bft  , \text{ and }   \theta = h  &\text{ on } N,\\ 
\end{aligned}
\right.
\end{equation}

\begin{remark}\label{rmk.mix.type.bc.hp}
 For the specific type of mixed boundary conditions set on $\partial \Pi$, it is not possible to use neither the Fourier transform in the vertical direction nor the Laplace transform in the horizontal direction in order to derive results established below. To our knowledge there are few results in the literature for unbounded domains with mixed boundary conditions on  noncompact boundaries. 
\end{remark}

\begin{thm}\label{half_space_stk_thm}
If the real $\alpha$ is such that $|\alpha| < 1$ and if
$$
{\bf f}\in \bws{\ud}{2}{\alpha}{D\cup N}, \quad h \in \ws{-\ud }{2}{\alpha}{N},
$$
 there exists a unique solution $(\bfw,\theta) \in \bws{1}{2}{\alpha}{\Pi}\times\ws{0}{2}{\alpha}{\Pi}$
solving problem \eqref{half_space_stk_eq}.
\end{thm}

Before giving the proof of the theorem, we need to prove two intermediate 
propositions. For this sake we define
$$
\begin{aligned}
&\xalpha := \{ \bfv \in \bws{1}{2}{\alpha}{\Pi}\text{ s.t. } \bfv = 0 \text{ on } D,\quad \bfv \cdot \bft = 0 \text{ on } N\} ,\quad Y_\alpha := \ws{0}{2}{\alpha}{\Pi}.
\end{aligned}
$$
First we show that the divergence operator is surjective from $\xalpha$ into $Y_\alpha$
\begin{prop} \label{surje}
For any given function $q \in Y_\alpha$ there exists a vector function $\bfv\in \xalpha$ 
such that
$$
\dive \bfv = q, \quad  \text{and}  \quad  \snrm{\bfv}{\xalpha}\leq k(\Pi,\alpha) \nrm{ q }{Y_\alpha}
$$
where the constant $k$ depends only on the geometry of the domain,
and on $\alpha$.
\end{prop}
\begin{proof}
We define a sequence of annular domains covering $\Pi$
$$
\begin{aligned}
& C_n := \{ y \in \Pi \text{ s.t. if } y=(r,\tteta) \text{ then }\quad r\in ]2^{n-1},2^n[\} ,\quad n \geq 1 ,\quad 
& C_0 := B(0,1) \cap \Pi .
\end{aligned}
$$
We decompose $q$ as $q=\sum_{n=0}^\infty q_n$, with $q_n := q \chiu{C_n}$. On each $C_n$ 
we solve the problem: find $\bfv_n  \in \bxan $ s.t. $\dive \bfv_n = q_n$ and $\snrm{\bfv_n}{\xalpha}\leq k(C_n,\alpha) \snrm{q_n}{Y_\alpha}$, where
$$
\begin{aligned}
\bxan:= \left \{ \bfv \in \xalpha \text{ s.t. } \right.  & \bfv = 0 \text{ on } \{ |y|=2^{n-1} \} \cup \{ |y|=2^{n} \} \cup (\ov{C}_n\cap D),  \\
\text{ and } & \left. \bfv \cdot \bft = 0 \text{ on } (\ov{C}_n \cap N) \right\} .
\end{aligned}
$$
But to solve the latter equation 
in a weak sense means 
$$
\int_{C_n} \dive \bfv_n \cdot \omega  \; rdr d\tteta  = \int_{C_n} q_n \cdot  \omega \; rdr d\tteta ,\quad \forall \omega \in \ws{0}{2}{-\alpha}{C_n},
$$
making the change of variables: $(\ti{r}=r / 2^{n-1},\tteta)$ and setting
$$
\tbfv_n(\ti{r},\tteta) := \bfv(2^{n-1} \ti{r} , \tteta),\quad \ti{q}_n(\ti{r},\tteta) := q(2^{n-1} \ti{r} , \tteta), \quad \ti{\omega} (\ti{r},\tteta) := \omega (2^{n-1} \ti{r} , \tteta),
$$
the problem becomes: find $\tbfv \in \bm{X}_{0,1}$ defined on $C_1$ s.t. 
$$
\int_{C_1} \ti{\dive} \tbfv_n \cdot \ti{\omega} \ti{r} d \ti{r} d \tteta 
= 2^{n-1} \int_{C_1}   \ti{q} \cdot  \ti{\omega} \; \ti{r}d\ti{r} d\tteta ,\quad \forall \omega \in L^2(C_1),
$$
the test space is defined on a compact fixed domain $C_1$,  weighted Sobolev
spaces coincide with the classical ones as soon as the weight is 
strictly positive and bounded. In this framework the operator 
$\ti{\dive}:\bm{X}_{0,1}\to Y_{0,1}$ is surjective thanks to  Lemma 4.9 p. 181 in \cite{ErGu.04.book}.
Thus there exists $\tbfv_n\in \bm{X}_{0,1}$ s.t. $\ti{\dive} \tbfv_n = \ti{q}_n 2^{n-1}$. Note that
there is no need of a compatibility condition on the integral of $\ti{q}_n$ 
as in Lemma 3.1 chap. III in \cite{Galdi} because $\tbfv_n \cdot \bfn \neq 0$ on $N\cap \ov{C}_1$. Moreover one has that
$$
\snrm{\tbfv_n}{\bH^1(C_1)} \leq k(C_1,0) \nrm{2^{n-1} \ti{q}_n}{L^2(C_1)} ,
$$
where $k$ depends only on the geometry of $C_1$ and is thus independent on $n$. Turning back 
to the original variables $(r,\tteta)$ one has then that $\dive \bfv_n = q_n$ and 
$$
\int_{C_n} \left| \nabla \bfv_n \right|^2 r \, dr d\tteta \leq k(C_1,0) \int_{C_n} |q_n|^2 r \, dr d\tteta .
$$
In order to recover the global weighted norm of $q$  in $\ws{0}{2}{\alpha}{\Pi}$, we multiply 
the inequality by $2^{2\alpha (n+1)}$ on both sides; we use that for $r \in [ 2^{n-1},2^n ]$,  
$\rho:= (1+r^2)^\ud$ can be estimated as 
$2^{2\alpha(n-1)} \leq \rho^{2\alpha} \leq 2^{2\alpha(n+1)}$ giving finally 
$$
\begin{aligned}
\int_{C_n} \left| \nabla \bfv_n \right|^2 \rho^{2\alpha} dy & \leq \int_{C_n} \left| \nabla \bfv_n \right|^2 2^{2\alpha(n+1)} dy \\
& \leq k(C_1,0) 2^{2\alpha(n+1)} \int_{C_n} q_n^2 dy \leq 2^{4\alpha} k(C_1,0) \int_{C_n} q_n^2 \rho^{2\alpha} dy .
\end{aligned}
$$
One defines  $\bfv:= \sum_n \bfv_n \chiu{C_n}$, because of the boundary conditions imposed on each of the $C_n$,  
$\bfv$ is  continuous on $\Pi$ and thus belongs to $\bws{1}{2}{\alpha}{\Pi}$. This gives the result.
\end{proof}

We lift problem \eqref{half_space_stk_eq} by subtracting to $\bfw$ a function $\cR(\bfw)$ satisfying:
$$
\cR(\bfw) \in \bws{1}{2}{\alpha}{\Pi}, \quad \cR(\bfw)={\bf f} \text{ on } D\cup N. 
$$
Such a lift exists (cf p. 249 \cite{Ha.71} for an explicit  form of $\cR(\bfw)$). 
We correct the divergence of $\cR(\bfw)$ by setting~:
$$
\cS (\bfw) \in \xalpha \text{  s.t. } \dive \cS(\bfw)= -\dive(\cR(\bfw)) \text{ and } \nrm{\cS(\bfw)}{\bws{1}{2}{\alpha}{\Pi}} \leq k  \nrm{\cR(\bfw)}{\bws{1}{2}{\alpha}{\Pi}},
$$
which is possible thanks to Proposition \ref{surje}.
The new variables 
$(\ti{\bfw}:=\bfw-\cR(\bfw)-\cS(\bfw),\theta)$ solve the homogeneous problem:
\begin{equation}\label{half_space_stk_eq_dirichlet_homo} 
\left\{ 
\begin{aligned} 
& -\Delta \ti{\bfw} + \nabla \theta = \Delta \cR(\bfw) + \Delta \cS(\bfw)  & \text{ in } \Pi, \\ 
& \dive \ti{\bfw} = 0 & \text{ in } \Pi, \\
&  \ti{\bfw} = 0 &\text{ on } D,\\ 
&  \ti{\bfw} \cdot \bft = 0 , \text{ and }   \theta = h & \text{ on } N.\\ 
\end{aligned}
\right.
\end{equation}
No we claim that $(\mF,\mG):=(\rho^\alpha \ti{\bfw}, \rho^\alpha \theta)$ solve in an equivalent way
the problem: find $(\mF,\mG)$ in $\bws{1}{2}{0}{\Pi} \times \ws{0}{2}{0}{\Pi}$ s.t.
\begin{equation}\label{rescaled} 
\left\{ 
\begin{aligned} 
& \begin{aligned}
&\cA \, \mF + \cB^T \mG = \rho^\alpha ( \Delta  \cR(\bfw) + \Delta  \cS(\bfw) ) \\ 
& \cB \, \mF  \hspace{1.2cm} = 0 \\
\end{aligned} & \text{ in } \Pi, \\
& \mF = 0 &\text{ on } D \\
& \left. \begin{aligned}
\mF \cdot \bft = 0 \\
\mG = \rho^\alpha h
\end{aligned} \right\} &\text{ on } N
\end{aligned}
\right.
\end{equation}
where
\begin{equation}\label{ops}
\begin{aligned}
&  \cA \mF  := - \Delta  \mF - 2 \ra \nabla \mF \cdot \nabla \frac{1}{\ra} -  \ra \Delta \frac{1}{\ra}\mF, \quad  \text{ and }
& \cB \mF := \dive \mF + \ra \nabla \left(\frac{1}{\ra}\right)\cdot \mF. \\
\end{aligned}
\end{equation}
Indeed, $\rho \in C^\infty(\Pi)$  thus  if $(\tw,\theta)$ solves \eqref{half_space_stk_eq_dirichlet_homo} in the distributionnal sense, then equivalently by its definition the pair $(\mF,\mG)$ solves \eqref{rescaled} also  in the distributionnal sense. 
Uniqueness is insured thanks to the onto mapping between  $\bws{1}{2}{\alpha}{\Pi} \times \ws{0}{2}{\alpha}{\Pi}$ and 
$\bws{1}{2}{0}{\Pi} \times \ws{0}{2}{0}{\Pi}$ (cf Theorem I.3 p. 243 in \cite{Ha.71}): if $(\tw,\theta)$ is a unique solution of \eqref{half_space_stk_eq_dirichlet_homo} then so is  $(\mF,\mG)$ for system \eqref{rescaled} and {\em vice versa}. The boundary conditions match between both problems by  similar onto trace mappings.
Note that the rhs in \eqref{rescaled} belongs to $\bws{-1}{2}{0}{\Pi}$ and the
boundary data to $\ws{-\ud}{2}{0}{N}$. 
We associate to \eqref{rescaled} the corresponding variational setting, namely we define: 
\begin{itemize}
\item   the velocity/pressure test space is $\xz\times Y_0$, 
\item the bi-continuous (resp continuous) forms $\aalpha,\balpha$ (resp. $\lalpha$) read
$$
\left\{ 
\begin{aligned}
&\aalpha ( \mF, \mV ) = ( \nabla \mF, \nabla \mV ) - 2 
\left( \rho^\alpha \nabla \mF  \nabla \left(
                             \frac{1}{\ra}
                      \right) ,\mV 
\right) 
- \left( 
        \rho^\alpha \Delta \left(
                      \frac{1}{\ra}
                \right)   \mF ,\mV 
\right)  ,  \quad   \forall \mF , \mV \in \xz,   \\
& \balpha (  \mF, \mq ) = -\left( \dive  \mF + \ra\nabla\left(\frac{1}{\ra}\right)\cdot \mF   , \mq \right), \quad  \forall \mF \in \xz , \forall \mq \in Y_0, \\
& \lalpha (  \mV ) =  ( \rho^\alpha ( \Delta  \cR(\bfw) + \Delta  \cS(\bfw) ) , \mV )_{ \Pi} - ( \ra h , \mV \cdot \bfn )_{N}, \quad  \forall \mV \in \xz. \\
\end{aligned}
\right.
$$
\item the variational problem:
the problem \eqref{rescaled} can then be restated in an equivalent way: find $(\mF,\mG) \in  \xz\times Y_0 $ s.t.
$$
\left\{ 
\begin{aligned}
  &\aalpha ( \mF , \mV ) &+ \balpha(\mV,\mG) &= \lalpha (\mV) & \forall \mV \in \xz,   \\ 
&\balpha(\mF,\mq) &                     &= 0 &  \forall \mq \in Y_0.
\end{aligned}
\right.
$$
\end{itemize}
We denote by $\tcA :\xz \to  \xz'$ the operator s.t.
$$
\aalpha(\mF,\mV)=<\tcA \mF,\mV>_{ \xz', \xz },\quad  \forall  \mV\in \xz. 
$$
The well-posedness of problem \eqref{rescaled} is equivalent
to  two conditions (Theorem A.56 p. 474 \cite{ErGu.04.book}):
\begin{enumerate}
\item[$(i)$] $\ti{P} \tcA : \ker(\cB) \to \ker(\cB)'$ is an isomorphism 
\item[$(ii)$] $\cB:\xz\to Y_0$ is surjective
\end{enumerate}
where $\ti{P}$ is the restriction of $\tcA$ to the kernel of $\cB$. Here we prove that these conditions are actually 
fulfilled.

\begin{proposition}\label{well_posed}
If $|\alpha|<1$ then  $\tcA$ satisfies condition $(i)$, whereas $\cB$ satisfies $(ii)$ without any restrictions on $\alpha$
\end{proposition}

\begin{proof}
We prove at first that condition $(i)$ is satisfied by showing that $\aalpha$ is a coercive bi-linear form. 
For every vector $\mF$ in $\xz$, one has after integration  by parts of the second term in the definition of $\aalpha$:
$$
\begin{aligned}
\aalpha(\mF,\mF) & =   
\left(\alpha \frac{(y\cdot \bfn)}{\rho^2} \mF , \mF 
\right)_{\partial \Pi} 
+ \snrm{\nabla \mF}{L^2(\Pi)}^2 \\
& \quad  
+ \int_{\Pi} \alpha \left( - \dive  \left(\frac{y}{\rho^2}\right) 
+ \left( \frac{2}{\rho^2}  - \frac{(\alpha+2)|y|^2}{\rho^4} \right) \right) \mF^2 dy 
 = \snrm{\nabla \mF}{L^2(\Pi)}^2 - \alpha^2 \int_\Pi  \frac{|y|^2}{\rho^4} \mF^2 dy,
\end{aligned}
$$
note that the boundary term on the first line above vanishes on $N$ because  $(y\cdot \bfn) = 0$, though $\mF_1\neq 0$ on this part of the boundary. Note also that this integration by part is justified for functions in $\xz$.
We use optimal Poincar{\'e}-Wirtinger estimates already presented in the proof of Theorem 5.3 p. 20 in \cite{MiVws}:
$$
 \int_\Pi \left| \frac{\mF  |y| }{\rho^2} \right|^2 dy \leq  \snrm{\frac{\mF}{\rho}}{\bL^2(\Pi)}^2 \leq \snrm{\nabla\mF}{\bL^2(\Pi)}^2 .
$$
Note that these Poincar{\'e}-Wirtinger estimates are possible because of the homogeneous Dirichlet conditions on $D$: they give stronger weights than the corresponding logarithmic weighted Hardy estimates available in the whole $\RR^2$ \cite{AmGiGiI.94}.
Finally one has
$$
(\tcA \mF, \mF)_\Pi \geq (1- \alpha^2) \snrm{\nabla \mF}{L^2(\Pi)^4},
$$
which implies coercivity of the operator if $|\alpha|<1$. 
Note that this result (also valid in the scalar case) improves Lemma 4.3 in 
\cite{BrBoMi}. This is essentially due to the integration by parts performed on
the term $(\nabla \mF y / \rho^2,\mF )_\Pi$  which avoids estimating this
term separately from  the others. 

We focus on the condition $(ii)$. For all $q \in Y_0$ we look for $\mF \in\xz$ s.t.
$$
\cB \mF = q, \text{ and } \nrm{ \mF }{\xz} \leq k \nrm{q}{\ws{0}{2}{0}{\Pi}},
$$
but this is equivalent to solve
$$
\dive\left( \frac{\mF}{\ra} \right) = \frac{q}{\ra }. 
$$
If $q \in Y_0$ then $q/ \ra \in \ws{0}{2}{\alpha }{\Pi}$ and by Proposition \ref{surje}
there exists $\bfv\in \xalpha$ such that 
$$
\dive \bfv = \frac{q}{\ra} \text{ and } \snrm{\bfv}{\bws{1}{2}{\alpha}{\Pi}} \leq k(C_1,\alpha) \snrm{\frac{q}{\ra}}{ \ws{0}{2}{\alpha }{\Pi}}.
$$
Set $\mF:=\ra \bfv$ thanks to the isomorphism between $\bws{1}{2}{\alpha}{\Pi}$ and $\bws{1}{2}{0}{\Pi}$ there exists
a constant s.t. 
$$
\snrm{\mF}{\bws{1}{2}{0}{\Pi}} \leq k_1 \snrm{\bfv}{\bws{1}{2}{\alpha}{\Pi}}  \leq k_1 k(C_1,\alpha) \snrm{\frac{q}{\ra}}{ \ws{0}{2}{\alpha }{\Pi}} = k' \snrm{q}{ \ws{0}{2}{0 }{\Pi}}.
$$
\end{proof}

\begin{proof}[Proof of Theorem \ref{half_space_stk_thm}]
Thanks to the equivalence between well-posedness and conditions $(i)$ and $(ii)$
one concludes the existence and uniqueness of a pair $(\mF,\mG)$ solving problem 
\eqref{rescaled}. Moreover one has the {\em a priori} estimates~:
$$
\begin{aligned}
 \nrm{\mF}{\bws{1}{2}{0}{\Pi}} + \nrm{\mG}{\ws{0}{2}{0}{\Pi}} \leq k' & \left( 
\nrm{\ra  \Delta (\cR(\bfw)+\cS(\bfw))}{\bws{-1}{2}{0}{\Pi}} 
+ \nrm{\ra h}{\ws{-\ud}{2}{0}{N}} \right)
\end{aligned}
$$
they are obtained similarly to those of Theorem 2.34 p. 100 in \cite{ErGu.04.book}.
The isomorphism between weighted spaces mentioned above and the equivalence of 
problems \eqref{rescaled} and \eqref{half_space_stk_eq_dirichlet_homo} gives existence and 
uniqueness of $(\ti{\bfw},\theta)$ solving problem \eqref{half_space_stk_eq_dirichlet_homo} and
{\em a priori} estimates 
$$
\begin{aligned}
 \nrm{\ti{\bfw}}{\bws{1}{2}{\alpha}{\Pi}} + \nrm{\theta}{\ws{0}{2}{\alpha}{\Pi}} &\leq  k'' \left\{ 
\nrm{\Delta (\cR(\bfw)+\cS(\bfw))}{\bws{-1}{2}{\alpha}{\Pi}} + \nrm{h}{\ws{-\ud}{2}{\alpha}{N}} \right\}, \\
& \leq  k''' \left\{ \nrm{ \cR(\bfw)+\cS(\bfw)}{\bws{1}{2}{\alpha}{\Pi}} + \nrm{h}{\ws{-\ud}{2}{\alpha}{N}} \right\} .
\end{aligned}
$$
This gives existence and uniqueness of $(\bfw,\theta)$ and due to the continuity 
of the lifts $\cR(\bfw)$ and $\cS(\bfw)$  with respect to  the data, one easily proves that
$$
\begin{aligned}
 \nrm{\bfw}{\bws{1}{2}{\alpha}{\Pi}} + \nrm{\theta}{\ws{0}{2}{\alpha}{\Pi}} \leq k''' & \left\{ 
\nrm{{\bf f}}{\bws{\ud}{2}{\alpha}{D\cup N}} 
+ \nrm{h}{\ws{-\ud}{2}{\alpha}{N}} \right\} 
\end{aligned}
$$
which ends the proof
\end{proof}

\begin{proof}[Proof of Theorem \ref{MiStokesUnbounded}]
We use Theorem \ref{half_space_stk_thm} to prove results for $(\bfw_i,\theta_i)$ for $i \in \{ \bfbeta,\bfupsi,\bfchi,\bfvarkappa \}$.
The data for these problems tends exponentially fast to zero: in terms of weights, the Dirichlet (resp Neumann) data is thus compatible with
$\bws{\ud}{2}{\alpha}{D}$ (resp. $\bws{-\ud}{2}{\alpha}{N}$) for any real $\alpha$. 
Setting $\bf{f}= \oobeta \lambda -\bfbeta $ (resp.  $\bf{f}= \ooupsi \lambda -\bfupsi $ and $\bf{f}= \oochi \lambda -\bfchi $ )  in \eqref{half_space_stk_eq}
is equivalent to the first (resp. second and third) problem in \eqref{vertical_micro_corr}. 
Applying Theorem \ref{half_space_stk_thm} gives then the claims.
\end{proof}

\section{Periodic boundary layers: proofs of Propositions \ref{cell.pbm.beta.prop} and \ref{prop.upsi} and of Corollary \ref{prop.press.infty}}
\label{periodic.bl}
\begin{proof}[of Proposition \ref{cell.pbm.beta.prop}]
We start by lifting the non-homogeneous Dirichlet boundary condition~:
we set $\cR(\bfbeta):= y_2 \phi(y_2) \eu \chiu{Z_+}$ and $\tbeta:=\bfbeta + \cR(\bfbeta)$, this  is still a divergence
free vector. Now, there exists a unique solution s.t.
$$
\left|\nabla \tbeta \right|_{L^2(Z)^4} \leq \left|\nabla \cR(\bfbeta)\right|_{L^2(Z)^4}.
$$
Indeed, in the space of divergence free $\bmdu$ functions vanishing 
on $P$, the gradient norm is a norm (via Wirtinger estimates), one
 proves existence and uniqueness of $\tbeta$ in $\bmdu$ by the Lax-Milgram Theorem.

We apply Lemma 3.4 and Proposition 3.5 of \cite{JaMiPise.96}
in order to recover the $L^2_\loc(Z)$ pressure solving:
$$
-\Delta \ti{\bfbeta} + \nabla \pi = -\Delta \cR(\bfbeta),
$$
and this gives existence and uniqueness of $(\ti{\bfbeta},\pi) \in \bm{D}^{1,2}_0 ( Z) \times L^2_\loc (Z)$. On the interface located above (resp. below) 
the obstacle $\js$ we apply 
the Fourier decomposition in modes as in Theorem 3 p. 10 \cite{JaMiNeTechRep}.
One obtains the exponential convergence towards the zero  modes  of $(\bfbeta,\pi)$
in an explicit way. To derive relationships between  constant values at infinity,
one has
\renewcommand{\labelenumi}{(\roman{enumi})}
\begin{enumerate}
\item by the divergence free condition that $\ov{\beta}_2 (\nu)= \ov{\beta}_2(\gamma) =0$ for all 
$\nu \geq \ydp$ and $\gamma \leq 0$.
\item integrating the first equation of \eqref{beta.cell} in every
transverse section $\{ y_2=\delta \}$ which does not cross the obstacle $\js$ gives
$$
\frac{{\rm d}^2}{{\rm d }{y_2^2}} \left( \int_{\{ y_2=\delta\} } \beta_1 (y_1,y_2) \,dy_1 \right)= \int_{\{ y_2=\delta \} } - \dd{^2\beta_1}{y_1^2} + \dd{\pi}{y_1}\; dy_1 = 0,
$$
by $y_1$-periodicity. This implies that $\obeta(\delta)$ is an affine 
function of  $\delta$. As the gradient rapidly goes to zero, the linear part is zero, 
we conclude that only the constant remains~: thus $\obeta_1(\delta)=\obeta_1(+\infty)$ for $\delta>y_{2,P}$, and  $\obeta_1(\nu)=\obeta_1(0)$ for $\nu<0$.
\item Set $\bG:= y_2 \chiu{Z^+}\eu$ and $F:=0$, they  satisfy:
\begin{equation}\label{sol.fund}
\left\{
\begin{aligned}
- &\Delta  \bG + \nabla F= - \delta_\Sigma\eu  & \text{ in } Z,\\ 
& \dive \bG= 0 & \text{ in } Z,
\end{aligned}
\right.
\end{equation}
we test the first equation in \eqref{beta.cell} by $G$ and the
first equation in \eqref{sol.fund} by $\bfbeta$, then we integrate
on $Z_{\nu,\gamma}$~:
$$
\begin{aligned}
& (\Delta  \bG - \nabla F,\bfbeta)-(\Delta  \bfbeta - \nabla \pi ,\bG) = (\delta_\Sigma \eu,\bfbeta)= \int_0^1 \beta_1(y_1,0) dy_1 = \ov{\beta}_1(0) \\
& = (\sigma_{\bG,F} \cdot \bfn, \bfbeta)_{\partial Z_{\nu,\gamma}} - (\sigma_{\bfbeta,\pi} \cdot \bfn, \bG)_{\partial Z_{\nu,\gamma}} 
 = (\sigma_{\bG,F} \cdot \bfn, \bfbeta)_P + \ov{\beta}_1(\nu)- (\sigma_{\bfbeta,\pi} \cdot \bfn, \bG)_P  \\
& = - ( \ddn{} (y_2 \eu), y_2 \eu )_P +\ov{\beta}_1(\nu) + ( \sigma_{\bfbeta,\pi} \cdot \bfn , \bfbeta)_P  \\
& = - ( \ddn{} (y_2 \eu), y_2 \eu )_P +\ov{\beta}_1(\nu)+ \snrm{\nabla \bfbeta}{L^2(Z)},
\end{aligned}
$$
where we neglected exponentially small terms on $\{ y_2=\nu \}$ and $\{ y_2=\gamma \}$. 
Now we explicit the physical meaning
of the constant 
 ${\cal Q}:=( \ddn{} (y_2 \eu), y_2 \eu )_P$ 
$$
{\cal Q} + ( \ddn{} (y_2 \eu), y_2 \eu )_{\{ y_2=\nu\} \cup \{y_2=\gamma\} }
=  (\Delta (y_2 \eu) , y_2 \eu )_{Z_{\nu,\gamma}} + |\nabla (y_2 \eu) |^2_{L^2(Z_{\nu,\gamma})^4},
$$
which  in turn gives~:
$$
{\cal Q} + \nu- \gamma = \nu - \gamma - | \js |,
$$
where $|\js|$  is the volume of the obstacle $\js$. The quantity 
$\cal Q$ represents the volume of fluid missing due to the
presence of the obstacle $\js$ above the limit interface $\Sigma$.
If we were to consider a straight channel without a collateral
artery but a roughness below the fictitious interface,  $\cal Q$ would
be a positive number equal to the volume of fluid present
below $\Sigma$.

Computations above are formal and can be rigorously derived by regularizing the obstacle $\js$ and then working on regular functions in order to obtain results stated above. None of the final quantities depending  on second order derivatives,  passing to the limit wrt to the regularization parameter, 
  extends results above to Lipschitz obstacles. 
\end{enumerate}
\end{proof}

\begin{proof}[of Proposition \ref{prop.upsi}]
The existence and uniqueness part follows exactly the same lines as in Proposition \ref{cell.pbm.beta.prop}, the exponential convergence is also proved the same way. We detail only  relationships between  horizontal averages.
\begin{itemize}
\item For $\ov{\Upsilon}_2$ one uses as for $\obeta_2$ the divergence free condition together with the boundary condition at infinity to obtain that $\ov{\Upsilon}_2(y_2)=0$ for all $y_2$ in $\RR \setminus ]0,\ydp[$.
\item Testing the first equation in \eqref{eq.upsi} by $\bfupsi$ and integrating on $Z_{\gamma,\nu}$ one gets when passing to the limit wrt $\nu \to \infty , \gamma \to - \infty $ that:
$$
 \ov{\Upsilon}_1(0)= \nrm{ \nabla \bfupsi } {L^2(Z)}^2.
$$
Testing the same equation again but with $\bG = y_2 \eu$ and integrating on $Z_{\nu,\gamma}$ gives:
$$
\ov{\Upsilon}_1(\nu)= \ov{\Upsilon}_1(\gamma) + (\sigma_{\bfupsi,\varpi}\cdot \bfn, y_2 \eu )_P ,\quad \forall \nu > \ydp ,\quad \forall \gamma <0.
$$
We compute:
$$
\begin{aligned}
& (\dive \sigma_{\bfupsi,\varpi} , \bfbeta )_{Z_{\nu,\gamma}}- (\dive \sigma_{\bfbeta,\pi} , \bfupsi )_{Z_{\nu,\gamma}} = - \obeta_1(0) \\
=&  (\sigma_{\bfupsi,\varpi} \cdot \bfn , -y_2 \eu )_P + (\sigma_{\bfupsi,\varpi} \cdot \bfn ,\bfbeta  )_{\{ y_2=\nu \} \cup   \{ y_2=\gamma \}} -(\sigma_{\bfbeta,\pi} \cdot \bfn , \bfupsi )_{\{ y_2=\nu \} \cup   \{ y_2=\gamma \}} 
\end{aligned}
$$
which gives after passing to the limit wrt  $\nu$ and  $\gamma$ 
$$
\obeta_1(0)= (\sigma_{\bfupsi,\varpi} \cdot \bfn , y_2 \eu )_P.
$$
\item Testing the first equation in  \eqref{eq.upsi} against $y_2 \chiu{Z^-}$ and integrating on $Z_{\nu,\gamma}$ one obtains easily 
that
$$
\oupsi_1 (y_2)= \oupsi_1(0),\quad \forall y_2 <0.
$$
Putting together equalities obtained above one concludes the proof.
\end{itemize}
\end{proof}

\begin{proof}[of Corollary \ref{prop.press.infty}]
Setting again $\tbeta:=\bfbeta+y_2 \eu \chiu{Z_+}$ and writing
$$
\begin{aligned}
&(\Delta \bfchi -\nabla \eta ,\tbeta)_{Z_{\delta,\nu}} + (\Delta \bfbeta-\nabla \pi ,\bfchi)_{Z_{\delta,\nu}}=0 \\
& = (\sigma_{\bfchi,\eta} \cdot \bfn , \tbeta)_{\{y_2=\delta\}\cup \{ y_2=\nu\}} - (\sigma_{\bfbeta,\pi} \cdot \bfn , \bfchi )_{\{y_2=\delta\}\cup \{ y_2=\nu\}} .
\end{aligned}
$$
When passing to the limits $\delta\to \infty$ and $\nu\to - \infty$ in the last expression, one obtains that:
$$
  - [\overline{ \eta \beta_2 } ]^{+\infty}_{-\infty}  +  [\overline{ \pi \chi_2 } ]^{+\infty}_{-\infty} = 0,
$$
now because $\beta_2 \to 0$ and $\chi_2 \to -1$, one gets the desired result at infinity.
As the pressure $\pi$ is harmonic in $Z$, the average 
$\ov{\pi}(\delta)$ is zero in $\RR_-\cup]y_{2,P},+\infty[$. The same proof holds for 
$\varpi$.
\end{proof}

\bibliographystyle{plain}
\bibliography{eqred}

\end{document}

%% file: stents_pbm.tex
\begin{picture}(0,0)%
\includegraphics{stents_pbm.ps_tex}%
\end{picture}%
\setlength{\unitlength}{4972sp}%
\begingroup\makeatletter\ifx\SetFigFont\undefined%
\gdef\SetFigFont#1#2#3#4#5{%
  \reset@font\fontsize{#1}{#2pt}%
  \fontfamily{#3}\fontseries{#4}\fontshape{#5}%
  \selectfont}%
\fi\endgroup%
\begin{picture}(4357,1530)(500,-1378)
\end{picture}%

%% file: domains.tex
\begin{picture}(0,0)%
\includegraphics{domains.ps_tex}%
\end{picture}%
\setlength{\unitlength}{3108sp}%
\begingroup\makeatletter\ifx\SetFigFont\undefined%
\gdef\SetFigFont#1#2#3#4#5{%
  \reset@font\fontsize{#1}{#2pt}%
  \fontfamily{#3}\fontseries{#4}\fontshape{#5}%
  \selectfont}%
\fi\endgroup%
\begin{picture}(8712,2866)(166,-2330)
\put(181,-241){\makebox(0,0)[lb]{\smash{{\SetFigFont{9}{10.8}{\rmdefault}{\mddefault}{\updefault}{\color[rgb]{0,0,0}$\Gin$}%
}}}}
\put(1576,389){\makebox(0,0)[lb]{\smash{{\SetFigFont{9}{10.8}{\rmdefault}{\mddefault}{\updefault}{\color[rgb]{0,0,0}$\Gun$}%
}}}}
\put(496,-1546){\makebox(0,0)[lb]{\smash{{\SetFigFont{9}{10.8}{\rmdefault}{\mddefault}{\updefault}{\color[rgb]{0,0,0}$\Gde$}%
}}}}
\put(1576,-691){\makebox(0,0)[lb]{\smash{{\SetFigFont{9}{10.8}{\rmdefault}{\mddefault}{\updefault}{\color[rgb]{0,0,0}$\Geps$}%
}}}}
\put(1621,-2266){\makebox(0,0)[lb]{\smash{{\SetFigFont{9}{10.8}{\rmdefault}{\mddefault}{\updefault}{\color[rgb]{0,0,0}$\Goutd$}%
}}}}
\put(1576,-1096){\makebox(0,0)[lb]{\smash{{\SetFigFont{9}{10.8}{\rmdefault}{\mddefault}{\updefault}{\color[rgb]{0,0,0}$\Gz$}%
}}}}
\put(1576,-196){\makebox(0,0)[lb]{\smash{{\SetFigFont{9}{10.8}{\rmdefault}{\mddefault}{\updefault}{\color[rgb]{0,0,0}$\Oeu$}%
}}}}
\put(1576,-1636){\makebox(0,0)[lb]{\smash{{\SetFigFont{9}{10.8}{\rmdefault}{\mddefault}{\updefault}{\color[rgb]{0,0,0}$\Ode$}%
}}}}
\put(2971,-1546){\makebox(0,0)[lb]{\smash{{\SetFigFont{9}{10.8}{\rmdefault}{\mddefault}{\updefault}{\color[rgb]{0,0,0}$\Gde$}%
}}}}
\put(2926,-241){\makebox(0,0)[lb]{\smash{{\SetFigFont{9}{10.8}{\rmdefault}{\mddefault}{\updefault}{\color[rgb]{0,0,0}$\Goutu$}%
}}}}
\put(3781,-241){\makebox(0,0)[lb]{\smash{{\SetFigFont{9}{10.8}{\rmdefault}{\mddefault}{\updefault}{\color[rgb]{0,0,0}$\Gin$}%
}}}}
\put(5176,389){\makebox(0,0)[lb]{\smash{{\SetFigFont{9}{10.8}{\rmdefault}{\mddefault}{\updefault}{\color[rgb]{0,0,0}$\Gun$}%
}}}}
\put(4096,-1546){\makebox(0,0)[lb]{\smash{{\SetFigFont{9}{10.8}{\rmdefault}{\mddefault}{\updefault}{\color[rgb]{0,0,0}$\Gde$}%
}}}}
\put(5221,-2266){\makebox(0,0)[lb]{\smash{{\SetFigFont{9}{10.8}{\rmdefault}{\mddefault}{\updefault}{\color[rgb]{0,0,0}$\Goutd$}%
}}}}
\put(5176,-1096){\makebox(0,0)[lb]{\smash{{\SetFigFont{9}{10.8}{\rmdefault}{\mddefault}{\updefault}{\color[rgb]{0,0,0}$\Gz$}%
}}}}
\put(5176,-196){\makebox(0,0)[lb]{\smash{{\SetFigFont{9}{10.8}{\rmdefault}{\mddefault}{\updefault}{\color[rgb]{0,0,0}$\Ou$}%
}}}}
\put(5176,-1636){\makebox(0,0)[lb]{\smash{{\SetFigFont{9}{10.8}{\rmdefault}{\mddefault}{\updefault}{\color[rgb]{0,0,0}$\Ode$}%
}}}}
\put(6571,-1546){\makebox(0,0)[lb]{\smash{{\SetFigFont{9}{10.8}{\rmdefault}{\mddefault}{\updefault}{\color[rgb]{0,0,0}$\Gde$}%
}}}}
\put(6526,-241){\makebox(0,0)[lb]{\smash{{\SetFigFont{9}{10.8}{\rmdefault}{\mddefault}{\updefault}{\color[rgb]{0,0,0}$\Goutu$}%
}}}}
\put(7606,-1816){\makebox(0,0)[lb]{\smash{{\SetFigFont{9}{10.8}{\rmdefault}{\mddefault}{\updefault}{\color[rgb]{0,0,0}$Z_-$}%
}}}}
\put(8146,-511){\makebox(0,0)[b]{\smash{{\SetFigFont{9}{10.8}{\rmdefault}{\mddefault}{\updefault}{\color[rgb]{0,0,0}$\js$}%
}}}}
\put(7741,254){\makebox(0,0)[lb]{\smash{{\SetFigFont{9}{10.8}{\rmdefault}{\mddefault}{\updefault}{\color[rgb]{0,0,0}$Z_+$}%
}}}}
\put(7696,-1276){\makebox(0,0)[lb]{\smash{{\SetFigFont{9}{10.8}{\rmdefault}{\mddefault}{\updefault}{\color[rgb]{0,0,0}$\Sigma$}%
}}}}
\put(8011,-826){\makebox(0,0)[lb]{\smash{{\SetFigFont{9}{10.8}{\rmdefault}{\mddefault}{\updefault}{\color[rgb]{0,0,0}$P$}%
}}}}
\put(8056,-106){\makebox(0,0)[b]{\smash{{\SetFigFont{9}{10.8}{\rmdefault}{\mddefault}{\updefault}{\color[rgb]{0,0,0}$\jm$}%
}}}}
\put(8371, 29){\makebox(0,0)[lb]{\smash{{\SetFigFont{9}{10.8}{\rmdefault}{\mddefault}{\updefault}{\color[rgb]{0,0,0}$\gamma_M$}%
}}}}
\end{picture}%

%% file: velo_horiz.tex
\begin{picture}(0,0)%
\includegraphics{velo_horiz.ps_tex}%
\end{picture}%
\setlength{\unitlength}{1973sp}%
\begingroup\makeatletter\ifx\SetFigFont\undefined%
\gdef\SetFigFont#1#2#3#4#5{%
  \reset@font\fontsize{#1}{#2pt}%
  \fontfamily{#3}\fontseries{#4}\fontshape{#5}%
  \selectfont}%
\fi\endgroup%
\begin{picture}(5780,3553)(1286,-3987)
\put(6375,-789){\makebox(0,0)[rb]{\smash{{\SetFigFont{5}{6.0}{\familydefault}{\mddefault}{\updefault}$\uepsu(x_1,\epsilon)$}}}}
\put(6375,-1070){\makebox(0,0)[rb]{\smash{{\SetFigFont{5}{6.0}{\familydefault}{\mddefault}{\updefault}$\ovuepsu(x_1,\epsilon)$}}}}
\put(6375,-1351){\makebox(0,0)[rb]{\smash{{\SetFigFont{5}{6.0}{\familydefault}{\mddefault}{\updefault}$\uepsu(x_1,0)$}}}}
\put(6375,-1632){\makebox(0,0)[rb]{\smash{{\SetFigFont{5}{6.0}{\familydefault}{\mddefault}{\updefault}$\ovuepsu(x_1,0^-)$}}}}
\put(1763,-3623){\makebox(0,0)[rb]{\smash{{\SetFigFont{5}{6.0}{\familydefault}{\mddefault}{\updefault}-0.005}}}}
\put(1763,-3115){\makebox(0,0)[rb]{\smash{{\SetFigFont{5}{6.0}{\familydefault}{\mddefault}{\updefault} 0}}}}
\put(1763,-2607){\makebox(0,0)[rb]{\smash{{\SetFigFont{5}{6.0}{\familydefault}{\mddefault}{\updefault} 0.005}}}}
\put(1763,-2099){\makebox(0,0)[rb]{\smash{{\SetFigFont{5}{6.0}{\familydefault}{\mddefault}{\updefault} 0.01}}}}
\put(1763,-1590){\makebox(0,0)[rb]{\smash{{\SetFigFont{5}{6.0}{\familydefault}{\mddefault}{\updefault} 0.015}}}}
\put(1763,-1082){\makebox(0,0)[rb]{\smash{{\SetFigFont{5}{6.0}{\familydefault}{\mddefault}{\updefault} 0.02}}}}
\put(1763,-574){\makebox(0,0)[rb]{\smash{{\SetFigFont{5}{6.0}{\familydefault}{\mddefault}{\updefault} 0.025}}}}
\put(1838,-3748){\makebox(0,0)[b]{\smash{{\SetFigFont{5}{6.0}{\familydefault}{\mddefault}{\updefault} 0}}}}
\put(2865,-3748){\makebox(0,0)[b]{\smash{{\SetFigFont{5}{6.0}{\familydefault}{\mddefault}{\updefault} 0.2}}}}
\put(3893,-3748){\makebox(0,0)[b]{\smash{{\SetFigFont{5}{6.0}{\familydefault}{\mddefault}{\updefault} 0.4}}}}
\put(4920,-3748){\makebox(0,0)[b]{\smash{{\SetFigFont{5}{6.0}{\familydefault}{\mddefault}{\updefault} 0.6}}}}
\put(5948,-3748){\makebox(0,0)[b]{\smash{{\SetFigFont{5}{6.0}{\familydefault}{\mddefault}{\updefault} 0.8}}}}
\put(6975,-3748){\makebox(0,0)[b]{\smash{{\SetFigFont{5}{6.0}{\familydefault}{\mddefault}{\updefault} 1}}}}
\put(4406,-3935){\makebox(0,0)[b]{\smash{{\SetFigFont{5}{6.0}{\familydefault}{\mddefault}{\updefault}$x_1$}}}}
\end{picture}%

%% file: velo_vert.tex
\begin{picture}(0,0)%
\includegraphics{velo_vert.ps_tex}%
\end{picture}%
\setlength{\unitlength}{1973sp}%
\begingroup\makeatletter\ifx\SetFigFont\undefined%
\gdef\SetFigFont#1#2#3#4#5{%
  \reset@font\fontsize{#1}{#2pt}%
  \fontfamily{#3}\fontseries{#4}\fontshape{#5}%
  \selectfont}%
\fi\endgroup%
\begin{picture}(5780,3553)(1286,-3987)
\put(6375,-789){\makebox(0,0)[rb]{\smash{{\SetFigFont{5}{6.0}{\familydefault}{\mddefault}{\updefault}$\uepsd(x_1,0)$}}}}
\put(6375,-1070){\makebox(0,0)[rb]{\smash{{\SetFigFont{5}{6.0}{\familydefault}{\mddefault}{\updefault}$\ovuepsd(x_1,0)$}}}}
\put(1763,-3623){\makebox(0,0)[rb]{\smash{{\SetFigFont{5}{6.0}{\familydefault}{\mddefault}{\updefault}-0.025}}}}
\put(1763,-3115){\makebox(0,0)[rb]{\smash{{\SetFigFont{5}{6.0}{\familydefault}{\mddefault}{\updefault}-0.02}}}}
\put(1763,-2607){\makebox(0,0)[rb]{\smash{{\SetFigFont{5}{6.0}{\familydefault}{\mddefault}{\updefault}-0.015}}}}
\put(1763,-2098){\makebox(0,0)[rb]{\smash{{\SetFigFont{5}{6.0}{\familydefault}{\mddefault}{\updefault}-0.01}}}}
\put(1763,-1590){\makebox(0,0)[rb]{\smash{{\SetFigFont{5}{6.0}{\familydefault}{\mddefault}{\updefault}-0.005}}}}
\put(1763,-1082){\makebox(0,0)[rb]{\smash{{\SetFigFont{5}{6.0}{\familydefault}{\mddefault}{\updefault} 0}}}}
\put(1763,-574){\makebox(0,0)[rb]{\smash{{\SetFigFont{5}{6.0}{\familydefault}{\mddefault}{\updefault} 0.005}}}}
\put(1838,-3748){\makebox(0,0)[b]{\smash{{\SetFigFont{5}{6.0}{\familydefault}{\mddefault}{\updefault} 0}}}}
\put(2865,-3748){\makebox(0,0)[b]{\smash{{\SetFigFont{5}{6.0}{\familydefault}{\mddefault}{\updefault} 0.2}}}}
\put(3893,-3748){\makebox(0,0)[b]{\smash{{\SetFigFont{5}{6.0}{\familydefault}{\mddefault}{\updefault} 0.4}}}}
\put(4920,-3748){\makebox(0,0)[b]{\smash{{\SetFigFont{5}{6.0}{\familydefault}{\mddefault}{\updefault} 0.6}}}}
\put(5948,-3748){\makebox(0,0)[b]{\smash{{\SetFigFont{5}{6.0}{\familydefault}{\mddefault}{\updefault} 0.8}}}}
\put(6975,-3748){\makebox(0,0)[b]{\smash{{\SetFigFont{5}{6.0}{\familydefault}{\mddefault}{\updefault} 1}}}}
\put(4406,-3935){\makebox(0,0)[b]{\smash{{\SetFigFont{5}{6.0}{\familydefault}{\mddefault}{\updefault}$x_1$}}}}
\end{picture}%

%% file: q2_colateral.tex
\begin{picture}(0,0)%
\includegraphics{q2_colateral.ps_tex}%
\end{picture}%
\setlength{\unitlength}{1579sp}%
\begingroup\makeatletter\ifx\SetFigFont\undefined%
\gdef\SetFigFont#1#2#3#4#5{%
  \reset@font\fontsize{#1}{#2pt}%
  \fontfamily{#3}\fontseries{#4}\fontshape{#5}%
  \selectfont}%
\fi\endgroup%
\begin{picture}(5927,3553)(1216,-3987)
\put(6375,-2845){\makebox(0,0)[rb]{\smash{{\SetFigFont{5}{6.0}{\familydefault}{\mddefault}{\updefault}direct}}}}
\put(6375,-3126){\makebox(0,0)[rb]{\smash{{\SetFigFont{5}{6.0}{\familydefault}{\mddefault}{\updefault}explicit}}}}
\put(6375,-3407){\makebox(0,0)[rb]{\smash{{\SetFigFont{5}{6.0}{\familydefault}{\mddefault}{\updefault}theoretical}}}}
\put(1888,-3623){\makebox(0,0)[rb]{\smash{{\SetFigFont{5}{6.0}{\familydefault}{\mddefault}{\updefault} 0.002}}}}
\put(1888,-3284){\makebox(0,0)[rb]{\smash{{\SetFigFont{5}{6.0}{\familydefault}{\mddefault}{\updefault} 0.004}}}}
\put(1888,-2945){\makebox(0,0)[rb]{\smash{{\SetFigFont{5}{6.0}{\familydefault}{\mddefault}{\updefault} 0.006}}}}
\put(1888,-2607){\makebox(0,0)[rb]{\smash{{\SetFigFont{5}{6.0}{\familydefault}{\mddefault}{\updefault} 0.008}}}}
\put(1888,-2268){\makebox(0,0)[rb]{\smash{{\SetFigFont{5}{6.0}{\familydefault}{\mddefault}{\updefault} 0.01}}}}
\put(1888,-1929){\makebox(0,0)[rb]{\smash{{\SetFigFont{5}{6.0}{\familydefault}{\mddefault}{\updefault} 0.012}}}}
\put(1888,-1590){\makebox(0,0)[rb]{\smash{{\SetFigFont{5}{6.0}{\familydefault}{\mddefault}{\updefault} 0.014}}}}
\put(1888,-1252){\makebox(0,0)[rb]{\smash{{\SetFigFont{5}{6.0}{\familydefault}{\mddefault}{\updefault} 0.016}}}}
\put(1888,-913){\makebox(0,0)[rb]{\smash{{\SetFigFont{5}{6.0}{\familydefault}{\mddefault}{\updefault} 0.018}}}}
\put(1888,-574){\makebox(0,0)[rb]{\smash{{\SetFigFont{5}{6.0}{\familydefault}{\mddefault}{\updefault} 0.02}}}}
\put(1963,-3748){\makebox(0,0)[b]{\smash{{\SetFigFont{5}{6.0}{\familydefault}{\mddefault}{\updefault} 0.05}}}}
\put(2520,-3748){\makebox(0,0)[b]{\smash{{\SetFigFont{5}{6.0}{\familydefault}{\mddefault}{\updefault} 0.1}}}}
\put(3077,-3748){\makebox(0,0)[b]{\smash{{\SetFigFont{5}{6.0}{\familydefault}{\mddefault}{\updefault} 0.15}}}}
\put(3634,-3748){\makebox(0,0)[b]{\smash{{\SetFigFont{5}{6.0}{\familydefault}{\mddefault}{\updefault} 0.2}}}}
\put(4191,-3748){\makebox(0,0)[b]{\smash{{\SetFigFont{5}{6.0}{\familydefault}{\mddefault}{\updefault} 0.25}}}}
\put(4747,-3748){\makebox(0,0)[b]{\smash{{\SetFigFont{5}{6.0}{\familydefault}{\mddefault}{\updefault} 0.3}}}}
\put(5304,-3748){\makebox(0,0)[b]{\smash{{\SetFigFont{5}{6.0}{\familydefault}{\mddefault}{\updefault} 0.35}}}}
\put(5861,-3748){\makebox(0,0)[b]{\smash{{\SetFigFont{5}{6.0}{\familydefault}{\mddefault}{\updefault} 0.4}}}}
\put(6418,-3748){\makebox(0,0)[b]{\smash{{\SetFigFont{5}{6.0}{\familydefault}{\mddefault}{\updefault} 0.45}}}}
\put(6975,-3748){\makebox(0,0)[b]{\smash{{\SetFigFont{5}{6.0}{\familydefault}{\mddefault}{\updefault} 0.5}}}}
\put(1331,-2037){\rotatebox{-270.0}{\makebox(0,0)[b]{\smash{{\SetFigFont{5}{6.0}{\familydefault}{\mddefault}{\updefault}$Q_{\Gz}$}}}}}
\put(4469,-3935){\makebox(0,0)[b]{\smash{{\SetFigFont{5}{6.0}{\familydefault}{\mddefault}{\updefault}$\epsilon$}}}}
\end{picture}%

%% file: error_l2.tex
\begin{picture}(0,0)%
\includegraphics{error_l2.ps_tex}%
\end{picture}%
\setlength{\unitlength}{1973sp}%
\begingroup\makeatletter\ifx\SetFigFont\undefined%
\gdef\SetFigFont#1#2#3#4#5{%
  \reset@font\fontsize{#1}{#2pt}%
  \fontfamily{#3}\fontseries{#4}\fontshape{#5}%
  \selectfont}%
\fi\endgroup%
\begin{picture}(5850,3553)(1216,-3987)
\put(3238,-789){\makebox(0,0)[rb]{\smash{{\SetFigFont{5}{6.0}{\familydefault}{\mddefault}{\updefault}$\ueps- \bfuz$}}}}
\put(3238,-1070){\makebox(0,0)[rb]{\smash{{\SetFigFont{5}{6.0}{\familydefault}{\mddefault}{\updefault}$\epsilon^{0.9}$}}}}
\put(3238,-1351){\makebox(0,0)[rb]{\smash{{\SetFigFont{5}{6.0}{\familydefault}{\mddefault}{\updefault}$\ueps- \ovueps$}}}}
\put(3238,-1632){\makebox(0,0)[rb]{\smash{{\SetFigFont{5}{6.0}{\familydefault}{\mddefault}{\updefault}$\epsilon^{1.5}$}}}}
\put(1888,-3623){\makebox(0,0)[rb]{\smash{{\SetFigFont{5}{6.0}{\familydefault}{\mddefault}{\updefault} 0.001}}}}
\put(1888,-2098){\makebox(0,0)[rb]{\smash{{\SetFigFont{5}{6.0}{\familydefault}{\mddefault}{\updefault} 0.01}}}}
\put(1888,-574){\makebox(0,0)[rb]{\smash{{\SetFigFont{5}{6.0}{\familydefault}{\mddefault}{\updefault} 0.1}}}}
\put(1963,-3748){\makebox(0,0)[b]{\smash{{\SetFigFont{5}{6.0}{\familydefault}{\mddefault}{\updefault} 0.01}}}}
\put(4469,-3748){\makebox(0,0)[b]{\smash{{\SetFigFont{5}{6.0}{\familydefault}{\mddefault}{\updefault} 0.1}}}}
\put(6975,-3748){\makebox(0,0)[b]{\smash{{\SetFigFont{5}{6.0}{\familydefault}{\mddefault}{\updefault} 1}}}}
\put(1331,-2037){\rotatebox{-270.0}{\makebox(0,0)[b]{\smash{{\SetFigFont{5}{6.0}{\familydefault}{\mddefault}{\updefault}$\nrm{\cdot}{\bL^2(\Oe)}$}}}}}
\put(4469,-3935){\makebox(0,0)[b]{\smash{{\SetFigFont{5}{6.0}{\familydefault}{\mddefault}{\updefault}$\epsilon$}}}}
\end{picture}%

%% file: error_hm1.tex
\begin{picture}(0,0)%
\includegraphics{error_hm1.ps_tex}%
\end{picture}%
\setlength{\unitlength}{1973sp}%
\begingroup\makeatletter\ifx\SetFigFont\undefined%
\gdef\SetFigFont#1#2#3#4#5{%
  \reset@font\fontsize{#1}{#2pt}%
  \fontfamily{#3}\fontseries{#4}\fontshape{#5}%
  \selectfont}%
\fi\endgroup%
\begin{picture}(5850,3553)(1216,-3987)
\put(3388,-789){\makebox(0,0)[rb]{\smash{{\SetFigFont{5}{6.0}{\familydefault}{\mddefault}{\updefault}$\peps- \pz$}}}}
\put(3388,-1070){\makebox(0,0)[rb]{\smash{{\SetFigFont{5}{6.0}{\familydefault}{\mddefault}{\updefault}$\epsilon^{1.15}$}}}}
\put(3388,-1351){\makebox(0,0)[rb]{\smash{{\SetFigFont{5}{6.0}{\familydefault}{\mddefault}{\updefault}$\peps- \ovpeps$}}}}
\put(3388,-1632){\makebox(0,0)[rb]{\smash{{\SetFigFont{5}{6.0}{\familydefault}{\mddefault}{\updefault}$\epsilon^{1.425}$}}}}
\put(1888,-3623){\makebox(0,0)[rb]{\smash{{\SetFigFont{5}{6.0}{\familydefault}{\mddefault}{\updefault} 0.001}}}}
\put(1888,-2098){\makebox(0,0)[rb]{\smash{{\SetFigFont{5}{6.0}{\familydefault}{\mddefault}{\updefault} 0.01}}}}
\put(1888,-574){\makebox(0,0)[rb]{\smash{{\SetFigFont{5}{6.0}{\familydefault}{\mddefault}{\updefault} 0.1}}}}
\put(1963,-3748){\makebox(0,0)[b]{\smash{{\SetFigFont{5}{6.0}{\familydefault}{\mddefault}{\updefault} 0.01}}}}
\put(4469,-3748){\makebox(0,0)[b]{\smash{{\SetFigFont{5}{6.0}{\familydefault}{\mddefault}{\updefault} 0.1}}}}
\put(6975,-3748){\makebox(0,0)[b]{\smash{{\SetFigFont{5}{6.0}{\familydefault}{\mddefault}{\updefault} 1}}}}
\put(1331,-2037){\rotatebox{-270.0}{\makebox(0,0)[b]{\smash{{\SetFigFont{5}{6.0}{\familydefault}{\mddefault}{\updefault}$\nrm{\cdot}{H^{-1}(\Ou'\cup\leps \cup \Ode)}$}}}}}
\put(4469,-3935){\makebox(0,0)[b]{\smash{{\SetFigFont{5}{6.0}{\familydefault}{\mddefault}{\updefault}$\epsilon$}}}}
\end{picture}%

%% file: hmesh.tex
\begin{picture}(0,0)%
\includegraphics{hmesh.ps_tex}%
\end{picture}%
\setlength{\unitlength}{1973sp}%
\begingroup\makeatletter\ifx\SetFigFont\undefined%
\gdef\SetFigFont#1#2#3#4#5{%
  \reset@font\fontsize{#1}{#2pt}%
  \fontfamily{#3}\fontseries{#4}\fontshape{#5}%
  \selectfont}%
\fi\endgroup%
\begin{picture}(5850,3553)(1216,-3987)
\put(6375,-3126){\makebox(0,0)[rb]{\smash{{\SetFigFont{5}{6.0}{\familydefault}{\mddefault}{\updefault}$h_{\max}$}}}}
\put(6375,-3407){\makebox(0,0)[rb]{\smash{{\SetFigFont{5}{6.0}{\familydefault}{\mddefault}{\updefault}$h_{\min}$}}}}
\put(1888,-3623){\makebox(0,0)[rb]{\smash{{\SetFigFont{5}{6.0}{\familydefault}{\mddefault}{\updefault} 0.001}}}}
\put(1888,-2098){\makebox(0,0)[rb]{\smash{{\SetFigFont{5}{6.0}{\familydefault}{\mddefault}{\updefault} 0.01}}}}
\put(1888,-574){\makebox(0,0)[rb]{\smash{{\SetFigFont{5}{6.0}{\familydefault}{\mddefault}{\updefault} 0.1}}}}
\put(1963,-3748){\makebox(0,0)[b]{\smash{{\SetFigFont{5}{6.0}{\familydefault}{\mddefault}{\updefault} 0.01}}}}
\put(4469,-3748){\makebox(0,0)[b]{\smash{{\SetFigFont{5}{6.0}{\familydefault}{\mddefault}{\updefault} 0.1}}}}
\put(6975,-3748){\makebox(0,0)[b]{\smash{{\SetFigFont{5}{6.0}{\familydefault}{\mddefault}{\updefault} 1}}}}
\put(1331,-2037){\rotatebox{-270.0}{\makebox(0,0)[b]{\smash{{\SetFigFont{5}{6.0}{\familydefault}{\mddefault}{\updefault}$h$}}}}}
\put(4469,-3935){\makebox(0,0)[b]{\smash{{\SetFigFont{5}{6.0}{\familydefault}{\mddefault}{\updefault}$\epsilon$}}}}
\end{picture}%